\input amstex
\documentstyle{amsppt}
\input xy
\xyoption{all}
\CompileMatrices
\xyReloadDrivers

\def\t#1 {\text{\rm #1}}
\def\b#1 {\text{\rm\bf #1}}


\topmatter
\title Controlled K-theory I: Basic theory
  \endtitle
\author Frank Quinn\endauthor
\date February 2004\enddate
\address Math, Virginia Tech, Blacksburg VA 24061-0123 USA\endaddress
\email quinn\@math.vt.edu\endemail
\thanks Partially supported by the National Science Foundation\endthanks
\abstract This paper provides a full controlled version of algebraic $K$-theory. There is a rich array of assembly maps; the
controlled assembly isomorphism theorem identifying the controlled group with homology; and the stability theorem
describing the behavior of the inverse limit as the control parameter goes to 0. There is also a careful treatment of spectral cosheaf homology and related tools, including an ``iterated homology identity''  giving a spectrum-level version of the Leray-Serre spectral sequence.  
\endabstract

\endtopmatter

\head 1. Introduction\endhead
Controlled algebra was originally
developed to encode obstructions in controlled topology \cite{Quinn 1, 2}.  At that time  attempts to develop controlled {\it algebraic\/}
$K$-theory were unsuccessful, and the topological obstructions were instead formulated using  pseudoisotopy.  For many
applications this is satisfactory, cf.
\cite{Quinn 2, 5, 7}, \cite{Hughes} but there are serious drawbacks: genuine $K$-theory is needed for a full controlled version
of surgery; and computations using the assembly map structure are greatly enhanced by access to the full range of
$K$-theory tools. For example the remarkable work of  Farrell and Jones and their coworkers \cite{Farrell-Jones 1, 2}
\cite{Berklove-Farrell-Juan-Pineda-Pearson}, \cite{Davis-L\"uck}, \cite{Roushon}, \cite{Aravinda-Farrell-Roushon} is done in the
context of pseudoisotopies in order to use the controlled pseudoisotopy theory developed in
\cite{Quinn 2}. From this they deduce information about Whitehead and lower $K$ groups of integral group rings. The $K$-theory
developed here should substitute for pseudoisotopy in many of these proofs, so similar conclusions should  hold for
higher
$K$-theory and arbitrary coefficient rings. Recently \cite{Bartels-Reich} have also made a start in this direction using continuously controlled $K$-theory. The continuous theory describes the inverse limit of finite-$\epsilon$ sets, but it is missing the ``stability'' theorem that describes how they approach this limit. The stability theorem is a crucial ingredient of many applications, c.f\. \cite{Quinn 2, 8}.  

This paper is intended as a foundation for applications. However the ones known to the author require a fair amount of additional development so will be undertaken in subsequent papers.  The first of these is \cite{Quinn 8}, describing a refinement of the Farrell-Jones ``fibered isomorphism conjecture''.    

Section 2 contains statements of results, with sufficient explanation to make them usable
without extensive reference to the body of the paper. The primary result in Part I is the Controlled assembly isomorphism theorem 2.2.1. A version of stability is included as Theorem 2.3.1.  

The body of the paper begins in Section 3 with a
version of uncontrolled $K$-theory using chain complexes. Most of the section is occupied by the proof that this agrees with
other definitions. Section 4 adds control to the definitions of \S3 and develops the elementary properties of the controlled
$K$-space. In Section 5 deeper properties (those requiring spacial localization techniques) are developed. In particular the
homology axioms needed in the following section are verified. Section 6 gives a development and axiomatic characterization of spectral cosheaf
homology and assembly maps. The Controlled Assembly Isomorphism theorem follows from this since controlled $K$-theory satisfies
the axioms.

\head 2. Results \endhead
The central construction is a simplicial set of $\epsilon$-controlled chain complexes over a metric space, briefly
described in 2.1. The first main result is the
controlled assembly isomorphism theorem of 2.2. This describes inverse limits of controlled $K$-spaces as homology with
spectral cosheaf coefficients.  The stability theorem of 2.3 asserts that the inverse system approaches the limit in a
nice way. 
 An application to the structure of topological
group actions is given in Section 2.4. 

\subhead 2.1 Controlled $K$-theory spaces\endsubhead

Our central object is a simplicial set $K^{lf}_1(X;p,R,\epsilon)$  defined given
\roster\item a locally compact metric space $X$, 
\item a map $p\:E\to X$,
\item a ring $R$, and
\item a real number $\epsilon>0$.
\endroster
This space is defined in \S4 using controlled geometric chain complexes. Very roughly these
are 2-complexes decorated with   $R$ data, mapping into $E$. The vertices of the 2-complex correspond to
$R$-modules, the edges are homomorphisms between the modules, and 2-cells encode relations among the morphisms. These are
``finitely generated'' in the sense that the map of the 2-complexes into $X$ are proper (thus locally finite), and
``controlled'' in the sense that each cell has size less than $\epsilon$ in
$X$. Straightforward properties  needed in formulating the results are given in 2.1.1--2.1.3. Inverse limits are
described in 2.1.4--2.1.5.

\subsubhead 2.1.1 Naturality\endsubsubhead A morphism between maps $(E@>p>> X) \to (F@>q>> Y)$ is a commutative
diagram
$$\CD E@>\hat f>>F\\
@VV{p}V@VV{q}V\\
X@>f>>Y\endCD$$
A morphism induces a map 
$$K^{lf}_1(X;p,R,\epsilon_1)@>>>K^{lf}_1(Y;q,R,\epsilon_2)$$
 defined simply by applying $\hat f$ to decorated 2-complexes mapping into $E$. These maps compose nicely, so
define a functor on an appropriate catgory.

\subsubhead 2.1.2 Restriction\endsubsubhead 
If $U\subset X$ is open then there is a restriction map
$$K^{lf}_1(X;p,R,\epsilon)\to K^{lf}_1(U;p|U,R,\epsilon).$$
In terms of decorated 2-complexes in $E$ this is given by restricting to the largest subcomplex mapping into $p^{-1}(U)$.
Getting this to make sense near the edge of $U$ requires a wrinkle in the definition that will be appreciated by experts,
see \S4.1.7.

\subsubhead 2.1.3 Other properties\endsubsubhead 
We mention several other basic aspects 
that will be explained in more detail when they are needed. 

The first such aspect is that the path
components of the space are essentially the $\epsilon-K_1$ obstructions arising from controlled h-cobordisms in
\cite{Quinn 1, 2}. This gives the connection to topology and also motivates the notation for the space. 

Another basic topic is the use of control {\it functions\/} $\epsilon\:X\to(0,\infty)$ rather than just constant
$\epsilon>0$. This only makes a difference when $X$ is noncompact, but is almost always required in
topological applications with noncompact $X$. The reasons we stick with constants here are: (1)
details are simpler; (2) the function case follows  from a relative form of the compact constant case by patching;
and (3) the constant-$\epsilon$ conclusions for noncompact $X$ are stronger than the function-$\epsilon$
version, and the difference should be useful in geometric applications.

 Finally, the  role of the space $E$ is to specify ``local
fundamental groups'' over
$X$. If
$E\to E'$ is a map of spaces over $X$ that (roughly) is an isomorphism of fundamental groups of inverse images of points
in $X$, then the induced map of $K_1$ spaces is (roughly) an equivalence. This is a consequence of the use of
2-complexes in the definition. This invariance property is extremely useful in  applications but plays a minor role here.

\subsubhead 2.1.4 The  inverse limit\endsubsubhead 
Define 
$$K^{lf}_1(X;p,R) = \t holim _{\epsilon\to0} K^{lf}_1(X;p,R,\epsilon).$$
The homotopy inverse limit is the space of paths to $\infty$ in the system. More precisely we define an inverse system
indexed on positive integers, with the $n^{th}$ space the inverval $[n,\infty)$ and maps the inclusions. Then the
homotopy inverse limit is the space of cofinal maps of inverse systems from this into the $K^{lf}_1(X;p,R,*)$ system. 

The on-the-nose (i.e\. not homotopy) inverse limit is the intersection 
$$\cap_{\epsilon}K^{lf}_1(X;p,R,\epsilon)=K^{lf}_1(X;p,R,0)$$ and is uninteresting. The homotopy version loosens this
up by allowing elements at different scales to be ``homotopic'' rather than equal. Behavior of the inverse system as it
approaches the limit is desribed in \S2.4.

\subhead 2.2 Controlled assembly \endsubhead
This theorem uses the ``spectral cosheaf homology''  described in \S8 of \cite{Quinn 2} and in more detail in \S6. The basic
idea is that a continuous spectrum-valued functor of spaces can be applied to point inverses of a map $E\to X$ to get a
``spectral cosheaf'' over $X$. Homology with coefficients in such a thing can be defined and satisfies appropriate versions
of the usual properties of homology. We think of homology as  accessible, at least  more so than $K$-theory, so a homological
description is good news.  A  version for compactly-supported homology is given in~2.2.4.

\proclaim{2.2.1 Controlled assembly isomorphism theorem} Suppose $X$ is a locally compact metric space, $p\:E\to X$ a map.
Then
\roster\item {\rm (metric independence)} up to homotopy $K^{lf}_1(X,\;p,R)$ is independent of the metric on $X$.
\item {\rm (spectrum structure)} $K^{lf}_1(X,\;p,R)$ has a natural $\Omega$-spectrum structure. Denote the first
deloop of this structure by
$\Bbb K^{lf}(X;p,R)$. 
\item {\rm (coefficients)} If\/ $E$ is a connected space with the homotopy type of a CW complex then choice of basepoint
$*\in E$ gives an equivalence of spectra with the (nonconnective) algebraic $K$-theory spectrum
$$\Bbb K(\t pt ;E\to\t pt , R)@>\simeq >> {\Bbb K}(R[\pi_1(E,*)]).$$
\item {\rm (assembly isomorphism)} If\/ $p$ is
a stratified system of fibrations and $X$ is an ANR then the assembly map
associated to the spectrum structure is an equivalence of spectra,
$$\Bbb H^{lf}(X;\Bbb K(p,R))@>\simeq >> \Bbb K^{lf}(X;p,R).$$
\endroster
\endproclaim
\subsubhead 2.2.2  Notes\endsubsubhead
\roster
\item Stratified systems of
fibrations are defined in \cite{Quinn 2} and described  in 6.2.1. The homotopy stratified group actions developed in \cite{Quinn 5} naturally give rise to these in profusion.
\item In (4) $\Bbb K(p,R)$ denotes the spectral cosheaf over $X$ obtained by applying $\Bbb K$ fiberwise to $p\:E\to X$.
This process is described in detail in \S6.7--6.8. Note that (3) identifies the fibers in this cosheaf  with standard 
$K$-theory of the fundamental group of the point inverses of $p$. Other assembly maps can be obtained by using naturality in
$p$.
\item Statement (2) can be thought of as saying this functor is
``self-computing'': not only is it homology, but it supplies the spectrum needed to define the homology and the assembly map
that makes the comparison. 
\item The significance of
(1) is illustrated by comparing
$R^n$ with its standard metric  to the metric induced by a homeomorphism with the open unit ball. The first metric is
complete so the typical
$K_1$ gadget is infinite, has
$\delta$-dense image, and satisfies the relations everywhere. The second metric is not complete, and in fact the
complement of any open neighborhood of the metric frontier (see 2.1.2) is compact. This means $K_1$ gadgets are finite and 
there is a big ``ragged edge'' near the metric frontier where structure relations fail. Nonetheless in the limit the ``ragged
edge'' shrinks away, and the map from the complete metric to the incomplete one does not lose information.
\endroster

\subsubhead 2.2.3 Compactly supported homology\endsubsubhead 
Ordinary homology can be obtained as the direct limit of homology of finite complexes mapping into the space. Since
finite complexes are compact, ``locally finite'' and ordinary homology coincide, so we define
$$K_1(X;p,R) = \t holim _{h\:Y\to X} K^{lf}_1(Y;h^*p,R),$$
where the limit is over the system of finite complex $Y\to X$, and $h^*p\:h^*E\to Y$ denotes the pullback (fiber product)
of
$p$ along
$h$. Morphisms in the inverse system are maps of finite complexes commuting with the maps to $X$. 

This version is defined for all $X, p$ (no local compactness requirement) and is functorial with respect to all
morphisms. The locally finite  theorem 2.2.1 gives, by  direct limits and a little fiddling to align filtrations:

\proclaim{Theorem} Suppose $X$ is a 
metric space and $p\:E\to X$ is a stratified system of fibrations. Then the limit
$K_1(X;p,R)$ has a natural spectrum structure, and  the assembly map 
$$\Bbb H(X;\Bbb K(p,R))@>\simeq >> \Bbb K(X;p,R).$$
 is an equivalence of spectra.
\endproclaim

\subhead 2.3 Stability\endsubhead
The $K$-space  $K^{lf}_1(X;p,R)$ is defined to be the homotopy inverse limit of spaces $K^{lf}_1(X;p,R,\epsilon)$, as
$\epsilon\to 0$. ``Stability'' makes
 properties of the limit available at finite-$\epsilon$ scales by showing the inverse system converges nicely. The version
given here is sufficient for most applications. To get a non-compact statement we 
extend the definition of $K^{lf}_1(X;p,R,\epsilon)$ in the evident way to permit $\epsilon$ to be a continuous function
$X\to(0,\infty)$. See 2.1.3 for a brief discussion of control functions.  

\proclaim{2.3.1 Proposition} Suppose $X$ is a finite dimensional locally compact metric ANR, $p\:E\to X$
is a stratified system of fibrations,  and $\epsilon\:X\to(0,\infty)$ is given. Then there is $\delta\:X\to (0,\infty)$
so that the relax-control map $\delta\to\epsilon$ homotopically factors through the limit in the sense that there is a map $S$
so that the diagram
$$\xymatrix{ K^{lf}_1(X;p,R,\delta)\ar[r]\ar[dr]^S& K^{lf}_1(X;p,R,\epsilon)\\
& K^{lf}_1(X;p,R)\ar[u]}$$
homotopy commutes.\endproclaim

This follows by plugging the lemmas of \S5 into the proof of stability for pseudoisotopy \cite{Quinn 2}. This proof can  be simplified using ideas in the proof of stability for surgery groups in \cite{Pedersen-Yamasaki}. Sharper versions are also possible, but we will wait for guidance from a significant application before trying to formulate this. 

\subhead 2.4 Group actions and induction\endsubhead
Here we restrict to homotopy stratified maps obtained from group actions. These are used in studying the local geometric
structure of group actions, and in induction approaches to describing $K$-theory.

We suppose $G$ is a (discrete) group and $X$ is a ANR with a homotopically stratified action of $G$. This class of actions
is developed in \cite{Quinn 5}, and includes  CW complexes
with a cellular action of $G$.  For homotopy purposes there is no loss in supposing $X$ is a $G$-CW complex. The
starting point is:
\proclaim{2.4.1 Lemma} Suppose $E\to X$ is $G$-equivariant and a fibration, $E$ is free and $X$ is an ANR with a homotopy
stratified action. Then the quotient map $E/G\to X/G$ is homotopy stratified over the orbit-type filtration of
$X/G$.\endproclaim It follows that spectral cosheaf homology is defined using this as reference map, and Theorem 2.2.1 applies.

 A consequence of  induction theory for finite groups
\cite{Dress 1} is that all finite groups satisfy hyperelementary induction. 

\subhead 2.5 An application to topology\endsubhead
This illustrates that  the formal framework already
has interesting consequences. 

 Suppose $f\:X\to Y$ is an equivariant map of homotopically stratified $G$ actions on ANRs, and
$p\:E\to Y$ is a fibration and equivariant with respect to a free action on $E$. Form the pullback
$$\CD f^*E@>>> E\\
@VV{f*^p}V@VV{p}V\\
X@>F>>Y\endCD$$
and consider the induced map on $K$:
$$\Bbb K(X/G;f^*p,R)@>>>\Bbb K(Y/G;p,E).$$
We say $f$ is a $K$-{\it theory equivalence\/} of $G$-spaces if this induced map is an equivalence of spectra. In some cases
we can detect $K$-equivalence using subgroups:
\proclaim{2.5.1 Proposition} Suppose $f$ as above is a $K$-theory equivalence of $H$-spaces for all hyperelementary subgroups 
$H\subset G$, and suppose all isotropy groups of the action on $Y$ satisfy hyperelementary induction. Then $f$ is a
$K$-theory equivalence of $G$-spaces.\endproclaim
This is proved with several applications of the iterated-homology identity 6.11 to relate it to a change-of-coefficients map in
which the $H$ versions appear in the coefficients. 

In particular $Y\times \Cal E_{G,\text{h-elem}}\to Y$ is an $H$-equivariant homotopy equivalence for any hyperelementary
$H\subset G$, so is a $K$-equivalence of $H$-spaces. It follows that it induces a $K$ equivalence as $G$ spaces, provided
the isotropy subgroups satisfy induction. 

Now suppose that $X$ is a compact  manifold with a  homotopy stratified action of a finite group $G$.
Let $X^*\subset X$
denote the set of points with nontrivial isotropy subgroups.  There is an obstruction to the existence of an equivariant
mapping cylinder neighborhood for
$X^*$ in
$X$, lying in a reduced controlled
$K_0$ of the quotient $X^*/G$, \cite{Quinn 2}. This group is $\pi_0$ (slightly reduced) of the homology of the quotient with
$\Bbb K$ coefficients. The Proposition above applies because the isotropy subgroups are finite, so we conclude the whole
spectrum is detected on quotients by finite hyperelementary groups. We only need  $\pi_0$, but whole-spectrum conclusions are easier than specializations to $\pi_0$. 
\proclaim{2.5.2 Corollary} Suppose $G$ acts on $X$ and $X^*\subset X$ is the non-free set as above. If\/
$X^*$ has $H$-equivariant mapping cylinder neighborhoods for all hyperelementary $H\subset G$ then it has a
$G$-equivariant mapping cylinder neighborhood.\endproclaim

\head 3. Chain complexes and $K_1(R)$\endhead
The space $K_1(R)$ is defined as the simplicial set whose simplices are diagrams of chain complexes
satisfying various conditions. Notations, suspensions, and mapping cones are described in 3.1, 
Triangular  simplices are defined in 3.4. $K_1(R)$ is finally defined in 3.5. One of the main technical tools, the
cancellation of inverses, makes its first appearence in 3.6. In 3.7 we give the first application: a
deformation retraction of the
$K_1$ space to a subcomplex equivalent to the Volodin description of $K$-theory. 

As mentioned above there is considerable similarity between this material and  \cite{Igusa-Klein 1} and
\cite{Igusa 1} on ``higher Reidemeister torsion''. Myriad details and their use of categories prevent us from making
use of this work.

\subhead 3.1 Chain complexes\endsubhead
This section mainly fixes notation and sign conventions. Our sign conventions differ from some
others, and this is explained in~3.1.5. 

\subsubhead 3.1.1 Notation\endsubsubhead
To simplify the formulae we use a uniform notation for the boundary homomorphisms and chain maps in the
structure. For example a chain complex will have the form $(C,c)$ where $C$ are the chain groups and $c\:C\to C$ denotes
the boundary homomorphism. Lower indices will be used to indicate different complexes, and upper indices degrees in the
complex. For instance  $C_i^j$ indicates the
degree-$j$ group of the chain complex
$C_i$. 
\subsubhead 3.1.2 Simplices of complexes\endsubsubhead
Elaborating on 3.1.1, an $n$-{\it simplex\/} of complexes
consists of 
 chain complexes $(C_i,c_i)$ for $i=0,\dots,n$ and for $i<j$, chain maps $c_{i,j}\:C_i\to C_j$. 

The face $\partial_i$ of a simplex is defined by omitting the $i^{th}$ vertex.
An $(n-1)$-simplex is supposed to be indexed on the set of integers $[0,n-1]$, so more precisely $\partial_i$ is 
 obtained from the order-preserving injection $[0,n-1]\to [0,n]$ onto  the complement of $i$. 

More generally if $J\subset [0,n]$ we define $\partial_J$ by composing with the order-preserving injection $[0,n-|J|]\to
[0,n]$ whose image is the complement of $J$. If we denote the injection by $\alpha$ then 
$$(\partial_J(C,c))_i = C_{\alpha(i)}.$$
This degree of precision in the notation is rarely needed. 

\subsubhead 3.1.3 Suspensions\endsubsubhead  
The suspension of a $Z$-graded module $C$ is denoted by $SC$, and is the same module with grading
shifted up by one. Explicitly, $(SC)^r=C^{r-1}$. Thus for instance the identity map $SC\to C$ has degree
$-1$. 

The suspension of
 an
$n$-simplex of complexes $(C,c)$ has graded module the suspension $SC$ and structure maps with signs according to the
parity of their degree:
$$\align sc_i&=(-1)c_j\\
sc_{i,j}&=(+1)c_{i,j}\endalign$$
 See 3.1.5 for comments on this choice of sign
conventions.
\subsubhead 3.1.4 Mapping cones\endsubsubhead  
Mapping cones of chain isomorphisms will be used to add trivial summands in $K$-theory.

Suppose $f\:C\to \bar C$ is a chain map. Then  
the {\it cone\/} $\bar C\oplus_{f}SC$ is the complex with graded modules $\bar C\oplus SC$, and 
boundary homomorphisms $\left(\smallmatrix \bar c&f\\ 0& -c\endsmallmatrix\right)$.

\proclaim{Lemma} If\/ $f$ is an isomorphism then $\bar C\oplus_{f}SC$ is contractible with canonical
contraction 
$\left(\smallmatrix 0&0\\ f^{-1}& 0\endsmallmatrix\right)$\endproclaim
Contractibility of the cone requires only that $f$ be a chain homotopy equivalence, though in this case the
contraction depends on the chain homotopies as well as the map itself. However the isomorphism case is the one needed
here.

\subsubhead 3.1.5 Sign conventions\endsubsubhead  
Sign adjustments are required in
suspensions and cones to get the various identities to work out. There are two standard choices. The one used here
has signs in the suspension structure: odd degree structure morphisms (boundary
homomorphisms and homotopies) pick up a $-1$ upon suspension. Suspensions of chain maps are unchanged because they have
degree 0. And if $f\:C\to D$ is a degree-0 chain map, then the same homomorphism is a degree-$(-1)$ chain
map $SC\to D$. The different sign requirements for the different degrees (see 3.1.3) comes from the signs
in the suspension.

The other convention does not change signs in 
suspensions, but then signs are required in morphisms.  If $f\:C\to D$ is a degree-0 chain map
then to get  a degree-$(-1)$ chain map $SC\to D$ one uses $(-1)^nf^n$ on the degree-$n$ summand of
$SC$. This is the convention used by Ranicki, for instance. The notation in this convention is more
elaborate because it explicitly uses the grading in $C$ rather than the structure degree. 

Topologically these conventions correspond to different orders in a product. If the suspension of
$\Delta^n$ is taken to be $I\times \Delta^n$  then the usual product sign conventions give the signed
suspension used here. If the suspension is taken to be $\Delta^n\times I$ then suspensions don't get
signs, but morphisms do.

\subhead 3.2 Triangular morphisms and simplices\endsubhead
Diagonal, triangular and increasing   morphisms of based modules are defined in 3.2.1--3. In 3.3.4 diagonal morphisms
on infinitely generated modules are shown to be well-behaved if the partial orders are ``locally bounded''.
\subsubhead 3.2.1 Based modules\endsubsubhead 
As usual, a {\it basis\/} for a module is a set of generators as a free module over the ring. In fact a basis renders the
algebraic structure of the module itself irrelevant,  but for the time being we retain the language as a comfortable
context. 

 If $C\supset D$
is a based submodule (i.e\. generated by a subset of the basis) then $C\perp D$ is the {\it perpendicular\/} submodule,
generated by the rest of the basis. There is a canonical sum decomposition $C=D\oplus(C\perp D)$, and isomorphism $C\perp
D\simeq C/D$. In controlled situations quotient modules behave badly, and perpendicular submodules will be used
instead. Finally if $f\:C\to \bar C$ is a homomorphism of based modules and $C\supset D$, $\bar C\supset \bar D$ are
based submodules, then $f$ induces homomorphisms $D\to \bar D$ and $C\perp D\to \bar C\perp \bar D$. If we represent $f$
as a matrix with respect to the bases then these induced homomorphisms are obtained by taking appropriate submatrices.
Algebraically they are described as compositions with the canonical inclusions and projections associated to the sum
decomopositions. 

\subsubhead 3.2.2 Diagonal homomorphisms\endsubsubhead 
Suppose $C$, $D$ are free  modules with specified bases denoted by $\t base (C)$, $\t base (D)$ respectively. The {\it
support\/},
$\t supp (f)$,  of a homomorphism
$f\:C\to D$ is the subset of the basis of $C$ on which $f$ is nontrivial. The homomorphism is {\it diagonal\/} if is
induced by an injection  $\t base (f)\:\t supp (f)\to \t base (D)$. 

More generally, if $U$ is a subset of $R$ then $f$
is $U$-diagonal if there is a diagonal map $\t base (f)$ so that for each basis element
$x\in\t base (C)$,
$f(x)$ is a multiple of
$\t base (f)(x)$ by an element of the subset $U$.  In these terms ``diagonal'' is the same as $\{1\}$-diagonal. Other
useful cases are $U=\{\pm 1\}$, and $U$ the units of $R$. In topology one often has a group ring $R=Z[\pi]$, and
$U=\pm\pi$.   Completely arbitrary  coefficients correspond to
$R$-diagonal maps.

In the definition of ``diagonal'' we are extending the convention that the free module with empty basis is 0, not the
empty module. Namely, extend the basis function of $f$ to $\t base (C)\to \t base (D)\cup\{\emptyset\}$ by taking all
elements not in the support of $f$ to the empty set. Then we take as a convention that the only possible multiple of the
empty basis element is 0. These conventions simplify the language in several places. For example the composition of
diagonal maps is diagonal, and with this convention if the coefficients are units then the basis function of the
composition is the composition of basis functions.

Note that  an $R$-diagonal homomorphism is an
isomorphism if and only if the underlying basis function is a bijection of bases, and the coefficient on each basis
element is a unit in the ring. In this case the inverse is evidently also diagonal, with inverse base function.

\subsubhead 3.2.3 Triangular homorphisms\endsubsubhead 
Fix a subgroup $U$ of the units of $R$. Suppose $f\:C\to D$ is a homomorphism of free based $R$-modules with  partial
orders on the bases.
$f$ is $U$-{\it triangular\/} with respect to the partial orders if
 there is a  decomposition $f= h+u$ where
\roster\item $h$ is $U$-diagonal; \item $\t base (h)\:\t supp
(h)\to
\t base (D)$ is order-preserving; and
\item $u$ is {\it increasing with respect to\/} $h$ in the sense that for each basis element $x\in C$, $u(x)$ lies in the
submodule of
$D$ generated by basis elements strictly greater than $\t base (h)(x)$ in the partial order.
\endroster

With appropriate hypotheses (see 3.4.4) the decomposition is unique. In this case the diagonal part is denoted by $\t diag
(f)$. We are using the convention that $\t base (h)(x)=\emptyset$ if $x$ is not in the support of $h$. In this case there
are no greater basis elements and so $u(x)=0$. In other words the support of $u$ is contained in the support of $h$. 

Note that the increasing condition in (3) only uses the partial order in $\t base (D)$, and in particular does not depend
on the underlying base function of
$h$ being order-preserving.  Both conditions are needed to show that composition of triangular functions is triangular.

When the domain and range are the same there is a natural notion of ``increasing'': $u(x)$ lies in the submodule
generated by basis elements strictly larger than $x$. In the terminology above this is the same as ``increasing with
respect to the identity'', and is equivalent to requiring $1+u$ to be triangular.

We explain why we use the term ``increasing''  instead of the more common ``upper triangular'' or ``lower
triangular''.  The usual practice is to refine the partial order to a total order, and then write the morphism as a
matrix.  We think of elements of the modules as column vectors, with morphisms acting as matrices on the right.
``Diagonal'' entries then occur in blocks along the diagonal (and can usually be arranged to be exactly on the
diagonal). ``Increasing'' entries lie below this (block) diagonal. There are two problems. First, we make extensive use
of matrix notations, but make no attempt to relate the implicit order in the matrix to the structural partial orders.
Thus the matrices often have entries on both sides of the diagonal. Second it seems to invite
confusion to try to use ``lower'' as a synonym for ``increasing''.

\subsubhead 3.2.4 Locally bounded partial orders\endsubsubhead
This is a condition to ensure that triangular morphisms of infinitely generated modules are well-behaved. It is
automatically satisfied for finitely-generated modules, and we use it here to simplify definitions. In controlled
situations definitions the condition may not hold, and we will have to assume some of its consequences for other
reasons.

A partial order is {\it locally bounded\/} if there is no element that is the
starting point of arbitrarily long increasing chains.  We say a {\it module\/} is locally bounded if the partial order on
the basis is locally bounded. The following is standard:

\proclaim{Lemma} In the following homomorphisms are understood to be triangular homomorphisms of modules with locally
bounded partially ordered bases.
\roster\item the diagonal-plus-increasing decomposition  is
well-defined, and is denoted by\quad $\t diag (f)+\t inc (f)$;
\item  the diagonal part is functorial: $\t diag (fg)=\t diag (f)\t diag (g)$,
\item if $f-g$ is increasing with respect to the diagonal part of one of $f$, $g$, then 
$\t diag (f)=\t diag (g)$,
\item A diagonal-plus-increasing homomorphism  is an isomorphism if and only if 
the diagonal part is an isomorphism. This is equivalent to the basis function being bijective and the coefficients being
units;
\item the inverse of a triangular isomorphism is triangular.  
\item there is a canonical (independent of the partial order) factorization $f=(1+\alpha_n)\dots (1+\alpha_1)\t diag (f)$ with $\alpha_i$ increasing and
$\alpha_i\alpha_j=0$ if\/ $j\leq i$. 
\endroster
\endproclaim
Suppose $f\:C\to D$ is triangular.
To see (1), define a decreasing filtration on the basis of $D$ with $n^{th}$ subset the elements that are the smallest in
an increasing chain of length at least $n$. Let $D_n$ be the submodule generated by the $n^{th}$ subset. Then $D_0=D$, and
$\cap_nD_n=0$ by the locally bounded condition. The increasing part of any  triangular map must lie in $D_1$, so
$f\:C\to D_0/D_1$ is induced by an injection on part of the basis and takes the rest to 0. This determines $h$ on the part
of the basis that does not map into $D_1$, and thus also $u$ on this subset. Restrict to the complement, then we have
$f_1\:E_1\to D_1$  triangular. The same argument determines the decomposition on the part of the basis that does not
map into
$D_2$. Continuing by induction shows the entire decomposition is determined. 

For  statement (2) let $f=h+u$, $g=j+v$ denote the decompositions. Then $fg=hj+(hv+uj+uv)$. $hj$ is diagonal, and since
$h,j$ are order-preserving and $u,v$ are increasing, the other terms are increasing. This gives a
diagonal-plus-increasing decomposition of $fg$. 

Statement (3) should be clear.

For (4) it is easy to see that if the diagonal part is not an isomorphism then the whole
homomorphism cannot be either. Conversely suppose the diagonal part is an isomorphism. Denote the decomposition by
$h+u$, then the inverse is given by 
$$(h+u)^{-1} = \Sigma_{i=0}^{\infty}(-h^{-1}u)^ih^{-1}.$$
The infinite sum is well-defined because it is finite when applied to any element of $D$. This follows from the bounded
hypothesis by essentially the argument used to show (1). The form of the inverse gives (5). For the second statement
in (4) a diagonal morphism decomposes as a direct sum of a zero part, and a sum of endomorphisms of $R$ indexed by the
underlying basis function. To be an isomorphism the basis function must be a bijection, and each 1-dimensional summand
must be an isomorphism. The latter condition is that the coefficients must be units in the ring.

Finally for (6) we set $f(\t diag f)^{-1}= 1+u$, with $u$ increasing.  We factor $1+u$ using essentially the opposite
of the argument for (1).  Let $D_0\subset D_1\subset \cdots$ be a maximal chain of based submodules so that $D_0=0$ and $u(D_i)\subset D_{i-1}$ for $i>0$. Since $u$ is increasing and the order is bounded this sequence terminates with $D_n=D$ some $n$. Let $p_i$ denote the based projection $D\to D_i$. Inductively define $\alpha_i$ by $\alpha_0=0$ and 
$$\alpha_i=(u-\Sigma_{j=0}^{i-1}\alpha_j)p_i.$$
It follows that $\alpha_ip_i=\alpha_i$ (definition); $\alpha_ip_j=0$ for $j<i$ (induction); $p_{i-1}\alpha_i=\alpha_i$ (definition of $D_i$); and $u=\Sigma_{j=0}^n\alpha_j$ (since $D_n=D$). From these it follows that $\alpha_i\alpha_j=0$ if $j\leq i$. Thus 
$$(1+\alpha_n)\cdots (1+\alpha_2)(1+\alpha_1)=1+\Sigma_{j=0}^n\alpha_j=1+u$$.

\subhead 3.3 Contractions and cancellation \endsubhead 
We use chain contractions to identify parts of a complex that can safely be omitted. Lemma 3.3.2 shows this is the
inverse of the usual stabilization operation, but with more flexibility in how the trivial summands are related to the
chosen bases. An explanation of the topological significance is given in 3.3.3.
\subsubhead 3.3.1 Definition \endsubsubhead
Suppose $C$ is a based complex with a partial order on the
basis in each degree, and $\xi$ is a chain contraction for $C$. We say $\xi$ 
{\it cancels the complement\/} of  a based subcomplex  $\hat C\subset C$ if
 \roster\item in each degree and when they are comparable, basis elements of $C\perp \hat C$ preceed those of $\hat C$ in the partial order, and 
\item there is a based decomposition $C\perp \hat C=D \oplus \bar D$ so that  $\xi$ 
has the form
$$\pmatrix \hat \xi&u&v\\ 0&0 &\delta\\ 0&0&0\endpmatrix\:\hat C\oplus D\oplus \bar D\longrightarrow \hat C\oplus D\oplus \bar D $$
with $\delta$ a
$\pm1$-triangular isomorphism.
\endroster
This is a very rigid condition.  If the partial order and $\xi$ are fixed then (1) implies that
the based subcomplexes for which this works are linearly ordered by inclusion, and the decomposition in (2) is unique
when it exists. The following is also clear:
\proclaim{3.3.2 Lemma} If\/  $\xi$ cancels the complement of\/ $\hat C\subset C$ then it preserves $\hat C$ and the
restriction ($\hat \xi$ in the matrix above) is a chain contraction of\/ $\hat C$.  
If\/  $\hat B\subset \hat C\subset C$ are based subcomplexes and\/ $\xi$  cancels
the complements in\/ $C$ of both,
  then $\xi|\hat C$  cancels the complement of\/ $\hat B$ in~$\hat C$.
\endproclaim

\subsubhead 3.3.3 Structure of the boundary \endsubsubhead

If $C = \hat C\oplus D\oplus \bar D$ is a cancellation decomposition as above then the
boundary homomorphism $c\:C\to C$ has the form 
$$\pmatrix \hat c& x&y\\ 0&d&z\\0&\delta^{-1}&\bar d\endpmatrix.$$
The lower-left 0 entries reflect the hypothesis that $\hat C$ is a subcomplex. The  $\delta^{-1}$ term comes from the
hypothesis that $\xi$ is a chain contraction.

\subsubhead 3.3.4 Standardizing cancellations \endsubsubhead The following lemma shows that cancellations
can be standardized, so the free entries in the matrix above reflect  flexibility in embedding the structure into
the complex rather than flexibility in the structure. Versions of this are used to standardize dimensions in 3.8 and
in the local cancellation of inverses in \S5.
\proclaim{Lemma } Suppose $C= \hat C \oplus D\oplus \bar D$ is a cancellation decomposition of $C$ and use  the
notation above for the matrix form of the boundary homomorphism. Then
$f$ is a chain isomorphism to the same graded module with boundary
homomorphism
$b$ and contraction $\beta$, to $C$,  where
$$f=\pmatrix 1&0&x\delta \\0&1&d\delta\\0&0&1\endpmatrix \text{,\quad} b=\pmatrix
\hat c&0&0\\0&0&0\\0&\delta^{-1}&0\endpmatrix\text{,\quad and \quad}\beta=\pmatrix
\hat \xi&u&-\hat c u\\0&0&\delta\\0&0&0\endpmatrix$$
\endproclaim
\demo{Proof} The chain-map claim  is that $fb=cf$, the contraction claim is that $f\beta=\xi f$. Multiply out the
matrices given, then this can be seen to follow from $c^2=0$. Alternatively,  $b$, $\beta$ are obtained by conjugating by $f$ and $f^{-1}$ is obtained by multiplying the off-diagonal terms of $f$ by $-1$. $f$ is  $\pm1$ triangular with
respect to a suitable partial order because $\delta$ is.
\enddemo
There are several useful remarks about this:
\roster\item The complement of $\hat C$ is also a subcomplex in this structure, so in particular the projection to $\hat C$ is a chain map.
\item $f$ is determined by $\hat C$ and $\xi$, and in particular is independent of a choice of partial order. The submodule $D$ is the image of $\xi$ intersected with the complement of $\hat C$, so it's complement $\bar D$ is also determined. Let $p$ denote the based projection $C \to \bar D$, then $f= (1-p)+cp\xi$ is a canonical description. 
\item The contraction $\beta$ can be further improved. The degree-2 homomorphism 
$$\pmatrix 0&0&a\\ 0&0&0\\0&0&0\endpmatrix$$ gives a chain homotopy to the same form with $a=0$.
\endroster

\subhead 3.4 Definition of $K_1(R)$\endsubhead
$K_1(R)$ is a simplicial set with simplices defined using triangular simplices of  chain complexes with contractions.
The basic idea is that a $K_1$ simplex is a sum of a nondegenerate part, in which the chain maps are isomorphisms, and a
cancellable part in which they may not be. Genuine
$K$-theory activity takes place in the nondegenerate part, while the  cancellations organized by the contractions
enable stabilization of dimension. Most of the complexity of the theory comes from the cancellations: not from
$K$-theory itself, but from the structures needed to keep it from escaping during stabilization.

\subsubhead 3.4.1 Simplices in $K_1(R)$\endsubsubhead 
An $n$-simplex of $K_1(R)$ is a finitely generated based chain $n$-simplex $C$ together with  chain contractions
$\xi_i$ for $(C_i,c_i)$, such that there are partial orders on the bases in each degree
so that for each $i<j$,
\roster
\item $c_{i,j}$  is a $\pm1$ triangular chain map  $C_i\to C_j$;
\item the elements in the image of the basis function of $c_{i,j}$ come after the elements in the complement, in the
partial order on each $C_j^n$;
\item if $i<j<k$ then $c_{j,k}c_{i,j}$ has the same underlying basis function as $c_{i,k}$;
\item $\xi$ cancels the  complement of the image $C_j\perp c_{i,j}(C_i)\subset C_j$.
\endroster
Note that the partial orders  are only required to exist and are not part of the data.  Diagrams of chain maps
 in the simplex, such as the triangle in (3), are not required to commute. However 
condition (3) is a weak version of ``commute up to increasing chain homotopy''.  These conditions are delicately
balanced in being strong enough to capture
$K$-theory but flexible enough to work in the cancellation-of-inverses construction of the next section and  various
localizations of it.  The following shows the images are based subcomplexes, which is required for condition (4) to
make sense.
\proclaim{Lemma} The image of  $c_{i,j}$ is the based subcomplex spanned by the image of the underlying basis function,
and if\/ $i<j<k$ then $c_{j,k}c_{i,j}$ has the same image as $c_{i,k}$.\endproclaim
The first part  follows from the triangularity of $c_{i,j}$ and the order condition (2). The second part follows from
this and the basis-function hypothesis (3). 

 Faces of such simplices are defined in the evident
way: the $j^{th}$ face omits the complex
$C_j$.  Faces of $K_1$ simplices are again in $K_1$, so $K_1(R)$ is a simplicial complex. More precisely
this gives it the structure of a
$\Delta$-set (omit degeneracies from the definition of simplicial set,
\cite{Rourke-Sanderson 1}). One can either work with these as $\Delta$-sets or introduce degeneracies by introducing
duplicates of the
$C_i$. In any case we proceed without mention of degeneracies and leave the choice of details to the reader.

\subsubhead 3.4.2 The image partial order\endsubsubhead We expand on the partial order condition (2) in the
definition of simplices in $K_1(R)$. Suppose $(C,c,\xi)$ is an $n$-simplex. Let $C_{i,j}$ denote the image of
$c_{i,j}$ in $C_j$, and for notational convenience we set
$C_{j,j}=C_j$. Recall that
$C^n_{i+1,j}\perp C^n_{i,j}$ denotes the based complement of the smaller image in the larger, in degree $n$. Since  
in  each degree the complement of $C_{i,j}$ preceeds it in the partial order we get 
$$(C^n_{j,j}\perp C^n_{j-1,j})<\cdots < (C^n_{i,j}\perp C^n_{i-1,j})<\cdots <(C^n_{1,j}\perp C^n_{0,j})<C^n_{0,j}.$$
We will refer to this as the ``image partial order'' on $C_j^n$. Any partial order making $(C,c,\xi)$ a $K_1$ simplex
must refine this partial order.

\subsubhead 3.4.3 Morphisms\endsubsubhead 
Morphisms of $K_1$ simplices are dictated in a straightforward way by the structure. The key property is that
triangulation gives $K_1$ simplices, so morphisms can be used to construct simplicial homotopies. Composition of
$K_1$ morphisms is generally {\it not\/} a morphism, because the partial orders will usually not be compatible.

A {\it morphism\/} of $K_1$ simplices $f\:(A,\alpha)\to(C,\xi)$ consists of chain maps $f_{i,j}\:A_i\to C_j$ for
$i\leq j$ such that
 there are choices of partial orders  satisfying
\roster\item $(A,\alpha)$ and $(C,\xi)$ are $K_1$ simplices with respect to the partial orders;
\item in each degree each $f_{i,j}$ is a $\pm1$ triangular morphism whose image basis elements follow those in the
complement;
\item if $i\leq j\leq k$ then the morphisms $c_{j,k}f_{i,j}$, $f_{i,k}$, and $f_{j,k}a_{i,j}$ all have the same
underlying set function; and
\item $\xi_j$ cancels the complement of the image of $f_{i,j}$.
\endroster
\subsubhead 3.4.4 Triangulation\endsubsubhead 
Suppose $f\:A\to C$ is a morphism of $K_1$ simplices. We think of this as a diagram modeled on $\Delta^n\times I$, and
triangulate it into $(n+1)$-simplices corresponding to the canonical triangulation of the model. 
Specifically we get $n+1$ such simplices, the $k$th one of which has vertices $A_0,\dots A_k,C_k, \dots,C_n$, and chain
maps the $a_{i,j}$, $c_{i,j}$ and $f_{i,j}$ whose domain and range are both vertices. 
$\partial_0$ of the 0th simplex is $C$, $\partial_{n+1}$ of the $(n+1)$st is $A$, and the $k$th and $(k+1)$th have
exactly one face in common. 

 The following is
straightforward:
\proclaim{Lemma} Simplices in the triangulation of a $K_1$ morphism are $K_1$ simplices.
\endproclaim
We describe how this is used to construct homotopies. 
Suppose $X$ is a simplicial complex, and $A, C\:X\to K_1(R)$ are simplicial maps. Suppose further that there are
morphisms $f_{\sigma}\:A(\sigma)\to C(\sigma)$ given for every $\sigma\in X$ so that
$f_{\partial_i\sigma}=\partial_if_{\sigma}$. Then simplices in the triangulations of the $f_{\sigma}$ fit together to
define a simplicial homotopy from $A$ to $C$. ``Simplicial homotopy'' means simplicial map from the standard simplicial
structure on
$X\times I$.

\subhead 3.5 Cancellation of inverses\endsubhead 
Direct sum defines a monoid structure on $K_1$. In this section we show the suspension is a homotopy inverse for this
structure by constructing a nullhomotopy of the sum
$\t id \oplus S \:K_1(R)\to K_1(R)$. This nullhomotopy is described in detail because elaborations of it are the bases
for most of the results in the paper. The first  application is in 3.7 to compress complexes into two degrees.

The construction is canonical, so it is sufficient to describe the effect on an arbitrary simplex of $K_1$ and they will
fit together to define a homotopy of the entire space. The  idea is that the chain
contraction provides an isomorphism to the cone on the identity: if $(C,\xi)$ is an $n$-simplex then
$$\left(\smallmatrix  1&-\xi\\0&1\endsmallmatrix\right)\:C\oplus SC@>\simeq>> C\oplus_{1} SC.$$
The cone on the identity has an obvious global cancellation structure ($\bar D=C$, $ D= SC$) so 
$0\to C\oplus_{1} SC$ is a morphism of
$K_1$ simplices. Subdividing these two morphisms (see 3.4.4) gives a homotopy to 0. In detail:

\proclaim{ 3.5.1 Lemma} Suppose $(C,c,\xi)$ is an $n$-simplex of $K_1(R)$.
\roster\item $C\oplus SC$ with contraction $\left(\smallmatrix 
\xi&0\\0&-\xi\endsmallmatrix\right)$ is a $K_1$ simplex;
\item the cone $C\oplus_1 SC$ with contraction $\left(\smallmatrix 
0&0\\1&0\endsmallmatrix\right)$ is a $K_1$ simplex; 
\item $f_{i,j}=\left(\smallmatrix 
1&-\xi_j\\0&1\endsmallmatrix\right)\left(\smallmatrix 
c_{i,j}&0\\0&c_{i,j}\endsmallmatrix\right)$ for $i\leq j$ is a morphism between these simplices; and 
\item the inclusion  $0\to C\oplus_1SC$ is a $K_1$ morphism. 
\endroster\endproclaim
\demo{Proof} The main thing to check is that there are appropriate partial orders in which contractions cancel
complements. Choose
$K_1$-simplex partial orders on bases, and as in 3.4.2 denote the image of
$c_{i,j}$ by
$C_{i,j}$ . Partial order the bases in the modules $C\oplus SC$ using the given partial order in each summand and by
shuffling the image partial orders (3.4.2) in the two pieces, with $SC$ terms before $C$ terms. More specifically the
degree
$n$ module in $C_j\oplus SC_j$ is partially ordered as follows (recall $(SC_j)^n=C_j^{n-1}$):
$$\align(C^{n-1}_{j,j}\perp C^{n-1}_{j-1,j})<(C^n_{j,j}\perp C^n_{j-1,j})&<\cdots \\ < (C^{n-1}_{i,j}\perp
C^{n-1}_{i-1,j})&<
 (C^n_{i,j}\perp C^n_{i-1,j})<\cdots < C^{n-1}_{0,j}<C^n_{0,j}.\endalign$$

\noindent $i\to j$ {\it maps}\quad In both the sum and the cone the $i\to j$ maps are
$\left(\smallmatrix  c_{i,j}&0\\0&c_{i,j}\endsmallmatrix\right)$. These are $\pm1$ triangular injections because
$c_{i,j}$ is. The images are the based subcomplexes $C_{i,j}\oplus S C_{i,j}$ and satisfy the image order condition
because the $c_{i,j}$ do. The
$i\to j$ map in the morphism is the composition of $\left(\smallmatrix  c_{i,j}&0\\0&c_{i,j}\endsmallmatrix\right)$
and 
$\left(\smallmatrix  1&-\xi\\0&1\endsmallmatrix\right)$. The first factor is $\pm1$ triangular as before. The second is
a
$\pm1$ triangular isomorphism because $\xi$ preserves images, by Lemma 3.4.1, and $SC$ terms in the image partial order
preceed $C$ terms. Specifically the summand $(C^{n-1}_{i,j}\perp C^{n-1}_{i-1,j})$ is taken  to $C^n_{i,j}=\oplus_{k\leq
i}(C^n_{k,j}\perp C^n_{k-1,j})$, and all parts of this are greater in the partial order. Thus the non-diagonal part of
$\left(\smallmatrix  1&-\xi\\0&1\endsmallmatrix\right)$ is increasing.
 Note  that since $\xi_j$ preserves
$C_{i,j}$, the morphism   $\left(\smallmatrix  1&-\xi\\0&1\endsmallmatrix\right)\left(\smallmatrix 
c_{i,j}&0\\0&c_{i,j}\endsmallmatrix\right)$  has image a based subcomplex satisfying the image order condition, namely
$C_{i,j}\oplus SC_{i,j}$.

Finally all these chain maps have underlying basis functions the product of basis functions on the summands. Since the
original simplex has commutative basis-function diagrams, so do the sums.

\medskip\noindent{\it Cancellation conditions}\quad In the direct sum the contractions cancel complements because they
do in the summands.
In the cone $C_j\oplus_1SC_j$ projection of the contraction into the complement of a subcomplex of the form  
$(C_{i,j}\oplus S C_{i,j})$
has a single nonzero entry, namely the identity from $(C_j\perp C_{i,j})$ to $(SC_j\perp SC_{i,j})$. The chain maps
in (3) have the same image subcomplexes as the $i\to j$ structure maps in the cone, so the same cancellations
work for them. This shows the $f_{i,j}$ define a $K_1$ morphism. 

Finally  (4) follows from the form of the contraction in the cone since it clearly cancels the entire complex. This
concludes the proof of the lemma. 
\enddemo

As explained above this implies
\proclaim{3.5.2 Corollary} Triangulation of the canonical morphisms of 3.5.1 gives a canonical homotopy of
$\text{id}\oplus S$ to 0.
\endproclaim

\subhead 3.6 Comparison with Volodin $K$-theory \endsubhead 
Volodin \cite{V} developed a version of $K$-theory using triangular matrices. This was shown to be the same as Quillen
$K$-theory by Vasserstein \cite{VS} and Wagoner \cite{W}, see also the  short (or at least highly compressed)
proof by Suslin \cite{S}. Here we  show there is a map from  Volodin's space to $K_1(R)$ so that the induced map on
loop spaces is a homotopy equivalence.  This is the principal ingredient of the coefficient part of the controlled
assembly isomorphism theorem 2.2.1, where a spectrum constructed from $K_1(R)$ is identified as the standard $K$-theory
spectrum. The spectrum-level statement follows from the result here, the relatively formal facts that the map is a map
of spectra, and isomorphism of higher homotopy implies isomorphism of lower homotopy by a suspension trick.

In this section we define Volodin's space and the map, and outline the proof. 

\subsubhead 3.6.1 Definition of $ V_k(R)$\endsubsubhead
We use  the description given by Suslin \cite{S}. For each $k\geq 0$ define a simplicial set $\bar V_k(R)$ whose
$n$-simplices are sequences $g_i$, $0\leq i\leq n$,
\roster\item each $g_i$ is an isomorphism $R^k\to R^k$ 
\item there exists a partial order on the basis so that each $g_ig_j^{-1}$ is  1-triangular with respect to this order.
\endroster
Recall that 1-triangular means triangular with diagonal entries 1. The face map $\partial_i$ is defined by
omitting the $i$ term in the sequence.

Next fix a standard countably generated free module $R^{\infty}$ and consider $R^k$ as the first $k$ coordinates. $R^k$
is then a based submodule of $R^{k+\ell}$, and isomorphisms of $R^k$ can be canonically extended to isomorphisms of
$R^{k+\ell}$ by the identity in the other coordinates. This gives  inclusions $\dots\subset V_k(R)\subset 
V_{k+1}(R)\subset\dots$. The Volodin space is the limit of this sequence:
$$ V(R)=\lim_{n\to\infty} V_k(R).$$

\subsubhead 3.6.2 The map $V(R)\to K_1^V(R)$\endsubsubhead
Suppose $(g_*)$ is an $n$-simplex of $V_k(R)$. Define an $n$-simplex of  $K_1(R)$ by:
\roster \item chain complexes with modules $R^k$ in degrees 1 and 0 and 0 otherwise, with boundary 
$g_i\:R^k\to R^k$ and contraction $g_i^{-1}$; and
\item chain maps 
$$\CD R^k@>1>> R^k\\
@VV{g_i}V @VV{g_j}V\\
R^k@>g_jg_i^{-1}>> R^k\endCD$$
\endroster
Note the cancellation  and image agreement conditions are automatically satisfied because the chain maps are
all isomorphisms. This construction  defines a simplicial inclusion. Further it clearly commutes with the stabilizations
$V_k\subset V_{k+1}$, so induces a map on the direct limit. The main result is:

\proclaim{3.6.3 Theorem} The map  $V(R)\to K_{1}(R)$ induces a homotopy equivalence of loop spaces.\endproclaim

Probably the most efficient way to prove this is to follow Suslin \cite{S} by constructing a space with special
$\pi_1$, showing it is acyclic, and appealing to a characterization theorem. We give a (rather elaborate) direct proof
because some parts of it are useful in other circumstances. 
 In the remainder of this section we set up notation and describe  steps in this proof. Proofs of individual steps are
given in the next four sections. 

The first step is ``compressing into two degrees'' in 3.7. Define $K^V_1(R)$ to be the subcomplex of $K_1(R)$ consisting
of complexes zero except in degrees 0 and 1, then the result is that
$K_1(R)$ deformation retracts to this. The map from the Volodin space maps into this subspace.

The second step in the proof, ``standardizing dimensions'', is given in 3.8. Define $K^V_{(\leq k)}(R)$ to be the
subcomplex of
$K^V_1(R)$ consisting of simplices whose chain modules have dimension
$\leq k$. Similarly $K^V_{(= k)}(R)$ is the subcomplex of chain modules with dimension exactly $k$. We want to show
these two subcomplexes are essentially the same. For convenience we work in the loop space. The result is that the loop
space
$\Omega K^V_{(\leq k)}(R)$ deformation retracts to
$\Omega K^V_{(= k)}(R)$, and the retraction homotopy commutes with loops on  stabilization:
$$\CD\Omega K^V_{(\leq k)}(R)@>{\t retract }>>\Omega K^V_{=  k}(R)\\
@VV{\t include }V@VV{\oplus\left(\smallmatrix 1\endsmallmatrix\right)}V\\
\Omega K^V_{(\leq k+1)}(R)@>{\t retract }>>\Omega K^V_{(=  k+1)}(R)\endCD$$
homotopy commutes. The subspaces $K^V_{(\leq k)}(R)$ form a direct system under inclusion, and $K_1^V(R)$ is the union.
The result shows it is also the limit of $K^V_{(= k)}(R)$ under stabilization. 
The map $V_k(R)\to K^V(R)$ has image in $K^V_{(= k)}(R)$ and these maps give a morphism of direct systems. Therefore to
prove Theorem 3.6.3 it is sufficient to prove the unstable analog: the map $V_k(R)\to K^V_{(= k)}(R)$
induces a homotopy equivalence of loop spaces.

The third step is ``rolling out'' in 3.9, and again takes place in  loop spaces. 1-simplices of $K^V_{(=
k)}(R)$ are triangular chain maps of complexes of length 2, so correspond to changing boundary homomorphisms
by composition on both the right and left by triangular isomorphisms. The Volodin space is constructed using cosets, so
uses composition  only on one side. To fix this we consider a loop $K^V_{(= k)}(R)$, which roughly corresponds to a
sequence of 2-sided compositions giving the identity. But then we can slide things past the identity to move them all
to one side. This deforms the loop space of $K^V_{(= k)}(R)$ to the loop space of $V^{\pm}_k(R)$, the subspace in which
the degree-1 part of chain maps in 1-simplices are identities. 

The final step, ``fixing signs'' in 3.10, shows the loop space of  $V^{\pm}_k(R)$ deformation retracts to the (image of)
loop space of the Volodin space $V_k(R)$. The difference between the two spaces is that $V_k(R)$ is defined using
$+1$-triangular morphisms, while $\pm1$-triangular ones are allowed in $V^{\pm}_k(R)$. The negative signs are a
consequence of using chain complexes, and cannot be avoided. However in the loop space they can be removed essentially
by taking absolute values on the diagonal.

Together these steps prove the comparison Theorem 3.6.3.

\subhead 3.7 Compressing into two degrees\endsubhead 
Define $K^V_1(R)$ to be the subcomplex of $K_1(R)$ consisting of complexes that are trivial in all degrees except 0 and
1. In the next section we analyse this subcomplex; here the objective  is to show
\proclaim{3.7.1 Theorem} There is a canonical deformation retraction of $K_1(R)$ to $K^V_1(R)$.\endproclaim
The proof uses ``degree-wise localization'' of the cancellation of inverses. ``Localization'' will refer to the use of
the steps of 3.4 to cancel some part of a complex  but leave other parts alone. In later sections we use reference maps
to a space to determine which parts to cancel. Here we use degrees, roughly speaking cancelling in positive degrees and
leaving degree 0 alone. 

The proof is a modification of the infinite swindle that contracts $K_1(R)$, but goes through infinite complexes that
are nontrivial in negative degrees. Recall that the suspension $S^jC$ of a complex is defined (3.1.3) by shifting the
grading by $j$. Consider the map that takes $C$ to infinite sum $\oplus_{j\geq 0}S^{-j}C$. We write this sum as $\cdots
\oplus S^{-j}C\oplus\cdots\oplus S^{-1}C\oplus C$ for later notational convenience.
If we associate this in pairs as $\cdots\oplus(S^{-2i-1}C\oplus S^{-2i})\oplus\cdots$ then applying the cancellation
homotopy to the pairs gives a homotopy of the map to the constant 0 map. If we associate as
$\cdots\oplus(S^{-2i}C\oplus S^{-2i+1})\oplus\cdots\oplus C$ and apply the cancellation homotopy then all the
nontrivial suspensions vanish and we get a homotopy to the identity map $C\mapsto C$. Putting these together gives a
nullhomotopy of the inclusion into a space of infinite complexes. Note however that since we used negative suspensions
the nonnegative-degree parts of the intermediate complexes are finite. The plan is roughly to get back into $K_1(R)$ by
discarding everything in negative degrees. We cannot literally do this because truncation does not ``commute'' with
cancellation of inverses. However we can limit the problems to degrees 0 and 1. The outcome will be a homotopy that
agrees with the swindle in degrees 2 and above, and stays in finitely generated nonnegative complexes. The end of the
homotopy is no longer completely trivial, but is concentrated in degrees 0 and 1. Finally we arrange that if a complex
is already concentrated in degrees 0 and 1 then the homotopy leaves it unchanged. This gives the deformation
retraction promised in the theorem. 

To  fill out this outline we must describe  how the cancellations in the infinite swindle are ``localized'' to
positive degrees.  To do this we write out low degree groups in each step, and display chain maps, contractions, etc\.
explicitly as matrices. We leave to the reader the verifications that these actually are chain maps and the
contractions cancel inverses.  As in 3.5 the process is canonical, so it is sufficient to describe the effect on a
single simplex
$(C,\xi)$.

\subsubhead 3.7.2 Stabilization\endsubsubhead 
The first step is to stabilize by including $C$ into the sum with
mapping cones of identity morphisms. The complexes involved are even negative suspensions with negative degrees omitted:
$S_0^{-2i}C$ denotes the complex with degree-$j$ group $C^{j-2i}$ if $j-2i\geq0$, and 0 otherwise. The mapping cone is
the sum of this with its suspension, which is the
$(-2i+1)$ suspension of $C$ truncated in degree 1. In the following diagram the group in a single degree is obtained
by summing a horizontal row, and arrows indicate nontrivial components of the boundary homomorphisms:

$$ \xymatrix@!C{... &C^{2r+2}\ar[d]&C^{2r+1}\ar[d]\ar[dl]&... &C^4\ar[d]&C^{3}\ar[d]\ar[dl]&C^{2}\ar[d]\\
... &C^{2r+1}\ar[d]&C^{2r}\ar[dl]&... & C^3\ar[d]&C^{2}\ar[dl]&C^{1}\ar[d]\\
... &C^{2r}&&... &C^2&&C^{0}
}$$
The boundary in this complex is $\left(\smallmatrix \ddots&&&\\ &c&1&\\&0&-c&\\&&&c\endsmallmatrix\right)$ in degrees
$\geq2$, and $\left(\smallmatrix \ddots&&&\\ &c&1&\\&&&c\endsmallmatrix\right)$ on degree 1. The contraction is
$\left(\smallmatrix \ddots&&&\\ &0&0&\\&1&0&\\&&&\xi\endsmallmatrix\right)$ in degrees $\geq1$ and $\left(\smallmatrix
\ddots&&\\ &0&\\&1&\\&&c\endsmallmatrix\right)$ in degree 0. 

\subsubhead 3.7.3 Reparameterization\endsubsubhead The next step is to use, as much as possible, the isomorphism of
3.5.1 from this to a complex with the same groups and a direct sum boundary structure. The isomorphism in degrees
$\geq1$ is $\left(\smallmatrix \ddots&&&\\ &1&-\xi&\\&0&1&\\&&&1\endsmallmatrix\right)$, and we use the identity in
degree 0. Specifying the chain maps determines the boundary structure in the range. This  is $\pm c$ on the
diagonal (as desired) in degrees
$\geq2$, and
$\left(\smallmatrix \ddots&&\\ &\xi c&&\\&c&&\\&&&c\endsmallmatrix\right)$ in degree 1. The contraction is $\pm \xi$ on
the diagonal in degrees $\geq2$. The contraction is $\left(\smallmatrix \ddots&&&\\
&c\xi&0&\\&0&c\xi&\\&&&c\xi\endsmallmatrix\right)$ in degree 1  and $\left(\smallmatrix \ddots&&\\
&\xi c\xi&\\&\xi c&\\&&\xi\endsmallmatrix\right)$ in degree 0. 

The second morphism goes {\it from\/} a complex with the same groups to this. It is defined like the first morphism but
with the association of terms shifted over by one: $\left(\smallmatrix \ddots&&&&\\
&1&-\xi&&\\&0&1&&\\&&&1&-\xi\\&&&0&1\endsmallmatrix\right)$ in positive degrees, and the identity in degree 0.  We
work out that the boundary homomorphism in the complex is a mapping cone in degrees
$\geq2$:
$\left(\smallmatrix \ddots&&&&\\
&-c&0&&\\&1&c&&\\&&&-c&0\\&&&1&c\endsmallmatrix\right)$, and in degree 1 is $\left(\smallmatrix \ddots&&&\\
&-\xi c\xi&&\\&c&\xi c&\\&&&-\xi c\xi\\&&&c\endsmallmatrix\right)$. Diagrammatically this is 

$$ \xymatrix@!C{...
&C^{6}\ar[d]_{c}\ar[dl]&C^{5}\ar[d]_{-c}&C^4\ar[d]_{c}\ar[dl]_1&C^{3}\ar[d]_{-c}&C^{2}\ar[d]_{c}\ar[dl]_{1}\\ ...
&C^{5}\ar[d]_{c}&C^{4}\ar[dl]_{c\xi}& C^3\ar[d]_{c}\ar[dll]^{-\xi
c\xi}&C^{2}\ar[dl]_{c\xi}&C^{1}\ar[d]_{c}\ar[dll]^{-\xi c\xi}\\ ... &C^{4}&&C^2&&C^{0} }$$
The contraction is the mapping cone contraction in degrees $\geq2$, is $\left(\smallmatrix \ddots&&&&\\
&0&0&&\\&c\xi&0&&\\&&&0&0\\&&&c\xi&0\endsmallmatrix\right)$ in degree 1, and $\left(\smallmatrix \ddots&&&&\\
&-\xi c\xi&0&\xi&\\&&&-\xi c\xi&0&\xi\endsmallmatrix\right)$ in degree 0. 

\subsubhead 3.7.4 Cancellation\endsubsubhead The final step is to note that the contractions allow us to
omit most of this. The residual subcomplex consists of the bottom two groups in even-numbered summands, and inclusion
 into the whole thing is a $K_1$ morphism. 

The residual complex is concentrated in degrees 0 and 1. Explicitly it is: in degree 0,
$\oplus_{* \text{even}} C^*$, in degree 1, $\oplus_{* \text{odd}} C^*$, the boundary is $c-\xi c\xi$, and the contraction
is $\xi c\xi-c\xi c$. We recognize this as a variant on the Whitehead formula for the torsion of a contractible complex.
The morphisms described above give a homotopy from the identity of
$K_1(R)$ to the map that takes
$C$ to the complex given by this formula, so into the subcomplex $K^V_1(R)$ as required. Finally note that if $C$ is
already concentrated in degrees 0 and 1 then it is left unchanged by all of this.

This completes the proof of the two-degree theorem.

\subhead 3.8 Standardizing dimensions \endsubhead In  $K_1$  simplices dimensions are allowed to increase
provided the contraction cancels the new part. Here we deform simplices to ones
where all dimensions are the same, roughly by using the cancellation information to fill out the small ones. 
Define (as in 3.6) $K^V_{(\leq k)}(R)$ to be the
subcomplex of
$K^V_1(R)$ consisting of simplices of complexes with dimension
$\leq k$ in each degree.  $K^V_{(= k)}(R)$ is the subcomplex of complexes with dimension exactly $k$ in
degrees 0 and 1.  

\proclaim{3.8.1 Proposition}  The loop space $\Omega K^V_{(\leq k)}(R)$ deformation retracts to $\Omega K^V_{(= k)}(R)$
so that the retraction homotopy commutes with inclusion and suspension:
$$\CD\Omega K^V_{(\leq k)}(R)@>{\t retract }>>\Omega K^V_{=  k}(R)\\
@VV{\t include }V@VV{\oplus\left(\smallmatrix 1\endsmallmatrix\right)}V\\
\Omega K^V_{(\leq k+1)}(R)@>{\t retract }>>\Omega K^V_{(=  k+1)}(R)\endCD$$\endproclaim

\subsubhead 3.8.2 Details on $K^V_1(R)$\endsubsubhead
The definition of $K_1$ simplices simplifies quite a bit for complexes concentrated in two dimensions. An $n$-simplex
consists of:
\roster\item (short contractible complexes) for $0\leq i\leq n$  an isomorphism $c_i\:C^1_i\to C^0_i$;
\item (chain maps) for $i<j$,  homomorphisms $c^*_{i,j}\:C^*_i\to C^*_j$, where
$*=0,1$, such that the diagram commutes:
$$\CD C^1_{i}@>{c^1_{i,j}}>>C^1_{j}\\
@VV{c_{i}}V@VV{c_{j}}V\\
C^0_{i}@>{c^0_{i,j}}>>C^0_{j}
\endCD$$
\item (partial orders) there exist partial orders on the bases of $C_i$ in each degree so that $c_{i,j}$ is $\pm1$
triangular, and the image of the basis function of $c_{i,j}$ is greater than its complement; 
\item (commutative basis functions) if $i\leq j\leq k$ then $c_{j,k}c_{i,j}$ has the same basis function as $c_{i,k}$;
and 
\item (cancellation condition) the restriction of $c_j$ to the based complements of $c^*_{i,j}$ is $\pm1$ triangular.
\endroster
This relates to definition 3.4.1 as follows: a complex of length 2 is contractible if and only if the boundary is an
isomorphism, and in this case there is a unique contraction (the inverse). Condition (1) thus encodes  contractible
complexes. Condition (2) is equivalent to the $c^*_{i,j}$ defining chain maps. (3) encodes the triangularity
hypothesis and the image order condition. (4) is the same as in 3.4.1. (5) encodes the cancellation property of the
contraction.  

The next result is a version of the ``standard coordinates'' for cancellations given in 3.3.2.

\proclaim{3.8.3 Lemma} Suppose $f\:B\to C$ is a $K^V_1$ morphism (1-simplex) and $\phi\:R^\ell\to C^1$ gives a based
isomorphism to the complement of the image $f(B^1)$. Then the isomorphism $\hat f$ given by $[f^1,\phi]$ in
degree 1 gives a $K^V_1$ commutative diagram (2-simplex)
$$\xymatrix{ B\ar[rr]^f\ar[dr]^{\t inclusion }&& C\\
& B\oplus 1_{R^{\ell}}\ar[ur]_{\hat f}}$$
\endproclaim
More precisely any partial orders on $B$, $C$ making $f$ a 1-simplex determines one on the
stabilization making the diagram a 2-simplex. In degree 0 $\hat f$ is given by $[f^0,c\phi]$.

To see $K^V_1$ structure we decompose $C$ as the image of $f$ plus a cancellation part: $C^1=\bar C^1\oplus D$ and
$C^0=\bar C^0\oplus \bar D$. The chain map $\hat f$ in these coordinates is:
$$\xymatrix{B^1\oplus R^{\ell}\ar[rr]^{\left(\smallmatrix 
f^1&0\\0&\phi\endsmallmatrix\right)}\ar[d]^{\left(\smallmatrix  b&0\\0&1\endsmallmatrix\right)}&&{\bar
C}^1\oplus D\ar[d]^{\left(\smallmatrix  \bar c&x\\0&\delta\endsmallmatrix\right)}\\
B^0\oplus R^{\ell}\ar[rr]^{\left(\smallmatrix  f^0&x\phi\\0&\delta\phi\endsmallmatrix\right)}&&{\bar C}^0\oplus {\bar
D}}$$
Since $\phi$ is a based isomorphism it is 1-triangular with respect to any partial order. Use the partial order on
$B^1\oplus R^{\ell}$ that makes the basis map of $\hat f^1$ order-preserving. This extends the partial order on $B^1$
and makes $\hat f^1$ $\pm1$ triangular because $f^1$ is. Further $\hat f^0$ is $\pm1$ triangular because the
basis of
$R^{\ell}$ preceeds the basis of
$B^0$ (by 3.8.2(3)) and $f^0$, $\delta$ and $\phi$ are all $\pm1$ triangular. This means $\hat f$ is a $K^V_1$ morphism
with respect to the given partial orders. Finally the diagram in the lemma commutes, so by 3.2.3 in each degree the
underlying basis functions commute. This means the diagram is a 2-simplex in $K^V_1$.

\subsubhead 3.8.4 Loop spaces \endsubsubhead The model we use for the loop space uses
 triangulations of $\Delta^n\times I$ such that the projection to $\Delta^n$ is simplicial.  There is a standard
minimal such triangulation described in 3.4.3. All such triangulations have  qualitative structure similar to the
standard one: there is a linear ordering of the $(n+1)$-simplices so adjacent ones share a face, and the projection of
each one to $\Delta^n$ collapses a single edge. 

If $K$ is a simplicial space define the loop space $\Omega( K,*)$ to be the simplicial space with  $n$-simplices 
simplicial maps of triangulations of $\Delta^n\times I$ as above, with $\Delta\times\{0,1\}$ mapping to the basepoint.
Boundaries of simplices are defined by restricting maps to faces
$(\partial_i\Delta^n)\times[0,1]$.

\subsubhead 3.8.5 Proof of Proposition 3.8.1\endsubsubhead
Take the identity map $1\:R^k\to R^k$ as the basepoint of both $K^V_{\leq k}$ and $K^V_{= k}$. We construct a
retraction $S\:\Omega K^V_{\leq k}\to K^V_{= k}$ and a morphism $s\:\t id \to S$. The deformation of the lemma is
obtained by subdividing the morphism. The retraction and morphism are constructed by induction on skeleta. More
precisely we construct them on the vertices, then extend to the 1-simplices. Higher dimensional simplices are defined to
be collections of 1-simplices satisfying some coherence conditions, so no further constructions are required, but
coherence must be checked.

A vertex of  
$\Omega K^V_{(\leq k)}$ is a simplicial map of a triangulation of $I$ into $K^V_{(\leq k)}$, with the ends going to
$1_{R^k}$. Such a map is a sequence of objects $C_0,\dots C_n$ with morphisms between adjacent ones: either $C_i\to
C_{i+1}$ or $C_i\leftarrow
C_{i+1}$. Define  objects $\hat C_i$ by stabilizing the $C_i$ to have dimension $k$: 
$\hat C_i=C_i\oplus (1\:R^\ell\to R^\ell)$, where  $\ell +\t dim (C_i^1)=k$. We want to fill in the dotted arrows in
the following diagram so the squares are $K^V_1$ morphisms (including coherent partial orders): 
$$\xymatrix{ C_0\ar[r]\ar[d]^{=}&{\dots}
\ar[r]&C_i\ar[d]^{\subset}&C_{i+1}\ar[l]\ar[r]\ar[d]^{\subset}&{\dots} &
C_n\ar[l]\ar[d]^{=}\\ {\hat C}_0\ar@{.>}[r]&{\dots} \ar@{.>}[r]&{\hat C}_i&{\hat
C}_{i+1}\ar@{.>}[l]\ar@{.>}[r]&{\dots} & {\hat C_n}\ar@{.>}[l]}$$
Consider a square in which the morphism goes from $C_i$ to $C_{i+1}$, and choose appropriate partial orders. Partial
order the bases in $\hat C_{i+1}$ by putting the new elements first, then the inclusion is a morphism and the
composition (diagonal in the following) is also a morphism.
$$\xymatrix{ C_i\ar[d]^{\subset}\ar[r]\ar[dr]&C_{i+1}\ar[d]^{\subset}\\ {\hat
C}_i\ar@{.>}[r]&{\hat C}_{i+1}}$$
According to Lemma 3.8.3 a based isomorphism from the stabilization factor $R^{\ell}$ in $\hat C^1_i$ to the complement
of the image in $\hat C^1_{i+1}$ determines a way to complete the lower triangle. There is such an isomorphism for
dimension reasons; choose one randomly and use that.

Squares in which the morphism goes from $C_{i+1}$ to $C_i$ are filled in with the mirror image of this argument.

Note the ends of the sequence already have
dimension $k$ so do not change. This means the lower row in the diagram gives a sequence that begins and ends with
identities, so defines a vertex in
$\Omega K^V_{(= k)}$. The rest of the diagram defines a morphism from the original vertex in $\Omega K^V_{(\leq
k)}$ to this. This defines the retraction
$S$ and morphism $s$ on vertices.

Now consider a 1-simplex in $\Omega K^V_{(\leq k)}$. 
This is a map of a triangulation of $\Delta^1\times I$ into $K^V_{(\leq k)}$, so consists of two vertices as above and a
sequence of 2-simplices filling in between them. Denote the vertices by $C_*$ and $D_*$, then a typical segment looks
like
$$\xymatrix{ {\dots}
\ar[r]&C_i\ar[d]\ar[dr]&C_{i+1}\ar[l]\ar[r]\ar[d] &
C_{i+2}\ar[l]\ar[d]\ar[r]\ar[dl]&{\dots}\ar[dl]\\ {\dots}
\ar[r]&D_j\ar[r]&D_{j+1}\ar[r]& D_{j+2}&{\dots}\ar[l]}$$
The retraction and morphism are already defined on vertices, so have sequences $\hat C_*$, $\hat D_*$ and morphisms
from the originals to these. To extend $S$ to the 1-simplices we need to fill in 2-simplices joining $\hat C_*$ and
$\hat D_*$. We do this by working along the loop coordinate. Note since loops begin at the identity the identity $\hat
C_0=\hat D_0$ fills in the beginning of the diagram. Suppose by induction that we have extended it up to a 2-simplex
with one vertex on $C_i$, so have everything but the dotted arrow in the following:
$$\xymatrix{ &&C_i\ar[dll]\ar[dl]\ar[d]\\
D_j\ar[r]\ar[d]&D_{j+1}\ar[d]&{\hat C}_i\ar[dll]\ar@{.>}[dl]\\
{\hat D}_{j}\ar[r]&{\hat D}_{j+1}&
}$$
Again 3.8.3 can be used to fill this in. Compose $C_i\to D_{j+1}\to \hat D_{j+1}$ to obtain the setting of 3.8.3. Now
however we use a particular basis-preserving isomorphism from the stabilization factor in $\hat C_i$ to the complement
in degree 1 in $\hat D_{j+1}$. Namely, the composition $\hat C_i\to \hat D_{j+1}\to \hat D_{j+1}$ is a $\pm1$
triangular isomorphism, and the underlying basis function preserves images, etc\. so gives a bijection. Applying 3.8.3
using this gives a dotted arrow so that the right square commutes and basis functions commute. This fills in this
triangle, and we move on to the next. Triangles with two vertices on the $C_*$ are handled similarly.

Eventually this process reaches the end of $\Delta^1\times I$. There nothing happens since the dimension is already
$k$, so the diagram ends with identities. The $\hat C$, $\hat D$ part of the diagram gives a 1-simplex in $\Omega
K^V_{(=k)}$ and we define this to be the image of  $S$. The morphisms going from the original to the new part give a
morphism from the original 1-simplex to the new one, and this is defined to be $s$.

Finally consider an $n$-simplex of $\Omega
K^V_{(\leq k)}$. The construction on 0- and 1-simplices gives all the objects and maps required to define $S$ and $s$
on this simplex, but we need to check $\pm1$ triangularity with respect to coherent partial orders, and commutativity
of underlying basis functions. The ability to choose coherent partial orders follows from the refinement noted after
the statement of 3.8.3 that any suitable partial orders on the input extends to ones on the output. To see the
underlying basis functions commute it is sufficient to check 2-simplices. There it is easily seen to hold by the choice
made in construction of 1-simplices. 

This $S$ maps to $\Omega
K^V_{(=k)}$. It is a retraction (identity on $\Omega
K^V_{(=k)}$) because nothing happens when dimensions are already $k$. This completes the proof of the
dimension-standardizing proposition.

\subhead 3.9 Rolling out\endsubhead
Recall (\S 3.6) that $V_k^{\pm}$ is defined like $V$ but with $\pm1$ triangular automorphisms of $R^k$. It maps to
$K^V_{(=k)}$, essentially to the subspace of short chain complexes in which the degree-1 part of chain maps are
base-preserving isomorphisms. The basic difference between $K^V_{(=k)}$ and $V_k^{\pm}$ is that in the former the
morphisms change objects by composition on both sides with triangular matrices, whereas in the latter compositions are
only on one side.

\proclaim{3.9.1 Proposition} The inclusion of loop spaces $\Omega V_k^{\pm}\to \Omega K^V_{(=k)}$ is a homotopy
equivalence of realizations. \endproclaim

The plan of the proof is to define a filtration $\Omega V_k^{\pm}=X_0\subset X_1\subset\dots$ with $\cup_iX= \Omega
K^V_{(=k)}$. We then construct a homotopy from the identity of $K^V_{(=k)}$ to a map $S$ so that $S$ decreases
filtration degrees (i.e\. $S(X_i)\subset X_{i-1}$) and the homotopy increases degrees in a controlled manner (the image
of the homotopy on $X_i$ is in $X_{2i}$). The proposition follows routinely from this.

\subsubhead 3.9.2 The filtration\endsubsubhead
An $n$-simplex of $\Omega K^V_{(=k)}$ is a special triangulation of $\Delta^n\times I$ and a simplicial map of this to
$K^V_{(=k)}$ that takes the ends to the basepoint $\t id \:R^k\to R^k$. For instance a vertex is a triangulation of $I$
and a map to $K^V_{(=k)}$, so a linear sequence of short chain complexes and maps. Since the maps are all invertible,
for notational convenience we can take them all to go in the direction of the loop coordinate. Explicitly we have
$k$-dimensional based modules and isomorphisms:
$$\xymatrix{R^k\ar@2{-}[r]\ar@2{-}[d]&C^1_0\ar[r]^{t_0}\ar[d]^{c_0}&C^1_1\ar[d]^{c_1}\ar[r]^{t_1}&{\dots}\ar[r]^{t_{n-1}}&C^1_n\ar@2{-}[r]\ar[d]^{c_n}&R^k\ar@2{-}[d]\\
R^k\ar@2{-}[r]&C^0_0\ar[r]^{t_0}&C^0_1\ar[r]^{t_1}&{\dots}\ar[r]^{t_{n-1}}&C^0_n\ar@2{-}[r]&R^k
}$$

The $n+1$-simplices in these triangulations of $\Delta^n\times I$ have a
linear order so that adjacent ones share a face, the first has $\Delta^n\times\{0\}$ as a face, and the last has
$\Delta^n\times\{0\}$ as a face. We define a map from this triangulation into $K^V_{(=k)}$ to be in the subspace $X_r$
if all but the last $r$ of these $n+1$-simplices map into $V^{\pm}_k$. 
\subsubhead 3.9.3 The homotopy\endsubsubhead
This construction is canonical so it is sufficient to describe the effect on a simplex. In fact we only describe it on
vertices. It is straightforward to extend this description to 1-simplices, but the diagrams are too complicated to be
useful and are best drawn for oneself. The general case follows from 1-simplices and easily-seen coherence of partial
orders.

Begin with a simplex, which is a diagram as in 3.9.2. The first step in the homotopy doubles the length of the diagram
by alternating the given chain maps with identities. The degree-1 terms of the chain maps are shifted one place to the
right by starting with an identity. There is a ``morphism'' (a 1-simplex) going from this new vertex to the old one,
giving the homotopy. In the diagram the front row is the original vertex, the back row is the new one, and the
morphisms joining them give the homotopy. Identity maps are denoted by long equality signs.
$$
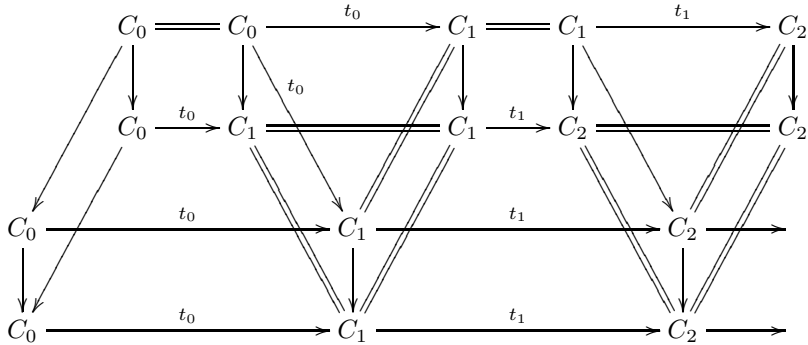
\xymatrix{
&C_0\ar@2{-}[r]\ar[ddl]\ar[d]&C_0\ar[rr]^{t_0}\ar[ddr]^<<<<<<<<{t_0}\ar[d]
&&C_1\ar@2{-}[r]\ar@2{-}[ddl]\ar[d]&C_1\ar[rr]^{t_1}\ar[ddr]\ar[d]&&C_2\ar[d]\ar@2{-}[ddl]
\\ &C_0\ar[r]^{t_0}\ar[ddl]&C_1\ar@2{-}[rr]\ar@2{-}[ddr]
&&C_1\ar[r]^{t_1}\ar@2{-}[ddl]&C_2\ar@2{-}[rr]\ar@2{-}[ddr]&&C_2\ar@2{-}[ddl]\\ 
C_0\ar[rrr]^{t_0}\ar[d]&&&C_1\ar[rrr]^{t_1}\ar[d]&&&C_2\ar[d]\ar[r]&\\
C_0\ar[rrr]^{t_0}&&&C_1\ar[rrr]^{t_1}&&&C_2\ar[r]&
}$$
Note that the stretching takes vertices in the filtration $X_r$ to $X_{2r-1}$. 

The next step is to apply the inverse of this move on the segment omitting the first and last squares. The result is a
sequence with length one more than the starting sequence and with the degree-1 terms of the chain maps shifted right one
position. This sequence begins with a chain map that is the identity in degree 1, and ends with one that is the
identity in degree 0. 

Repeat this sequence of moves. This takes us to a sequence of length greater by 2, with the last two squares the
identity in degree 0, explicitly the sequence ends with 
$$\xymatrix{{\dots}\ar[r]^{t_{n-3}}&C^1_{n-2}\ar[rr]^{t_{n-2}}\ar[d]^{t_{n-1}t_{n-2}c_{n-2}}&&C^1_{n-1}\ar[d]^{t_{n-1}c_{n-1}}\ar[rr]^{t_{n-1}}&&C^1_n\ar[d]^{c_n}\\
{\dots}\ar[r]^{t_{n-1}}&C^0_n\ar@2{-}[rr]&&C^0_n\ar@2{-}[rr]&&C^0_n
}$$
This vertex is in the same filtration degree as the original one because the filtration is defined using the end of the
degree-1 part of the chain maps, and this has not changed.

The final step in the homotopy shortens this diagram by 1 by moving the end degree-1 term down to degree 0. In the
diagram the beginning sequence is in the back, the end in the front.
$$\xymatrix{&{\dots}\ar[r]^{t_{n-3}}&C^1_{n-2}\ar[rr]^{t_{n-2}}\ar[d]\ar@2{-}[ddl]&&C^1_{n-1}\ar[d]_{t_{n-1}c_{n-1}}\ar@2{-}[ddr]\ar[rr]^{t_{n-1}}&&C^1_n\ar[d]^{c_n}\ar[ddl]^>>>>>>{b(t_{n-1})^{-1}}\\
&{\dots}\ar[r]&C^0_n\ar@2{-}[rr]\ar@2{-}[ddl]&&C^0_n\ar@2{-}[rr]\ar[ddr]_>>>>>>>{(c_nt_{n-1})^{-1}}&&C^0_{n}\ar[ddl]^{b(c_nt_{n-1})^{-1}}\\
{\dots}\ar[r]^{t_{n-3}}&C^1_{n-2}\ar[rrrr]^{t_{n-2}}\ar[d]^{t_{n-1}t_{n-2}c_{n-2}}&&&&C^1_{n-1}\ar@2{-}[d]&\\
{\dots}\ar[r]^{t_{n-1}}&C^0_n\ar[rrrr]^{(t_{n-1}c_{n-1})^{-1}}&&&&C^1_{n-1}&
}$$
On the right-hand slanted arrows ``$b(*)$'' denotes the underlying basis function of the triangular map $*$. Note $t_*$
is triangular because it is part of a chain map in a 1-simplex, and $c_n$ is a based isomorphism (in fact identified
with $\t id \:R^k\to R^k$). The diagrams all commute except for the triangles on the right which are off by the
triangular parts of the inverses of $t_{n-1}$ and $c_nt_{n-1}$ respectively. Recall, however, that 2-simplices in
$K^V_{(=k)}$ are not required to commute. The requirements are that triangularity conditions should all hold with
respect to fixed choices of partial order, and the underlying basis functions should commute. These conditions clearly
do hold. The diagram therefore defines a 1-simplex. The new vertex at the end of the simplex has one fewer
non-basis-preserving map in degree 1, so is in lower filtration degree than the original. This accomplishes the
objective of the construction.

As noted above it is straightforward to extend this to 1-simplices, but the diagrams have twice as many terms so we
will not attempt them here. Higher simplices are defined by collections of 1-simplices satisfying coherence
conditions on partial orders and basis functions. Since the construction is canonical this follows formally. 

\subhead 3.10 Fixing signs\endsubhead
This section gives the last step in the proof of Theorem 3.6.3:
\proclaim{3.10.1 Proposition} The loop space $\Omega V^{\pm}_k$ deformation retracts to $\Omega
V_k$. \endproclaim
Recall that an $n$-simplex in $V^{\pm}_k$ is a sequence of automorphisms of $R^k$ (matrices) $(g_0,\dots,g_n)$ so that
all products $g_jg_i$ with $j>i$ are all $\pm1$ triangular with respect to some partial order of the basis of $R^k$.
The next lemma asserts that we can change the signs on the diagonal arbitrarily.
\proclaim{3.10.2 Lemma} Suppose $(g_i)$ is an $n$-simplex of $V^{\pm}_k$ and $(\epsilon_0, \dots ,\epsilon_n)$ are
$\pm1$ diagonal matrices. Then $(\epsilon_ig_i)$ is also an $n$-simplex and $(g_i)\to(\epsilon_ig_i)$ is a morphism of
simplices.\endproclaim
Here ``morphism'' of simplices should be interpreted using the images in $K^V$. More directly it means it gives a
homotopy of the simplex defined on the standard triangulation of $\Delta^n\times I$, with $n+1$-simplices of the form
$(g_0,\dots,g_{i-1},\epsilon_ig_i,\dots,\epsilon_ng_n)$. To check these are simplices use a partial order that
makes $(g_i)$ a simplex. The products in one of these $n+1$-simplices are of the form $g_jg_i^{-1}$,
$\epsilon_jg_jg_i^{-1}$, and $\epsilon_jg_jg_i^{-1}\epsilon_i^{-1}$. These are all $\pm1$ triangular because the
$g_jg_i^{-1}$ are and the $\epsilon_j$ are $\pm1$ triangular with respect to any partial order. 

\subsubhead 3.10.3 Proof of the Proposition\endsubsubhead
We begin by defining the map $S$ and homotopy $s$ on vertices. A vertex of $\Omega V^{\pm}_k$ is a sequence $(g_0,\dots
g_r)$ of $k\times k$ matrices with $g_0=\t id =g_r$ and for each $i$ there is a partial order in which
$g_ig_{i-1}^{-1}$ is
$\pm1$ triangular. Inductively define $\epsilon_i$ by $\epsilon_0=\t id $, and $\epsilon_i=\t diag
(g_ig_{i-1}^{-1}\epsilon_{i-1}^{-1})$. The matrix on the right is $\pm1$ triangular with respect to some partial
order, so the diagonal part is a well-defined $\pm1$ diagonal matrix (Lemma 3.2.4). Define $S$ on this path by 
$$S(g_*)=\cases (\epsilon_0g_0,\dots,\epsilon_rg_r)&\text{\quad if }\epsilon_n=\t id \\
(\epsilon_0g_0,\dots,\epsilon_rg_r,\t id )&\text{\quad if not}\endcases$$
According to the lemma this is a path in $V^{\pm}_k$, and there is a morphism (standard homotopy) from the original
path to the new one. (The second case requires sticking on a trivial 2-simplex at the end of the path). 

The choice of the $\epsilon_i$ implies the image path is actually in $V_k$, and is unchanged if the original path were
already in this subspace.

Lemma 3.10.2 implies that this map and homotopy extend to higher simplices of $\Omega V^{\pm}_k$. This is because the
higher simplices are determined by their vertices, on which the map is defined, and triangularity of products with
respect to a single partial order. The lemma asserts that this triangularity holds. 

The extension to all simplices of $\Omega V^{\pm}_k$ gives a map into  $\Omega V_k$ that is the identity on $\Omega
V_k$, and a homotopy of this to the identity in $\Omega V^{\pm}_k$. This is the conclusion of Proposition~3.10.1.

\head 4. Controlled $K$-theory\endhead
In this section we give the definition and naturality properties that do not require spacial localization. General
definitions for controlled algebra are in \S4.1, and differ in many details from versions used earlier. These differences
do not change the final outcomes but make some arguments easier. The controlled version of triangularity is described in 4.2. The
controlled $K_1$ space is defined in \S4.3, essentially by inserting many $\epsilon$s into the uncontrolled definition.
Naturality with respect to restriction and uniformly continuous functions are immediate consequences of the definition. Dependence
on the map $E\to X$ is also analysed in this section. The inverse limit as $\epsilon\to0$ is introduced in Section 4.4.

\subhead 4.1 Controlled algebra\endsubhead
The definitions we use for controlled algebra  are  different from other versions, particularly in the systematic use
of paths. An earlier version of path-based definitions is in \cite{Quinn 6}, see also \cite{Anderson-Munkholm}.
Metric notations
are given in
\S4.4.1, modules described in 4.1.2, and homomorphisms defined in \S4.1.3.

\subsubhead 4.1.1 Metric notions \endsubsubhead
Suppose $X$ is a metric space with metric $d$. We say a path $\gamma\:[0,1]\to X$ {\it has radius\/} $<\epsilon$ if
$d(\gamma(s),\gamma(t))<\epsilon$ for all $s,t\in [0,1]$.  

Now suppose $\epsilon>0$ and 
$Y\subset X$.  Define the ``$\epsilon$ (path) enlargement''
$Y^{\epsilon}$ to be all points that can be reached from $Y$ by a path of radius $<\epsilon$. Previously we have used the
notation $Y^{\epsilon}$ to denote the points with distance less than $\epsilon$ from $Y$, but the truth is that we don't care
how close a point is if you can't get there by a small path. It is possible to redefine the metric using paths so distances
encode path properties and the two definitions coincide. We prefer not to do this because it complicates consideration of
subspaces.

 In the opposite direction define the ``$\epsilon$ (path) reduction'' to be points that cannot be reached from the
complement of
$Y$ by a path of radius $<\epsilon$. Note $Y^{-\epsilon}=X-(X-Y)^{\epsilon}$. 

The ambient space $X$ is not displayed in this notation. It is usually clear from the context but occasionally, as in the
 lemma below, will require comment. The benefit of using path enlargements is that they are less dependent
on the ambient space than distance enlargements, as the following ``excision'' lemma shows:
\proclaim{Lemma} Suppose $U$ is open in $X$. Then $(X-U)^{\epsilon}\cap U = (\bar U- U)^{\epsilon}\cap U$, where $\bar U$
denotes the closure in $X$, the path enlargement on the left is taken in $X$, the one on the right in
$\bar U$. 
\endproclaim
\demo{Proof} Since $(\bar U, \bar U-U)\subset (X,X-U)$ the containment $(X-U)^{\epsilon}\cap U \supset (\bar U-
U)^{\epsilon}\cap U$ is clear. To show the other direction suppose $x\in (X-U)^{\epsilon}\cap U$. This means there is a path
$\gamma\:[0,1]\to X$ of radius $<\epsilon$ starting at $x$ and ending in $X-U$. Let $t_0$ be the minimum of $\{t\mid
\gamma(t)\notin U\}$. Since $U$ is open $\gamma(t_0)\notin U$, but it is in the closure. Thus $\gamma|[0,t_0]$ is a path of
radius $<\epsilon$ in $\bar U$ starting at $x$ and ending in $\bar U-U$. This shows $x\in (\bar U-U)^\epsilon$, as required.
\enddemo
\subsubhead 4.1.2 Metric frontiers\endsubsubhead
The {\it metric frontier\/} of a metric space is the complement of $X$ in its completion, and is denoted by $\t Fr (X)$.
Technically this is the set of limit points of Cauchy sequences in $X$ that do not converge in $X$ itself. We use the shorthand
notation $\t Fr ^{\epsilon}X$ for $X\cap(\t Fr (X))^{\epsilon}$, where the
enlargement $(*)^{\epsilon}$ is taken in the metric closure.  $\t Fr
^{\epsilon}X$ can also be described as starting points of  non-convergent Cauchy arcs $\gamma\:[0,1)\to X$ of
radius
$<\epsilon+\delta$ for some positive $\delta$. Note that if $X$ is complete (e.g\. compact) then the metric frontier is
empty.

\proclaim{Lemma} \roster\item
If\/ $U\subset X$ is open then $(\t Fr ^{\epsilon}X\; \cup (\bar U-U)^{\epsilon})\cap U = \t Fr ^{\epsilon}U$.
\item if $f\:X\to Y$ is a map of metric spaces such that $d(x,x')<\delta$ implies $d(f(x),f(x'))<\epsilon$, then
$f(W^{\delta})\subset (f(W))^{\epsilon}$.
\item if $f\:X\to Y$ is proper and uniformly continuous then it extends to a map of completions with $\bar f(\t Fr
X)\subset \t Fr Y$.
\endroster
\endproclaim Note that if the hypotheses of both (2) and (3) hold then $f(\t Fr ^{\delta}X)\subset\t Fr ^{\epsilon}Y)$.

Statement (1) is the Lemma of the previous section applied in the metric closure of $X$. 
Statement (2) is clear because $f$ takes paths of radius $<\delta$ to paths of radius $<\epsilon$.

In (3) uniform continuity implies Cauchy paths $[0,1)\to X$ are taken to Cauchy paths in $Y$, so the map extends to a map
of completions. If the image of a Cauchy path converges in $Y$ then by properness the inverse image of the closure is
compact in $X$. But a Cauchy path in a compact space converges. Thus if the original path does not converge in $X$, neither
does the image, and so completion points in $X$ map to completion points in $Y$.
\subsubhead 4.1.3 Geometric modules and morphisms\endsubsubhead
The setting is a ``coefficient map'' $p\:E\to X$ with $X$ a metric space. Distances are measured in $X$, while $E$ serves
to record fundamental group data.

 A {\it geometric $R$-module\/} on $p\:E\to X$ is a free module with basis
$S$, and a function $s\:S\to E$. 

A {\it morphism\/} of geometric modules is a formal linear combination of
paths in $E$. Specifically $f\:R[S_1]\to R[S_2]$ is of the form $\Sigma_i
r_i\gamma_i$, where 
\roster\item  $r_i\in R$;
\item $\gamma_i$ consists of points $\gamma_i(0)\in S_1$,
$\gamma_i(1)\in S_2$, and a path $\gamma\:[0,1]\to E$ from the image of $\gamma_i(0)$
to the image of $\gamma_i(1)$ such that
\item only finitely many of the $\gamma_i$ begin at any given element of
$S_1$.
\endroster 
We consider two morphisms to be the same if they differ by the operations
\roster\item omit a path with zero coefficient, or
\item if a path occurs multiple times, replace with a single copy and add the coefficients.
\endroster

Note that a geometric morphism induces a homomorphism on the free modules
generated by the basis sets. In some cases the geometric morphism can be
essentially reconstructed  from the homomorphism, but even in these cases
the geometric version works better in constructions.

Geometric morphisms can be composed in the evident way: compose paths
when possible, and multiply coefficients. This composition can be made
associative by interpreting ``path'' to mean a map of $[0,t]$ for some
$t\geq0$ rather than restricting to $[0,1]$. Composed paths are then
defined on longer intervals, and do not require  reparameterization
to get back to $[0,1]$.

A geometric morphism is said to have {\it radius less than\/} $\epsilon$
if all paths with nonzero coefficients have images in $X$ of radius $<\epsilon$. 

Morphism radius is subadditive in compositions: if $f\:A\to B$ has radius
$\leq\epsilon$ and $g\:B\to C$ has radius $\leq \delta$ then $gf$ has
radius $\leq \epsilon+\delta$. In particular $\epsilon$ morphisms do not form a category because compositions are usually
bigger than
$\epsilon$. The inability to use category methods makes
controlled algebra considerably more complicated than uncontrolled
versions. Some inverse limit versions can be made categorical
\cite{Pedersen}, but they are not well-adapted to applications.

If $Y\subset X$ then the {\it restriction\/} $f|Y$ consists of the paths of $f$ that lie in $p^{-1}(Y)$. Restriction is
functorial in the sense that composition of restrictions is the restriction of the composition. 

\subsubhead 4.1.4 Homotopy of morphisms \endsubsubhead
A {\it homotopy\/} of a geometric morphism $\Sigma r_i\gamma_i$ consists
of homotopies (rel ends) of all the paths $\gamma_i$. A homotopy
determines another morphism by using the same coefficients and the paths
at the end of the homotopies. Two morphisms are homotopic if (after combinations and omission of paths
with zero coefficient) they can be related in this way.

A homotopy of geometric morphisms has {\it radius less than\/} $\epsilon$
if each of the homotopies individually has radius less than $\epsilon$.
Generally a homotopy $Y\times I\to X$ has radius less than $\epsilon$ if
the restriction to each arc $\{y\}\times I$ has radius
$<\epsilon$. Explicitly this means $h\:Y\times I\to X$ satisfies $d(h(y,s),h(y,t))<\epsilon$ for all $y\in Y$, $s,t\in I$.  

Note that a $\delta$ homotopy of an $\epsilon$ morphism gives a morphism
of radius less than $\epsilon+2\delta$.  Reversing an $\epsilon$ homotopy from $f_0$  to
$f_1$ gives an $\epsilon$ homotopy from $f_1$ to $f_0$. Homotopies can be composed, and
radii add. The additivity means ``$\epsilon$ homotopy'' is not an equivalence relation, and we cannot work with
 ``homotopy classes''. As with the failure to be a category
this requires us to be more direct and explicit in constructions.

The {\it restriction\/} of a homotopy to $Y\subset X$  consists of all the
individual path homotopies that lie in
$Y$. Restriction and homotopy do not play well together, but  homotopies of radius $<\epsilon$ give predictable results over the
reduced set $Y^{-\epsilon}$. 

\subsubhead 4.1.5 Uncontrolled morphisms \endsubsubhead
Note that if no size conditions are imposed then morphisms on $E\to X$ are a category, and homotopy does give an
equivalence relation. In this case $X$ becomes irrelevant and choice of a basepoint gives an equivalence between homotopy
classes of $R$-morphisms and homomorphisms over the ring $R[\pi_1E]$ \cite{Quinn 6}. ``Forgetting control'' therefore gives
a map of controlled algebra into a classical algebraic setting.

\subsubhead 4.1.6 Main example: geometric cellular chains \endsubsubhead
Let $X$ be a CW complex with skeleta $X^{(n)}$. We suppose $X$ is slightly generic in the sense that attaching maps for
$n+1$-cells
$S^n\to X^{(n)}$ are differentiable with invertible derivative at points mapping to centers of $n$-cells. This means these
points are isolated, and the sign of the derivative associates  $\pm1$ to each one.  An arbitrary CW complex can be made
generic by arbitrarily small perturbations of the attaching maps, so the conclusions apply generally.

For each integer $n\geq0$ define a geometric $Z$-complex $C_n$ by taking the basis to be the centers of the $n$-cells. 

Define a geometric morphism $\partial\:C_n\to C_{n-1}$ by: the paths from the center of an $n$-cell are rays to points on
$S^{n-1}$ that map to centers of $n-1$-cells under the attaching map. The coefficient on a path is the sign of the determinant
of the derivative of the attaching map at the endpoint. The groups and homomorphisms obtained from these geometric objects are
exactly the usual cellular chain complex for $X$. However in the geometric version $\partial^2$ is plainly not trivial.
\proclaim{Lemma} The composition $\partial^2$ in this complex is homotopic to $0$. If each cell has
radius
$<\epsilon$ then $\partial$ has radius $<\epsilon$ and the nullhomotopy has radius $<2\epsilon$.\endproclaim
\demo{Indication of proof} (See \cite{Quinn 6} for details.) Suppose $\sigma$ is an $n$-cell and $\tau$ an $n-2$-cell. The usual
geometric proof that algebraically
$\partial^2=0$ shows the composed paths from center of $\sigma$ to center of $\tau$ occur in pairs with opposite signs. It also
shows that these pairs are homotopic, by a homotopy that is very close to the union of $\sigma$, $\tau$, and the intermediate
$n-1$-cells. These homotopies give a homotopy of $\partial^2$ to a morphism with half as many paths, and coefficient $1+(-1)$
on each one. We can omit all these to get the 0 geometric morphism. The estimates follow from the locations of the paths and
homotopies.
\enddemo

\subsubhead 4.1.7 Controlled chain notions\endsubsubhead
The definition abstracts what we would get by restricting the example above to an
open set. Alternatively we ``add control'' to the usual definitions following what will become a familiar pattern: modules
become geometric modules, homomorphism become geometric morphisms, identities (such as
$\partial^2=0$) become homotopies that are only assumed to exist, and identities are allowed to fail near the metric
frontier. In the following we fix 
$p\:E\to X$ and $\epsilon>0$.

An $\epsilon$ {\it chain complex\/} over $p$ consists of
\roster\item geometric modules $C^i$ for $i$ nonnegative integers;
\item geometric $R$-morphisms $c\:C^i\to C^{i-1}$ with radius $<\epsilon$; and
\item there exist (for each $i$) a homotopy  of radius $<\epsilon$ from $c^2$ to a morphism $C^i\to C^{i-2}$ that is
trivial on basis elements outside $\t Fr ^{3\epsilon}X$.
\endroster

Similarly an $\epsilon$ {\it chain map\/} $f\:(C,c)\to (D,d)$ consists of 
\roster \item for each $i$ a geometric morphism $f\:C^i\to D^i$ of radius $<\epsilon$;
\item for each $i$ there exist homotopies of radius $<\epsilon$ from $d\,f\:C^i\to D^{i-1}$ to morphisms that agree with $fc$ on basis
elements outside $\t Fr ^{\epsilon}X$.
\endroster

Finally an $\epsilon$ {\it chain contraction\/} $\xi$ for a complex $(C,c)$ is
\roster\item geometric morphisms $\xi\:C^i\to C^{i+1}$ of radius $<\epsilon$ so that
\item there exist $\epsilon$ homotopies starting with $c\xi+\xi c\:C^i\to C^i$ and going to morphisms that are the identity on basis
elements not in $\t Fr ^{3\epsilon}$.
\endroster

\subsubhead 4.1.8 $\epsilon$ isomorphisms \endsubsubhead
 The definition here differs from the usual one by allowing
misbehavior near the metric frontier. 

Suppose $f\:A\to B$ is a geometric morphism over $p\:E\to X$. $f$ is an $\epsilon$ {\it isomorphism\/} if
\roster\item $f$ has radius $<\epsilon$;
\item there is $f^-\:B\to A$ with radius $<\epsilon$ such that
\item there are $\epsilon$ homotopies $h_A$, $h_B$ starting with $f^-f$ and $ff^-$ respectively and ending with morphisms that
are identities when restricted to $A|(X-\t Fr ^{3\epsilon}X)$, $B|(X-\t Fr ^{3\epsilon}X)$ respectively.
\endroster
\proclaim{Lemma} Suppose $f\:A\to B$ is an $\epsilon$ isomorphism over $p\:E\to X$, and $U\subset X$ is open. Then the
restriction
$f|U\:A|U\to B|U$ is an $\epsilon$ isomorphism over $U$.\endproclaim
Previous versions of this lemma asserted that $f|U$ is an ``isomorphism over $U^{-3\epsilon}$''. This however depends on the
situation of $Y$ in $X$ while the lemma above uses $U-\t Fr ^{3\epsilon}U$, which is intrinsic to $U$. 

Let $f^-$, $h_A$, $h_B$ be data showing that $f$ is an $\epsilon$ isomorphism on $X$. We claim the restrictions to $U$ give data
showing $f|U$ is an $\epsilon$ isomorphism on $U$. Fix a basis element $x\in A$, then the first assertion is that
$f(x)=(f|U)(x)$ unless $x\in \t Fr ^{\epsilon}U$. If the two morphisms differ on $x$ then there is a path $\gamma$ of $f$ that
starts at
$x$ and intersects $X-U$. Let $t_0$ be the minimum of $\{t\mid \gamma(t)\notin U$. Since $U$ is open this point is in the
closure of $U$ but not in $U$ itself. It thus gives a point in the metric frontier. Since the path has radius  $<\epsilon$, $x$ is in
$\t Fr ^{\epsilon}U$.

Similarly $f^-|U$ equals $f^-$ except within $\epsilon$ of the frontier. Composing we see that $(f^-|U)(f|U)$ is equal to
$f^-f$ on $A|(U-\t Fr ^{2\epsilon}U)$. 

Finally the same argument shows that the homotopy $h_A|U$ is equal to $h_A$ on any path that does not intersect $\t Fr
^{\epsilon}U$. Therefore on $A|(U-\t Fr ^{3\epsilon}U)$ the morphisms, composition, and homotopy are all unchanged. By
hypothesis the original homotopy ends with the identity over $X-\t Fr ^{3\epsilon}X\supset (U-\t Fr ^{3\epsilon}U)$. This gives the
conclusion needed for the lemma.

\subhead 4.2 Triangularity\endsubhead

\subsubhead 4.2.1 Diagonal and triangular morphisms \endsubsubhead
The main modification in the uncontrolled definition is to relax conditions near the metric frontier.

Suppose $A, B$ are geometric modules over $E\to X$ and $U$ is a subgroup of the units of $R$ ( in this paper $\{1\}$ or $\{\pm 1\}$). 

A geometric morphism $d\:A\to B$  is $\epsilon, U$ {\it diagonal\/} over $X$ if there is at most one path beginning on each basis
element of
$A$, there is such a path beginning at $x$ if $x\notin \t Fr ^{\epsilon}X$, the endpoints of these paths are distinct,  the
coefficients are elements of the subgroup
$U$, and the paths all have radius $<\epsilon$.

Compositions of $\epsilon$ diagonal morphisms are $2\epsilon$ diagonal. More generally the composition of an $\epsilon$ and a $\delta$ diagonal morphism is $\epsilon+\delta$ diagonal.

A geometric morphism is $\epsilon$ {\it triangular\/} if it is an $\epsilon$ morphism and also triangular in the uncontrolled sense. It follows for instance that the composition of two $\epsilon$ triangular morphisms is $\epsilon$ triangular. 
 Additional leverage comes from hypotheses on the partial order rather than modification of ``triangularity.''

\subsubhead 4.2.2 $\epsilon$-bounded partial orders \endsubsubhead
Suppose $X$ is a metric space and $S\to X$ is a function. A partial order on $S$ is $\epsilon$ {\it bounded\/} if for each $s\in S$ there is $n$ so that any increasing chain starting with $s$ has length bounded by $n$ and its image in $X$  lies in the $\epsilon$ ball about the image of $s$. 

The following is the analog of Lemma 3.2.4:

\proclaim{Lemma} 
Suppose  $f\:A\to B$ is $\epsilon$ triangular with respect to an $\epsilon$ bounded partial order on the base of $B$. Then
\roster\item the diagonal part  of $f$ is well-defined independently of the choice of such a partial order;
\item $f$  has a  $3\epsilon$ triangular inverse;
\item if\/ $U\subset X$ is open then the restriction to $U$ is $\epsilon$ triangular over $U$;
\item if\/ $g\:B\to C$ is $\delta$ triangular with respect to the same partial order on $B$ then $gf$ is 
$\epsilon+\delta$ triangular;
\item there is a factorization $f=\cdots (1+\alpha_n)\cdots(1+\alpha_1)\t diag (f)$ with $\alpha_i$ increasing,
$\alpha_i\alpha_j=0$ if\/ $j\leq i$, and starting at any basis element the composition is finite with paths of length $<\epsilon$ . 
\endroster\endproclaim
The proof is a straightforward extension of the proof of 3.2.4. For instance if we write $d+u=(1+ud^{-1})d$ then the inverse of the first factor is $\sum_{i=0}^{\infty}(-ud^{-1})$ and the inverse for the composition is $d\sum_{i=0}^{\infty}(-ud^{-1})$. Note the sum has only finitely many nonzero terms because the partial order is locally bounded. The terms in the sum are arbitrarily long compositions but have paths of radius $<3 \epsilon$. This is  because they are compositions of paths of radius $<\epsilon$ joining points that lie within $\epsilon$ of their starting point, since the partial order is $\epsilon$ bounded. 

Similarly the factorization is obtained as in 3.2.4. Again we observe the partial order condition shows the factorization exists, but the construction does not depend on the partial order.

\subsubhead 4.2.3 Examples\endsubsubhead
These illustrate the force of the $\epsilon$ hypothesis on the partial order. 
Let $X$ be the real line, and $A_{i,j}$ the
geometric module generated by the integers
$i,\dots,j$. Define
$u\:A_{i,j}\to A_{i,j}$ by: if $i\leq k<j$ then $u(k)$ is the segment from $k$ to $k+1$ with coefficient 1, and $u(j)=0$. Then $\t
id +u$ is triangular with respect to the standard linear order. It has radius $<1+\delta$ for any $\delta>0$. It is an isomorphism.
However the inverse has radius $j-i$. 

By mapping this example into a space $X$ we get examples of $\epsilon$  morphisms that are isomorphisms and triangular in the
uncontrolled sense, but whose inverse  has  radius arbitrarily large up to the diameter of $X$. 

Setting $i=-\infty$, $j=\infty$ (so $A$ is generated by all the integers, and $u$ is the upward shift) we get $\epsilon$ 
morphisms triangular in the uncontrolled sense and for which $d$ is an isomorphism, but $d+u$ is not even  an uncontrolled
isomorphism.

\subsubhead 4.2.4 Cancellation  \endsubsubhead
Suppose $C$ is an $\epsilon$ chain complex and $\xi$ an $\epsilon$ chain contraction. As in  as in 3.1.1 we say
$\xi$  {\it cancels the complement\/} of a based subcomplex $\hat C$ if a based decomposition  $C\perp \hat C=
D\oplus
\bar D$  such that the $C\to C\perp \hat C$ component of the boundary homomorphism has the form
$$\pmatrix \hat \xi&u&v\\0&0 &\delta\\
0&0&0\endpmatrix\:\hat C\oplus D\oplus \bar D\longrightarrow \hat C\oplus D\oplus \bar D $$ with $\delta$ a
$\pm1$ $\epsilon$ triangular morphism.

\subhead 4.3 Definition of $K^{lf}_1(X;p,R,\epsilon)$\endsubhead
In this section we define a space of $\epsilon$ controlled complexes over $p\:E\to X$, essentially by adding ``$\epsilon$
geometric'' and a local finiteness hypothesis to
 everything in definition 3.4.1. 
Basic naturality properties (those not needing
spacial localization) are also given in this section: restrictions in 4.3.2, functional images in 4.3.3, and dependence on
the reference map in 4.3.4. 

\subsubhead 4.3.1 Simplices in $K^{lf}_1$\endsubsubhead
Fix  $\epsilon>0$. An $n$-simplex of $K^{lf}_1(X;R,\epsilon)$ consists of:
\roster\item $\epsilon$ chain complexes $C_i$ over $p$ for $0\leq i\leq n$ so that the  basis elements have no points of
accumulation in the metric completion;
\item contractions $\xi_i$ for $C_i$ so that $(\xi_i)^k$ has radius $<\epsilon$ for all $k$;
\item for  $i<j$ an $\epsilon$ chain map $c_{i,j}\:C_i\to C_j$;
\item there exist$\epsilon$ bounded partial orders on the bases of the $C_i$ in each degree so that 
\itemitem{i)} the chain maps $c_{i,j}$ are $\epsilon$ triangular isomorphisms onto their images in each degree, and if $i<j<k$ then the
underlying basis functions for $c_{i,k}$ and $c_{j,k}c_{i,j}$ are equal on elements not taken to 0 by either;
\itemitem{ii)} in each degree, basis elements in $C_j$ not in the image of $c_{i,j}$ and outside $\t Fr ^{\epsilon}X$ preceed those in
the image of
$c_{i,j}$; and
\itemitem{iii)} the
contraction $\xi_j$ cancels the complements of the images of $c_{i,j}$.
\endroster

As in 3.4 we define faces of $K_1$ simplices by omitting one of the complexes. The collection of all such simplices defines a
simplicial complex (or a
$\Delta$-set), and this (or its geometric realization) is defined to be $K^{lf}_1(X;p,R,\epsilon)$.

\subsubhead 4.3.2 Restrictions \endsubsubhead
Suppose $U\subset X$ is an open set. All the definitions have been arranged so restrictions of $\epsilon$ objects to $U$ give
$\epsilon$ objects on $U$. Therefore we have contrived:
\proclaim{Lemma} If\/ $U\subset X$ is open then restriction defines a natural map 
$$K^{lf}_1(X;p,R,\epsilon)@>>>K^{lf}_1(U;p,R,\epsilon)$$
\endproclaim
This is an important part of the structure used to identify the limit as homology.

\subsubhead 4.3.3 Functional images \endsubsubhead
Suppose 
$$\xymatrix{ E\ar[d]^p\ar[r]^{\hat f}& F\ar[d]^q\\
X\ar[r]^f& Y}$$
is a commutative diagram with $X$, $Y$ metric spaces. Applying $(f,\hat f)$ to geometric gadgets over $p$ gives geometric
gadgets over 
$q$. A size estimate implies it takes controlled gadgets to controlled gadgets:
\proclaim{Lemma}Suppose $f\:X\to Y$ is a proper map, extends to a proper map of completions,  and $d(x,x')<\delta$ implies
$d(f(x),f(y))<\epsilon$. Then composition induces a natural map 
$$K^{lf}_1(X;p,R,\delta)@>f_*>>K^{lf}_1(Y;q,R,\epsilon)$$
\endproclaim
The properness is needed to see that the frontier of $X$ maps to that of $Y$, and that images of bases do not acquire accumulation
points in $Y$.  Otherwise this should be clear.
\subsubhead 4.3.4 Dependence on the  data\endsubsubhead
Definitions of geometric algebra over $p\:E\to X$ do not use very much of $E$ and $X$. All the paths and homotopies
essentially lie in the 2-skeleton, so other data with  similar 2-skeleta should have the same algebra and
$K$-theory. The next definition and lemma make this precise.
\subsubhead Definition\endsubsubhead a commutative diagram 
$$\CD E@>>>F\\
@VVV@VVV\\
X@>>>Y\endCD$$
is $(\delta,1)$-connected  if for every relative 2-complex $(K,L)$ and commutative diagram 
$$\xymatrix{
L\ar[d]\ar[r]^{\subset}& K\ar[d]\ar@{.>}[dl]_{g}\\
E\ar[r]^{\hat f}&F
}$$
there is a map $g$ (dotted arrow) so that the upper triangle commutes, the lower triangle commutes up to $\delta$, and on
the 1-skeleton of $K$ it $\delta$ homotopy commutes ($\delta$ measured in $X$).
Previous versions \cite{Quinn 1, 2} did not include the 1-skeleton homotopy  condition and we had to reconstruct it
using a $(\delta,1)$-connectedness hypothesis on $X$ itself.  Putting the condition in this definition makes the following
easier and more general:
\proclaim{Lemma} 
Suppose $\epsilon, \delta>0$ and $\hat f\:p\to q$ is a $(\delta,1)$-connected map of spaces over an embedding $X\subset Y$.
Then there is a map
$f^*$ as shown that makes the diagram homotopy commute:
$$\xymatrix{K^{lf}(X;p,R,\epsilon)\ar[rr]^{f_*}\ar[d]^{\t relax }&&K^{lf}(Y;q,R,\epsilon)\ar[d]^{\t relax }\ar[dll]_{f^*}\\
K^{lf}(X;p,R,\epsilon+3\delta)\ar[rr]^{f_*}&&K^{lf}(Y;q,R,\epsilon+3\delta)
}$$
\endproclaim
We have assumed $X$ is a  subspace (i.e\. $X\to Y$ is an isometry) to avoid having to spell out metric conditions on the
map. Note that $\epsilon$ algebra only uses an $\epsilon$ dense 2-skeleton of $X$. In the limit as $\epsilon\to 0$ all of
$X$ gets used, but 1-connectedness is still sufficient at the $E$ level.
\demo{Proof} 
Begin with a geometric module $A$ over $q$. Apply 1-connectedness with $K$ the basis of $A$ and $L$ empty. This gives a
$\delta$ homotopy of the basis to the image of a factorization through $\hat f$. Define $f^*A$ to be a module over $p$
obtained this way, then the basis homotopy gives a basis-preserving isomorphism $\hat f_*f^*A\to A$ of radius $<\delta$. 

Now suppose $c\:A\to B$ is a geometric morphism of modules over $q$ with radius $<\epsilon$. Composing with the
basis-preserving isomorphisms above gives a morphism $\hat f_*f^*A\to \hat f_*f^*B$. Apply the 1-connectedness hypothesis
again, this time with $K$ the paths in this morphism and $L$ the bases of $f^*A$ and $f^*B$. Use the resulting paths in $E$
and the coefficients in $c$ to define $f^*c\:f^*A\to f^*B$. The size constraint on the lift shows these paths have radius
$<\epsilon+\delta$. The homotopy part of 1-connectedness gives
$\delta$ homotopies from the $\hat f$ images of these paths to the paths in the morphism, and we conclude that $\hat
f_*f^*c$ is
$\delta$ homotopic to $c$. Moreover since there is a one-to-one correspondence between paths in the two morphisms, $f^*c$ is
 triangular if  $c$ is. 

Finally we consider relations. We need that if the composition $A@>a>>B@>b>>C$ is $\epsilon$ homotopic to $c\:A\to C$ then
$(f^*b)(f^*a)$ is $\epsilon+3\delta$ homotopic to $f^*c$. First note we get an $\epsilon+2\delta$ homotopy of images
$( \hat f_*f^*b)( \hat f_*f^*a)$ and $ \hat f_*f^*c$ by composing the homotopies from these to $ba$ and $c$ with the given
homotopy from $ga$ to $c$. Apply 1-connectedness with $K$ the domain of these homotopies (a union of copies of $I\times I$)
and
$L$ the boundaries. We conclude there are homotopies in $E$ whose images in $F$ are within $\delta$ of the input data. The
distance estimate shows these homotopies have radius $<\epsilon+3\delta$, as required.  There is no conclusion about the
images of these lifts being homotopic to the input data, but we do not need it: relations require only  existence. 

The fact that homotopy relations lift shows that if $c$ is an $\epsilon$ isomorphism or triangular isomorphism then $f^*c$
is an $\epsilon+3\delta$ isomorphism or triangular isomorphism. similarly chain complexes lift to chain complexes, chain
maps to chain maps, etc. This shows $f^*$ defines $K^{lf}_1(X;q,r,\epsilon)\to K^{lf}_1(X;q,r,\epsilon+3\delta)$ as
required. The basis-preserving maps $\hat f_*f*A\to A$ constructed in the first paragraph give the homotopy making the
lower triangle in the diagram in the lemma commute. If $B$ is geometric over $p$ then a similar construction gives a
basis-preserving isomorphism $B\to f^*\hat f_*B$, and this gives a homotopy making the upper triangle commute.

This completes the proof of the lemma.
\enddemo
\subhead 4.4 The limit $K^{lf}_1(X;p,R)$\endsubhead
The main object of the paper is defined in 4.4.1 as the inverse limit of the $\epsilon$ spaces of the previous section.
Functoriality is described in 4.4.2.
\subsubhead 4.4.1 Definition of $K^{lf}_1$  \endsubsubhead
If $\epsilon>\delta>0$ then relaxing control gives an inclusion
$K^{lf}_1(X;p,R,\delta)\subset K^{lf}_1(X;p,R,\epsilon)$. This defines an inverse system as $\epsilon\to 0$, and as
indicated in 2.1.4 we define 
$$K^{lf}_1(X;p,R) = \t holim _{\epsilon\to0} K^{lf}_1(X;p,R,\epsilon).$$ 
Homotopy inverse limits (the ``holim'' in the definition) are discussed in 2.1.4. Using the path model for linearly ordered
homotopy limits we can describe these as maps $(0,1]\to K^{lf}_1(X;p,R,\gamma)$, some $\gamma$, that ``converge in $X$'' in
the limit
$t\to 0$. More explicitly this means a triangulation of $(0,1]$ (usually infinite near 0) and a map so that for every
$\epsilon>0$ there is $\delta>0$ so that on $(0,\delta]$ the image lies in $K^{lf}_1(X;p,R,\epsilon)$. 

The naturality of the $\epsilon$ versions passes to the limit:
\proclaim{4.4.2 Proposition} $K^{lf}_1(X;p,R)$ has natural restrictions to open sets, and is functorial with respect to
 morphisms $(p\:E\to X)\longrightarrow(q\:F\to Y)$ over proper uniformly continuous 
$X\to Y$. Finally if\/ $X=Y$ and  $E\to F$ 
$(\delta_i,1)$-connected over $X$ for all 
$\delta_i> 0$ then the induced map on $K^{lf}_1$ is a homotopy equivalence.
 \endproclaim
Note the uniform continuity condition on naturality is too restrictive for the characterization theorem. It is extended to
naturality over all proper maps in the next section.

\head 5. Spacial localization and the axioms\endhead
Spacial localization refers to factorization of controlled activity, usually homotopies, into pieces constant over various pieces of the control space. For instance the cancellation of inverses is a canonical nullhomotopy of the map that takes a complex to the complex plus its suspension. In 5.1  the controlled version is ``localized'' to a subset $Y$ in the sense that the homotopy is described as the composition of two homotopies, one that keeps things constant over $Y$ and cancels away from $Y$, and another that finished the cancellation over $Y$. This depends  on special properties of triangular isomorphisms so is less formal than the properties developed in \S4. Section 5.2 extends this to a ``skew inverse'' version: a $K_1$ morphism $B\to A$ determines  a cancellation of $A\oplus SB$. Again this can be localized to a subset. Finally these are used in \S5.3 to localize general homotopies.

The constructions of \S\S5.1--5.3 are used to verify the homology axioms of 6.1 for $K$-theory. The homotopy axiom is verified in 5.4 using the simplest case of cancellation. Localization of homotopies is used to verify the exactness axiom in 5.4 and independence of the metric in 5.5. The stability theorem in \S7 is also proved with these tools.

\subhead 5.1 Standard inverses\endsubhead
The cancellation of inverses is a canonical nullhomotopy  $H\:\t id \oplus S\sim 0$ as  maps from $K_1$ to itself, see \S3.5. The objective here is to show that the process preserves $\epsilon$ control in an appropriate sense, and to localize it: factor it into a part that cancels the sum away from $Y$ but doesn't change it over $Y$, and a part that finishes the cancellation over $Y$. 

\proclaim{5.1.1 Proposition} Suppose $X, p, R$ as usual, $X\supset Y$, and  $\epsilon>0$ is fixed. Then there are
\roster\item a canonical map $T\:K_1^{lf}(X;p,R,\epsilon)^{(1)}\to K_1^{lf}(X;p,R,7\epsilon)$, where the superscript $(1)$ indicates the first derived subdivision;
\item canonical homotopies $H_1\:\t id \oplus S \sim T$ and $H_2\: T\sim 0$; and
\item a canonical homotopy $G\:H_1H_2\sim H$, where $H$ is the standard cancellation (i.e\. $H_2$ for the case $Y$ empty).
\endroster
These satisfy:
\roster\item restricted to $Y$, $H_1$ is constant equal to the identity, $H_2=H$ and $G$ is constant;\item restricted to  $X-Y^{20\epsilon}$, $H_1=H$, $H_2$ is constant equal to 0, and $G$ is constant;
\item If\/ $Y^{50\epsilon}-Y^{-25\epsilon}\subset U$ and $U$ is open in $X$ then $T, H_i$ and $G$ are natural with respect to restriction to $(U,Y\cap U)$; and
\item on $K_1^{lf}(X;p,R,\delta)$ $H_i$ and $G$ have radius $<7\delta$.
\endroster
\endproclaim
The $Y$-dependent parts of the construction are actually determined by data over $Y^{20\epsilon}-Y$. A more generous region is used in (3) to avoid ambiguity with a detail of the definition of $K$, namely that structural hypotheses on $\gamma$ objects are allowed to fail within $3\gamma$ of the metric frontier of the space. This means for $\gamma=7\epsilon$ on $Y^{50\epsilon}-Y^{-25\epsilon}$ there may be some question about what happens outside $Y^{29\epsilon}-Y^{-4\epsilon}$, but this is safely away from $Y^{20\epsilon}-Y$.

\subsubhead 5.1.2 Notation and outline of the proof\endsubsubhead
Fix $X, p, R$ and shorten $K_1^{lf}(X;p,R,\epsilon)$ to $K_1(\epsilon)$. The constructions are canonical (after choice of $Y$ and $\epsilon$) so need only be described for a typical simplex of $K_1(\epsilon)$.
We write an $n$-simplex as 4-tuple $(C,c,\xi,c_{i,j})$, where $C$ is a sequence of graded modules $C_i$ for $0\leq i\leq n$, $c_i$ is a boundary morphism for $C_i$, $\xi_i$ is a contraction for the complex $(C_i,c_i)$, and $c_{i,j}\: C_i\to C_j$ are chain maps, all in the controlled geometric senses defined in Section 4. 

In these terms the starting point, $(\t id \oplus S)(C)$ is given by $$(C\oplus SC,\left( \smallmatrix c&0\\0&-c\endsmallmatrix\right), \left(\smallmatrix \xi&0\\0&-\xi\endsmallmatrix\right),\left(\smallmatrix c_{i,j}&0\\0& c_{i,j}\endsmallmatrix\right)).\tag{$*$}$$
 Recall that the negative signs in the $SC$ component of the boundary and contraction come from the sign conventions of 3.1.3. Namely, the suspension functor changes signs of maps by $(-1)^{\t degree }$, and the boundary and contraction have degrees $-1$ and $1$ respectively. There is no sign change in the chain map piece because chain maps have degree 0. 

 The construction proceeds in three steps, each of which produces $K_1$ morphisms. The final homotopies are obtained by subdividing morphisms as described in 3.4.4. The first two steps standardize structure of objects but do not change underlying graded modules. The third step truncates modules.  We show how to factor each step, and assemble the three steps  by arranging that factorizations of the pieces can be commuted. The first derived subdivision arises in factoring the truncation step.
 
\subsubhead 5.1.3 Morphism to cones\endsubsubhead
This is an elaboration of the fact that $\left(\smallmatrix 1&-\xi\\0&1\endsmallmatrix\right)$ gives a chain isomorphism from $C\oplus SC$ to the mapping cone of the identity. 

Conjugating the starting data $(*)$ with $\left(\smallmatrix 1&-\xi\\0&1\endsmallmatrix\right)$ gives a  morphism to the simplex with data 
$$(C\oplus SC,\left( \smallmatrix c&\xi c+c \xi\\0&-c\endsmallmatrix\right), \left(\smallmatrix \xi&2 \xi^2\\0&-\xi\endsmallmatrix\right), \left(\smallmatrix c_{i,j}&c_{i,j} \xi_i -\xi_j c_{i,j}\\0& c_{i,j}\endsmallmatrix\right)).$$
There is a homotopy of this to 
$$(C\oplus SC,\left( \smallmatrix c&1\\0&-c\endsmallmatrix\right), \left(\smallmatrix \xi&0\\0&-\xi\endsmallmatrix\right),\left(\smallmatrix c_{i,j}&0\\0& c_{i,j}\endsmallmatrix\right))\tag{$**$}$$
in the following sense:
\roster \item By the hypothesis that $\xi$ is a contraction there is an $\epsilon$ homotopy of geometric morphisms $\xi c+c \xi\sim 1$. 
\item The degree-2 homomorphism $\left( \smallmatrix 0&2\xi^3\\0&0\endsmallmatrix\right)$ is a chain homotopy from $\left(\smallmatrix \xi&2 \xi^2\\0&-\xi\endsmallmatrix\right)$ to $\left(\smallmatrix \xi&0\\0&-\xi\endsmallmatrix\right)$ (after homotopy of geometric morphisms), and
\item  $ \left(\smallmatrix 0&\xi_j c_{i,j} \xi_i\\0& 0\endsmallmatrix\right)$ is a chain homotopy from $\left(\smallmatrix c_{i,j}&\xi_j c_{i,j}+c_{i,j} \xi_i\\0& c_{i,j}\endsmallmatrix\right)$ to $\left(\smallmatrix c_{i,j}&0\\0& c_{i,j}\endsmallmatrix\right)$, again after homotopy of geometric morphisms. 
\endroster
If the original simplex has radius $\leq \delta$ then all of these homotopies have radius less than $5\delta$. The largest potential excursions come from using the contraction identity to simplify the result of applying the boundary to the $2\xi^3$ term in the homotopy of the contraction. 

We claim the morphism $\left(\smallmatrix 1&-\xi\\0&1\endsmallmatrix\right)$ gives a $K_1(7\epsilon)$ morphism from the starting data $(*)$ to $(**)$. The content of this assertion is that there are $7\epsilon$-bounded partial orders in which the chain maps and morphism are triangular, and complements of images cancel. These partial orders are obtained from partial orders assumed to exist in the original simplex by a modification of ``shuffling the image partial orders'', see the proof of 3.5.1. In detail, the basis of $C_i\oplus SCi$ in degree $k$ is the union $\t base (C_i^k)\cup \t base (C_i^{k-1})$. On each subset the partial order is the given one. Suppose $s\in \t base (C_i^k)$ and $t\in \t base (C_i^{k-1})$. If they are comparable then $s>t$. Define them to be comparable if 
\roster\item the distance between them is less than $5\epsilon$, and
\item for every $j<i$, if $t\in \t im (SC_j)$ then $s\in \t im (C_j)$.
\endroster
An increasing chain in this partial order consists of a chain in $SC$ followed by a chain in $C$. The sub-chains have size bounded by $\epsilon$ and the beginning of the second is within $5\epsilon$ of the end of the first. Thus the maximum size is $7\epsilon$. The total length is bounded because the lengths of the subchains are. This therefore defines a $7\epsilon$ bounded partial order.  $\left(\smallmatrix 1&-\xi\\0&1\endsmallmatrix\right)$ and the homotopies used to modify the structure are all triangular in this partial order because the maps on the individual summands are, they preserve images (see 3.3), and the inter-summand maps have size $<5\epsilon$.

We now localize this. Let $U_1$ be a subset of  $X$ (to be specified later in terms of $Y$). Denote by $\rho$ the projection of a geometric module to the submodule supported by $U_1$.  Denote projection to the complement of $U_1$ by $(1-\rho)$. This is a slight abuse of notation since $(1-\rho)$ technically defines a homomorphism homotopic to this projection.  Factor the basic morphism by 
$$\left(\smallmatrix 1&-\xi\\0&1\endsmallmatrix\right)=\left(\smallmatrix 1&-\xi(1-\rho) \\0&1\endsmallmatrix\right)\left(\smallmatrix 1&-\xi(\rho)\\0&1\endsmallmatrix\right).$$
We extend this to a factorization of the whole $K_1$ morphism. Conjugating by the first (right) factor takes the starting data to 
$$(C\oplus SC,\left( \smallmatrix c&\xi\rho c+c \xi\rho\\0&-c\endsmallmatrix\right), \left(\smallmatrix \xi\rho& \xi\rho\xi+\xi^2\rho\\0&-\xi\endsmallmatrix\right),\left(\smallmatrix c_{i,j}&c_{i,j} \xi_i\rho-\xi_j\rho c_{i,j}\\0& c_{i,j}\endsmallmatrix\right)).$$
Modify the contraction and chain maps by homotopies $\left( \smallmatrix 0&2(\xi)^3\rho\\0&0\endsmallmatrix\right)$ and $ \left(\smallmatrix 0&\xi_j c_{i,j} \xi_i\rho\\0& 0\endsmallmatrix\right)$ respectively.  Finally note that the $\epsilon$ homotopy of homomorphisms $\xi c+c\xi\sim 1$ provides a homotopy of $\xi\rho c+c\xi\rho$ to a morphism that is $1$ over $U_1^{-\epsilon}$ and is constant over $X-U^{\epsilon}$ where $\xi\rho c+c\xi\rho =0$. Putting these together gives a simplex structure on the modules $C\oplus SC$ that agrees with the cone structure over $U_1^{-5\epsilon}$ and the sum structure over $X-U_1^{5\epsilon}$. The homomorphism $\left(\smallmatrix 1&-\xi(\rho)\\0&1\endsmallmatrix\right)$ provides a $K_1$ morphism from the sum structure to this, and is the identity over $X-U_1^{5\epsilon}$. 

This defines the first factor. The second factor is obtained similarly: consider $\left(\smallmatrix 1&-\xi(1-\rho) \\0&1\endsmallmatrix\right)$ as a map from the object just constructed to the object with structure obtained by conjugating. Modify this by the rest of the chain homotopies, namely $\left( \smallmatrix 0&2(\xi)^3(1-\rho)\\0&0\endsmallmatrix\right)$ and $ \left(\smallmatrix 0&\xi_j c_{i,j} \xi_i(1-\rho)\\0& 0\endsmallmatrix\right)$. The result is homotopic to the cone structure. Note this construction gives the identity over $U_1^{-5\epsilon}$. 

\subsubhead 5.1.4 Morphism to trivial complexes\endsubsubhead
This step is an elaboration of the fact that $\left( \smallmatrix 1&0\\c&1\endsmallmatrix\right)$ is a chain isomorphism of complexes 
$$(C\oplus SC, \left( \smallmatrix c&1\\0&-c\endsmallmatrix\right))\to (C\oplus SC, \left( \smallmatrix 0&1\\0&0\endsmallmatrix\right)).$$
This is also used in the trivialization of cancellable complexes in 3.3.3. 

Conjugating the output of the first step, 
$(**)$,
by $\left( \smallmatrix 1&0\\ c&1\endsmallmatrix\right)$ gives 
$$(C\oplus SC,\left( \smallmatrix 0&1\\c^2&0\endsmallmatrix\right), \left(\smallmatrix \xi&0\\c\xi+\xi c&-\xi\endsmallmatrix\right),\left(\smallmatrix c_{i,j}&0\\ cc_{i,j}-c_{i,j}c& c_{i,j}\endsmallmatrix\right)).$$
This is homotopic to 
$$(C\oplus SC,\left( \smallmatrix 0&1\\0&0\endsmallmatrix\right), \left(\smallmatrix 0&0\\ 1&0\endsmallmatrix\right),\left(\smallmatrix d_{i,j}&0\\0& d_{i,j}\endsmallmatrix\right))\tag{$***$}$$
where $d_{i,j}$ is the diagonal part of $c_{i,j}$. To see this first note the chain and contraction identities and $c^2\sim 0$ give geometric homotopies to
$$(C\oplus SC,\left( \smallmatrix 0&1\\0&0\endsmallmatrix\right), \left(\smallmatrix \xi&0\\1&-\xi\endsmallmatrix\right),\left(\smallmatrix c_{i,j}&0\\ 0& c_{i,j}\endsmallmatrix\right)).$$
The degree 2 homomorphism $\left(\smallmatrix0&0\\\xi&0\endsmallmatrix\right)$ provides a homotopy from this contraction to 
 $\left(\smallmatrix 0&0\\1&0\endsmallmatrix\right)$. Finally if the diagonal-plus-increasing decomposition of $c_{i,j}$ is $d_{i,j}+u_{i,j}$ then $\left(\smallmatrix0&0\\-u_{i,j}&0\endsmallmatrix\right)$ gives a chain homotopy from the chain map to $\left(\smallmatrix d_{i,j}&0\\0& d_{i,j}\endsmallmatrix\right)$. 
 
This makes sense because the diagonal-plus-increasing decomposition is determined by the given data, and in particular does not depend on a particular partial order. 
 
 As in the previous step this defines a $K_1$ morphism of $K_1$ simplices. Again appropriate partial orders are obtained by shuffling image partial orders and making points comparable only if they are within $5\epsilon$, but this time $SC$ comes {\it after\/} $C$ in each sub-piece to reflect the fact that  $\left( \smallmatrix 1&0\\ c&1\endsmallmatrix\right)$ is {\it lower\/} triangular. 
 
This morphism is factored in the same way as in the previous step. Let $U_2$ be a subset of $X$ and let $\rho$ denote projection to $U_2$. Factor the morphism by 
$$\left( \smallmatrix 1&0\\ c&1\endsmallmatrix\right)=\left( \smallmatrix 1&0\\ c(1-\rho)&1\endsmallmatrix\right)\left( \smallmatrix 1&0\\c\rho&1\endsmallmatrix\right).$$
Write the homotopies of the contraction and chain maps as 
$$\align \left(\smallmatrix0&0\\\xi&0\endsmallmatrix\right)=&\left(\smallmatrix0&0\\\xi(1-\rho)&0\endsmallmatrix\right)+\left(\smallmatrix0&0\\\xi\rho&0\endsmallmatrix\right)\\
\left(\smallmatrix0&0\\ -u_{i,j}&0\endsmallmatrix\right)=&\left(\smallmatrix0&0\\ -u_{i,j}(1-\rho)&0\endsmallmatrix\right)+\left(\smallmatrix0&0\\ -u_{i,j}\rho&0\endsmallmatrix\right)\endalign$$
Define a $K_1$ morphism from the cone structure by conjugating by the right factor, changing the contraction by the right summand, and applying homomorphism homotopies over $U_2^{-3\epsilon}$. The image of this morphism has the cone structure over $X-U_2^{3\epsilon}$ and the trivial structure over $U_2^{3\epsilon}$. 

Finally, again as in the previous step, the remainder of the morphism and homotopy give the second factor in the localization.

\subsubhead 5.1.5 Factoring trivial complexes\endsubsubhead
The contraction of a trivial complex (i.e\. of the form $(***)$) clearly cancels the whole complex, so the inclusion of the 0 complex is a morphism. The final step is to factor this using the fact that trivial complexes  have lots of subcomplexes whose complements cancel. 

A single complex of the form $(C\oplus SC,\left( \smallmatrix 0&1\\0&0\endsmallmatrix\right), \left(\smallmatrix 0&0\\ 1&0\endsmallmatrix\right))$ is a large sum of 2-dimensional subcompexes each located over a point. Thus if  $U_3$ is a subset of $X$ the restriction to $U_3$  is a subcomplex of the same form and the contraction cancels the complement so the inclusion into the whole complex is a morphism. A complication arises with simplices of complexes. The chain maps in the form  $(***)$ are diagonal so they take one 2-dimensional summand to another, but they may have nonzero radius. Thus the image of something in $U_3$  may not lie in $U_3$. Restricting both the range and domain to $U_3$ may give a chain map that is not injective, so cannot be a $K_1$ morphism. We deal with this by subdividing: restrict both the range and domain of a chain map, but introduce a new vertex between them. The new vertex will be the subcomplex of the range obtained by adding the restriction of the range and the image of the restriction of the domain. Both vertex restrictions restrictions inject into this.

The first derived subdivision of a simplicial complex is the complex with $n$-simplices monotone sequences of length $n+1$ of subsimpices of a simplex of the original. Faces are defined by omission, and degeneracies by duplication of an element in the sequence. 

Regard an $n$-simplex $(C,c,\xi,c_{i,j})$ of $K_{1}(\epsilon)$ as the image of a map $\Delta^n\to K_{1}(\epsilon)$. We describe how to get a map on the first derived subdivision. A vertex of the subdivision is a subsimplex $\tau$ (monotone sequence of length 1). Take this to $(C_i,c_i,\xi_i)$, where $i$ is the largest vertex in $\tau$. Edges should go to chain maps. If $(\sigma, \tau)$ is a monotone sequence with  largest vertices $i,j$  then $i\leq j$. Define the chain map $c_{\sigma, \tau}$ to be $c_{i,j}$, where we understand this to be the identity if $i=j$. 

We now define the restriction of the subdivision of a ``trivial'' simplex. Take a vertex $\tau$ to the sum of the images in its largest vertex, of the restrictions to $U_3$ of all  vertex complexes. Edges go to the evident restrictions and inclusions. Since all the complexes are trivial, the contractions cancel complements of all images, and this defines a simplicial map of $(\Delta^n)^{(1)}\to K_1$. Further, inclusions define a $K_1$ morphism of this to the unrestricted subdivision. This is our factorization of the inclusion of 0 into the original trivial simplex. The inclusion is the identity over $U_3$ and 0 outside $U_3^{\gamma}$, where $\gamma$ is the radius of the trivial simplex. The trivial simplices produced by earlier constructions have radius $<4\epsilon$.

\subsubhead 5.1.6 Proof of Proposition 5.1.1\endsubsubhead
Denote the morphisms described above by
$$\xymatrix{(*)\ar[r]^F&(**)\ar[r]^G &(***)&0.\ar[l]_H}$$ 
Let $Y$ be the subset specified in the statement of the proposition, and let $U_1=Y^{5\epsilon}$. Then the construction of 5.1.3 factors $F$ as $F_2F_1$ with $F_1$ constant on $Y$ and $F_2$ constant over $X-Y^{10 \epsilon}$. Next let $U_2=Y^{13\epsilon}$. Then the construction of 5.1.4 factors $G$ as $G_2G_1$ with $G_1$ constant on $Y^{10\epsilon}$ and $G_2$ constant on $X-Y^{16\epsilon}$. Finally let $U_3=Y^{16\epsilon}$, then 5.1.5 factors $H$ as $H_1H_2$ with $H_1$ constant over $Y^{16\epsilon}$ and zero outside $Y^{20\epsilon}$. 

Next we observe that two morphisms commute if one is constant where the other is nontrivial. Thus we can rearrange the sequence of morphisms to 

$$\xymatrix{ (*)\ar[r]^{F_1}&\ar[r]^{G_1}&&T \ar[l]_{H_1}\ar[r]^{F_2}&\ar[r]^{G_2}&&0\ar[l]_{H_2}}$$
where $T_0$ is the result of the first three operations. 

All of these morphisms have size $< 7\epsilon$ (the limiting factor being the bound on the partial orders in the first step). The intermediate functions therefore define  maps into $K_1(7\epsilon)$. Subdividing the morphisms gives homotopies $H_1$, $H_2$ as required for the Proposition. Note that the difference between the composition (in the concatenation sense) of these homotopies and the standard nullhomotopy comes from rearranging commuting $K_1$ morphisms. ``Subdividing'' the resulting commutative diagrams of morphisms gives the canonical homotopy $G$ between these two homotopies.

The other assertions in the Proposition are supposed to be easily seen from the form of the
construction.

\subhead 5.2 Skew inverses\endsubhead
In this section we use a $K_1$-morphism $f\:B\to A$ to give a cancellation of $A\oplus SB$. The standard cancellation of inverses corresponds to the identity map $A=A$.  As in the standard case this cancellation can be localized.

\proclaim{5.2.1 Proposition} Suppose $X, p, R$ as usual, $X\supset Y$, and $ \epsilon$ is fixed. Suppose $M$ is a simplicial set, $A,B\:M\to K_1^{lf}(X;p,R,\epsilon)$ are  simplicial and $f\:B\to A$ is an $\epsilon$ $K_1$ morphism. Then there are
\roster\item a canonical map $T\:M^{(1)}\to K_1^{lf}(X;p,R,9\epsilon)$, where the superscript $(1)$ indicates the first derived subdivision;
\item canonical homotopies $H_1\:A \oplus SB \sim T$ and $H_2\: T\sim 0$; and
\item a canonical homotopy $G\:H_1H_2\sim H$, where $H$ is the standard cancellation (i.e\. $H_2$ for the case $Y$ empty).
\endroster
These satisfy:
\roster\item restricted to $Y$, $H_1$ is constant equal to the identity, $H_2=H$, and $G$ is constant;\item restricted to  $X-Y^{30\epsilon}$, $H_1=H$, $H_2$ is constant equal to 0, and $G$ is constant;
\item $T, H_i$ and $G$ are natural with respect to simplicial maps of $M$ and restriction to $U,U\cap Y$, provided $Y^{60\epsilon}-Y^{-40\epsilon}\subset U$;  and
\item if  on a subcomplex $W\subset M$, $A, B, f$ have radius $<\delta$, then $T$, $H_i$, $G$ have radius $<9\delta$. 
\endroster
\endproclaim
As in 5.1 a generous margin is used  in (3) to avoid conflicts with details of the definition. Conclusion (4) can be thought of as a naturality with respect to scale. The original $\epsilon$ is used to choose subsets used in the construction, but otherwise does not enter, so radius of the output depends only on radius of the input. In fact the $\delta$ in (4) does not even have to be measured with the same metric, so (4) gives a naturality with respect to metric.

Again we shorten the notation for $K$ space to $K_1(\epsilon)$, and restrict attention to a single simplex of $M$ because the construction is canonical. The image of the simplex under $A$ is denoted $(A_i,a_i,\alpha_i,a_{i,j})$ as usual, and simlarly for $B$. We begin by describing the full construction $A\oplus SB \sim 0$, then show how to localize it. Finally most of the construction (including estimates) follows the standard case, and we focus on new features.

\subsubhead 5.2.2 Improvement of the image of $f$\endsubsubhead
Since $f_i\:B_i\to A_i$ is a $K_1$ morphism the contractions in $A$ cancel the complement of the image. According to 3.3.4 there is a canonical triangular endomorphism $h_i$ of the graded module $A_i$ that is the identity on the image of $f$ and conjugating by $h_i$ changes the boundary maps so the result decomposes (as a chain complex) as the sum of the image of $f$ and its based complement.  Moreover the contractions can be changed by homotopy, without disturbing the cancellation properties, to also decompose as sums. 

After this modification the chain maps can be improved too. Decomposing $A_i$ and $A_j$ as image of $f$ plus based complement, the fact that $a_{i,j}$ carries the image of $f_i$ into the image of $f_j$ means it has the form $\left( \smallmatrix x&y\\0&z\endsmallmatrix\right)$. Since the boundary homomorphisms in the two complexes are diagonal it follows that $\left( \smallmatrix x&0\\0&z\endsmallmatrix\right)$ (i.e\. omit the upper right term) is also a chain map. The order hypotheses on cancellations implies the difference is increasing with respect  to any admissible partial order. Finally recall that diagrams of chain maps in $K_1$ simplices are only required to commute up to increasing homomorphisms. Putting these together shows that omitting the off-diagonal terms from the $a_{i,j}$ gives another $K_1$ simplex, and the identity maps on the $A_i$ from the original simplex to this one is a $K_1$ morphism. 

The first step in the main construction is the $K_1$ morphism described above, from $A$ to the modification in which the image of $B$ splits. We proceed assuming this splitting condition.

Define $p_i\:A_i\to B_i$ to be the based projection of $A_i$ to the image of $f$, composed with the inverse of $f$. Thus $pf$ is the identity of $B$, and $fp$ is the projection of $A$ to the image of $f$. After the modifications in the first step $p$ is a chain map, $fp\alpha=\alpha fp$, and $fp\, a_{i,j}= a_{i,j}fp$, where we have omitted evident subscripts on $f$ and $p$. 
\subsubhead 5.2.3 Morphism to  cones\endsubsubhead
The sum $A\oplus SB$ has  boundary maps, contractions and chain maps 
$$(A\oplus SB,\left( \smallmatrix a&0\\0&-b\endsmallmatrix\right), \left(\smallmatrix \alpha&0\\0&-\beta\endsmallmatrix\right),\left(\smallmatrix a_{i,j}&0\\0& b_{i,j}\endsmallmatrix\right)).\tag{$*$}$$
The graded endomorphism $\left( \smallmatrix 1&-f\beta\\0&1\endsmallmatrix\right)$ defines a chain map to the structure conjugated by the endomorphism,
$$(A\oplus SB,\left( \smallmatrix a&af\beta +f\beta b\\0&-b\endsmallmatrix\right), \left(\smallmatrix \alpha&\alpha f\beta+f\beta^2\\0&-\beta\endsmallmatrix\right),\left(\smallmatrix a_{i,j}&a_{i,j} f\beta+f\beta b_{i,j}\\0& b_{i,j}\endsmallmatrix\right)).$$
Use the identity $af\beta +f\beta b \sim f$ in the boundary homomorphism. The degree-2 map 
$$\pmatrix f(\beta-p\alpha f)\beta p&\alpha f\beta^2\\ \beta p&\beta^2\endpmatrix$$
gives a homotopy from the given contraction to $\left( \smallmatrix \alpha(1-fp)&0\\p&0\endsmallmatrix\right)$. This takes us to
$$(A\oplus SB,\left( \smallmatrix a&f\\0&-b\endsmallmatrix\right), \left(\smallmatrix \alpha(1-fp)&0\\p&0\endsmallmatrix\right),\left(\smallmatrix a_{i,j}&a_{i,j} f\beta+f\beta b_{i,j}\\0& b_{i,j}\endsmallmatrix\right)).\tag{$**$}$$
At this point the contraction already gives a cancellation of the entire complex. We elaborate on this to make it easier to localize.
\subsubhead 5.2.4 Morphism to  trivial complexes\endsubsubhead
Apply the endomorphism $\left( \smallmatrix 1&0\\bp&1\endsmallmatrix\right)$ to get a morphism to a new simplex,
$$(A\oplus SB,\left( \smallmatrix a(1-fp)&f\\0&0\endsmallmatrix\right), \left(\smallmatrix \alpha(1-fp)&0\\p&0\endsmallmatrix\right),\left(\smallmatrix *&*\\0& *\endsmallmatrix\right)).\tag{$***$}$$
The triviality of the lower-left entry in the chain map is a consequence of the first step, where it was arranged that $a_{i,j}f =fp a_{i,j}f$. We will not need to track the other entries explicitly.

The simplex $(***)$ is based isomorphic to one of the form
$$(D\oplus SD,\left( \smallmatrix 0&d\\0&0\endsmallmatrix\right), \left(\smallmatrix0&0\\d^{-1}&0\endsmallmatrix\right),\left(\smallmatrix *&0\\0& *\endsmallmatrix\right)).\tag{$4*$}$$
where $d$ is a triangular endomorphism of the based graded module $D$.
The $1-fp$ summand of $A$ has this structure by the first step. The $fp$ summand of $A$, plus $SB$, has this structure after the last step. Combining these gives the structure for the whole complex. 

We simplify this further. Conjugate by $\left( \smallmatrix d&0\\0&1\endsmallmatrix\right)$ to get a morphism to a complex of the form 
$$(D\oplus SD,\left( \smallmatrix 0&1\\0&0\endsmallmatrix\right), \left(\smallmatrix0&0\\1&0\endsmallmatrix\right),\left(\smallmatrix d_{i,j}&0\\0& d_{i,j}\endsmallmatrix\right)).\tag{$5*$}$$
The final modification is to replace the $d_{i,j}$ terms in the chain maps by their diagonal maps. This is still a simplex, and the conjugation map from $(4*)$ to this diagonalized version is still a $K_1$ morphism because morphisms are only required to commute modulo increasing homomorphisms. 

Finally we complete the process $A\oplus SB\sim 0$ using the $K_1$ morphism from 0 to the complex $(5*)$.
\subsubhead 5.2.5 Localization\endsubsubhead
As explained in 5.1.6 it is sufficient to show how to  factor the morphisms. Morphisms of the form $\left( \smallmatrix 1&*\\0&1\endsmallmatrix\right)$ or alternatively of the form $1+u$ with $u^2=0$, can be factored over $Y$ simply by decomposing $u$ as a sum of a piece 0 over $Y$ and one 0 over $X-Y$ (see 5.1.3). Similarly modifications of contractions by homotopy can be factored by decomposing the homotopy as a sum. 

The end of the argument is also easy to localize. The final complex is $(5*)$, after adjustment to have diagonal chain maps. This decomposes as a huge sum of 2-dimensional complexes, and the cancellation is localized after subdivision as in~5.1.5.
This explains how to localize all the steps except the morphism from $(***)$ to $(4*)$ defined by $\left( \smallmatrix d&0\\0&1\endsmallmatrix\right)$. This step does not arise in the standard-inverse cancellation because in that case $d=1$. In general we do this by factoring $d$ as $d_2d_1$ with $d_1=1$ over $Y$ and $d_2=1$ over $X-Y^{4\epsilon}$. It is sufficient to do this when $\t diag (d)=1$, so $d=1+u$ with $u$ increasing. Set $d_1=1+u|Y$, and $d_2=(1+u)(1+u|Y)^{-1}$. The inverse is $3\epsilon$ triangular, and over $Z=X-Y^{3\epsilon}$ the inverse of the restriction is the restriction of the inverse: $(1+u|Y)^{-1}|Z=(1+u)^{-1}|Z$. Thus the composition defining $d_2$ reduces to $1$ over a slightly smaller set. 

As in 5.1.6 factoring the morphisms into commuting pieces gives a homotopy of the standard (unlocalized) homotopy to the composition of two partial homotopies with the properties claimed in Proposition~5.2.1.

\subhead 5.3 Localizing homotopies\endsubhead
 A homotopy into a $K$ space is factored into two homotopies, one that changes things away from $Y$ and one that changes things near $Y$. The thing appearing halfway through (between the two homotopies) is a spliced version of the things at the beginning and end. 

\proclaim{5.3.1 Lemma} Suppose $M$ is a simplicial complex, $A,B\:M\to K^{lf}_1(X;p,R,\epsilon)$ are simplicial, and $H\:A\sim B$ is a homotopy. Then there are
\roster\item a subdivision $M'$ and simpicial map $C\:M'\to K_1(9\epsilon)$;
\item homotopies $H_1\: A \sim C$ whose restriction to $Y$ is constant, and $H_2\:C\sim B$ whose restriction to $X-Y^{-60\epsilon}$ is constant; and 
\item a homotopy $H_1 H_2\sim H$.
\endroster
The map and homotopies constructed have the property that if $W\subset M$ is a subcomplex and the restrictions of $A,B$ and $H$ to $W$ are in $K_1(\delta)$ then the output map and homotopies are in $K_1(9\delta)$.\endproclaim
This is proved in \S\S5.3.2--5.3.4 by describing the homotopy $A\sim B$ as obtained from a sequence of $K_1$ morphisms, then changing them to a ``zigzag'' pattern to which 5.2 applies.

The  final conclusion holds even if the $\delta$ size condition uses a different metric than the one used for the $\epsilon$ conditions; see the remark after 5.2. This is used to show metric independence in \S5.6, though it could be avoided with more elaborate estimates in a single metric. 

The final conclusion also shows that the construction passes to homotopy inverse limits, though it still depends on an initial choice of $\epsilon$. Explicitly, we can regard a map of $M$ into the homotopy inverse limit as a map $M\times [n,\infty)\to K_1(1/n)$ that for $k>n$ takes $M\times [k,\infty)$ into $K_1(1/k)$. Applying the Lemma to homotopies of such maps factors them into pieces that still have size going to 0 in the $[1/n,\infty)$ coordinate. This gives:
\proclaim{Corollary} Suppose $H$ is a homotopy between maps $A,B\:M\to K^{lf}_1(X;p,R)$, $Y\subset X$ and $\epsilon>0$. Then $H$ is homotopic to a composition of $H_1\:A\sim C$ and $H_2\:C\sim B$ with the restriction of $H_1$ to $Y$ constant, and the restriction of $H_2$ to $X-Y^{60 \epsilon}$ constant.\endproclaim

\subsubhead 5.3.2 Zigzag homotopies \endsubsubhead
Suppose $A\leftarrow B\rightarrow C$ are $K_1$ morphisms. The ``zigzag homotopy'' from $A$ to $C$ is obtained by composing
\roster\item the homotopy from $A$ to $A\oplus SB\oplus C$ obtained by adding $A$ to the skew cancellation homotopy of 5.2 for $SB\oplus C$ using the second morphism, and
\item the homotopy from $A\oplus SB\oplus C$ to $C$ obtained by adding $C$ to the skew cancellation homotopy for $A\oplus SB$.\endroster
 The name comes from the picture:
$$
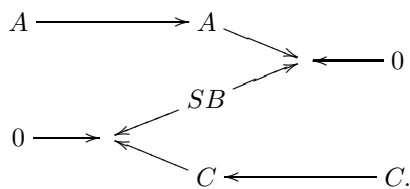
\xymatrix@R2pt{
A\ar[rr]&&A\ar[dr]\\
&&&&0\ar[l]\\
&&SB\ar[dl]\ar[ur]\\
0\ar[r]&\\
&&C\ar[ul]&&C.\ar[ll]
}$$
Variations on the following will be used in several places:
\proclaim{Lemma} Suppose $A\leftarrow B\rightarrow C$ are $K_1$ morphisms. Then the homotopy from $A$ to $C$ obtained by triangulating morphisms is homotopic to the zigzag homotopy described above. \endproclaim
\subsubhead Proof\endsubsubhead 
First, it is sufficient to show this for $A=B=C$. The diagram 
$$\xymatrix{B\ar[d]&B\ar[r]^=\ar[l]_=\ar[d]&B\ar[d]\\
A&B\ar[l]\ar[r]&C
}$$
gives a homotopy between the lower homotopy $A\sim C$ and the composition of the upper (identity) homotopy $B\sim B$ with the edge homotopies $A\sim B$ and $B\sim C$. Further this diagram
induces a diagram relating the morphisms used to define the zigzag homotopies. This similarly relates the two zigzag homotopies. Therefore a homotopy in the special case composes with these to give a homotopy in the general case.

Next we describe an intermediate step in the special case. Define the ``eye'' homotopy  beginning the same way as the zigzag: by adding $B$ to the cancellation homtopy $0\sim SB\oplus B$. But then finish by canceling the same pair rather than switching copies of $B$. The name comes from the picture:
$$\xymatrix@R2pt{
B\ar[rr]&&B&&B\ar[ll]\\
&&SB\ar[dl]\ar[dr]\\
0\ar[r]&&&&0\ar[l]\\
&&B\ar[ul]\ar[ur]
}$$
First note that this is homotopic to the triangulation homotopy (omit the eye): introducing and then removing the $SB\oplus B$ pair in exactly the same way is the composition of a homotopy $0\sim SB\oplus B$ and its reverse, so is homotopic to the constant homotopy. on the other hand there is a basis and structure preserving isomorphism between the eye and the zigzag given by the identity on the left half and the basis bijection that interchanges the $B$ terms on the right half. This isomorphism triangulates to give a homotopy between the zigzag and the eye. Combining these observations gives a homotopy from the zigzag to the constant homotopy. This completes the Lemma.

\subsubhead 5.3.3 The bounded case of Lemma 5.3.1\endsubsubhead
Suppose that there is a subdivision $M'$ and a triangulation $I'$ of $[0,1]$ so $H$ is simpicial with respect to the standard triangulation of the product $M'\times I'$. This is essentially a boundedness hypothesis on the number of vertices needed in the triangulation of $I$: in general for infinite $M$ there may not be such a bound, and we explain what to do about that in 5.3.4. 

When $H$ is simplicial on the product then restrictions to $M'\times\{v_i\}$, where $v_i$ are vertices of $I'$, gives a sequence of maps $M'\to K_1(\epsilon)$. The families of chain maps associated to vertices in $M'$ times edges in $I'$ fit together to give morphisms between these maps. Triangulations of these morphisms gives a simpicial map on the same triangulation as $H$. This is not quite $H$:  in some simplices some structure maps $a_{i,j}$ are replaced by compositions $a_{k,j}a_{i,k}$, but the identity on vertices defines a morphism between this new map and $H$.  They are therefore homotopic.

Denote the sequence of maps and morphisms by 
$$\xymatrix{A_1&B_1\ar[l]_{f_1}\ar[r]^{g_1}&A_2&B_2\ar[l]_{f_2}\ar[r]^{g_2}&{\cdots }\ar[r]&A_n&B_n\ar[l]_{f_n}\ar[r]^{g_n}&A_{n+1}}$$
where $A_1$ is the subdivision of $A$ and $A_{n+1}$ is the subdivision of $B$. 

The next step is a variation on the zigzag homotopy of 5.3.2. First add to $A_1$ the skew-inverse cancellations of pairs $SB_i\oplus A_{i+1}$ defined using the morphisms $g_{i}$. This gives a homotopy $A_1\sim A_1\oplus \Sigma_{i=1}^{n}(SB_i\oplus A_{i+1})$. 
Next reassociate the sum to $\Sigma_{i=1}^{n-1}(A_i\oplus SB_i)\,\oplus A_n$. Canceling the indicated pairs using the morphisms $f_i$ gives a homotopy 
$\Sigma_{i=1}^{n}(A_i\oplus SB_i)\,\oplus A_n\sim  A_n$. 
Joining these two homotopies gives a zigzag homotopy $A_1\sim A_n$. As in Lemma 5.3.2 this homotopy is homotopic to the original. It has a lot more kinks than 5.3.2, but it can be seen as a composition of a lot of simple zigzags or the proof of 5.3.2 can easily be generalized.

This zigzag homotopy can be localized to give the statement of Lemma 5.3.1. Specifically  the localized version of 5.2 splits the homotopy $A_1\sim A_1\oplus \Sigma_{i=1}^{n}(SB_i\oplus A_{i+1})$ into a homotopy constant over $Y$ and one constant over $X-Y^{-30\epsilon}$.  Similarly split the homotopy $\Sigma_{i=1}^{n}(A_i\oplus SB_i)\,\oplus A_n\sim  A_n$ into a piece constant over $Y^{30\epsilon}$ and one constant over $X-Y^{60 \epsilon}$. The middle pieces in the composition of these four homotopies commute (each is constant where the other is not). Interchanging them gives the desired factorization.

\subsubhead 5.3.4 The general case\endsubsubhead
By subdivision and straightforward manipulations we can arrange $M'$ to be a union of closed subcomplexes $N_i$ with 
\roster\item $N_i\cap N_j$ is empty unless $i,j$ are equal or differ by one;
\item $N_i\cap N_{i+1}$ is collared in each larger complex (e.g\. like the boundary of a regular neighborhood);
\item restricted to $N_i\times I$, $H$ is simplicial on the standard triangulation of $N$ times a triangulation $I_i$ of $I$, and the triangulation $I_{i+1}$ is a subdivision (i.e\. obtained by adding vertices) of $I_i$.
\endroster
As in 5.4.2 the vertices of $I_i$ give maps $N_i\to K_1(\epsilon)$, and the edges give morphisms between these maps. On $N_i\cap N_{i+1}$ the longer sequence is obtained by inserting identity morphisms into the shorter.

We can use the bounded-length construction to produce morphisms, homotopies etc\. on each $N_i$. These are patched together using the method of the  zigzag homotopy lemma. We outline the argument.

Construct an intermediate object by modifying the $i+1$ procedure:  introduce the alternating sum as before, but then reassociate differently. On the terms in the $i$ sequence associate as in the $i$ sequence. Leave the new (identity) terms alone. Now cancel the paired terms to get a sequence of morphisms. 

There is a basis and structure preserving isomorphism that matches up the associations used in the original $i+1$ and modified $i$ cancellations. This extends to a commutative diagram between the two output sequences of morphisms. This gives a homotopy between the $i+1$ sequence and the intermediate one.

On the other hand the intermediate sequence is the sum of the $i$ sequence and a lot of terms obtained by  introducing and then canceling in the same way, identity morphisms. On the homotopy level this corresponds to composing a homotopy with its inverse, so the result is canonically homotopic to the constant homotopy.  

The conclusion is that there is a canonical homotopy between the homotopies on $N_i\cap N_{i+1}$ produced by the $i$ and $i+1$ constructions. Define  maps on the union by using this homotopy in a collar on one side.

\subhead 5.4 The homotopy axiom\endsubhead
The homotopy axiom for controlled $K$-theory is the assertion that that the inclusion $X\times\{0\}\subset X\times I$ induces a homotopy equivalence of $K^{lf}_1$ spaces. We reduce this to the more primitive statement 5.4.1, then prove that using cancellations. The proof is similar to that of 5.3, but here we are producing homotopies rather than modifying them.

First note projection $X\times I\to X$ also defines a map of $K$ spaces, and the composition $X\times\{0\}\to X\times I\to X$ is the identity. Thus projection provides a homotopy inverse if we can show the other composition induces a map homotopic to the identity. The other map is the endomorphism of $K_1^{lf}(X\times I; p\times I, R)$ defined by projecting objects to $X\times \{0\}$. 

The second observation is that the $K$ spaces involved here are homotopy inverse limits  as $\epsilon \to 0$, of $\epsilon$ spaces (see 4.4). Thus it is sufficient to do an $\epsilon$ version for each $\epsilon$ and show they are compatible up to homotopy. The statement we need is:

\proclaim{5.4.1 Proposition} Suppose $X, p, R$ are as usual and $\epsilon>0$, then there is a standard homotopy between the relax-control inclusion
$$K^{lf}_1(X\times I;p\times I,R,\epsilon)\longrightarrow K^{lf}_1(X\times I;p\times I,R,8\epsilon)$$
 and the map induced by projection to $X\times \{0\}$, followed by inclusion. If\/ $\epsilon >\delta$ then there is a homotopy between the $\delta$ homotopy and the composition of $K_1(\delta)\to K_1(\epsilon)$ and the $\epsilon$ homotopy. \endproclaim

The proof uses only the unlocalized standard-inverse version of cancellation. The constructions are canonical so it is sufficient to describe the effect on a single simplex $(C,c,\xi,c_{i,j})$ of $K^{lf}_1(X\times I;p\times I,R,\epsilon)$.

Suppose $\epsilon =1/n$, and choose maps $p_i$, $0\leq i\leq n$ that push $X\times I$ toward the $0$ end, so that 
\roster\item $p_n$ is the identity;
\item $p_0$ is the projection to $X\times \{1\}$; and
\item the distance from $p_i$ to $p_{i+1}$ is less than $\epsilon$.
\endroster
For instance $p_i(x,t)=(x,i \epsilon t)$. In these terms the object is to define a homotopy from $p_n(C)$ to $p_0(C)$.

The first step is to use 5.1 to get a $7\epsilon$ homotopy 
$$p_n(C)\sim \big(\Sigma_{i=0}^{n-1} p_i(C)\oplus S\,p_i(C) \big)\;\oplus p_n(C).$$
Reassociate this to 
$$p_0(C)\oplus\big( \Sigma_{i=0}^{n-1}S\,p_i(C)\oplus p_{i+1}(C)\big).$$
Next use $d(p_i,p_{i+1})<\epsilon$ to change the $p_i(C)$ terms to $p_{i+1}(C)$ terms. More formally we can think of this as an $\epsilon$ morphism of $K_1$ objects, so it subdivides to give a homotopy. The final step uses 5.1 backwards to cancel the resulting pairs $S\,p_{i+1}(C)\oplus p_{i+1}(C)$, leaving $p_0(C)$. This gives a ``zigzag'' homotopy in the sense of 5.3 from $p_n(C)$ to $p_0(C)$.

Note that on $K_1(\delta)$ the cancellations gives $7\delta$ homotopies, but moving $p_i$ to $p_{i+1}$ is still an $\epsilon$ move. The size of the final homotopy is therefore limited by the number of intermediate pieces used, not the size of the input.

Now we must relate the homotopies obtained at different levels.
Suppose  $k\delta =\epsilon$. The prescription for $\delta$ calls for $k$ times as many shifted copies of $C$ as for $\epsilon$, and we can think of the $\epsilon$ process as involving the copies with index a multiple of $k$. This identifies the $\delta$ homotopy as obtained by putting extra zigzags in the $\epsilon$ homotopy, and therefore homotopic by Lemma 5.3.2. We go through  the argument to clarify where the size estimates enter.

  As an intermediate step in the comparison we introduce the additional copies in the $\epsilon$ process as well, but then cancel them back out without rearranging them. On the homotopy level this corresponds to composing a homotopy and its inverse, so this modified process gives a homotopy homotopic to the $\epsilon$ homotopy. 

The modified $\epsilon$ process starts out introducing the same $p(C)\oplus S\,p(C)$ pairs as the $\delta$ process, but they are grouped differently for the second cancellation. The differences separate into blocks of length $k+1$. Expicitly, in a sequence 
$$(p_{ki}(C)\oplus S\,p_{ki}(C))\oplus(p_{ki+1}(C)\oplus S\,p_{ki+1}(C))\cdots (p_{k(i+1)}(C)\oplus S\,p_{k(i+1)}(C))$$
 the $\delta$ process reassociates linearly, while the $\epsilon$ process associates $S\,p_{ki}(C)$ with $p_{k(i+1)}(C)$ and leaves the terms in between as they are. Note there is a permutation of the $S\,p_{*}(C)$ terms that takes the associate of each $p_{j}(C)$ in one association to its associate in the other. We can think of this as giving a basis-preserving endormorphism of the sum that preserves all $K_1$ structure. The cancellation construction is natural with respect to basis-preserving isomorphisms, so this isomorphism extends to a commutative diagram of morphisms between the morphisms used to define the cancellations. Triangulating this commutative diagram gives a homotopy between the homotopies. 
 
Finally since the permutation preserves the blocks of length $k+1$ and the shifted copies in such a block differ by less than $\epsilon$, it follows that the based isomorphisms and therefore the homotopies obtained from them do not increase radius more than $\epsilon$. This means the $\epsilon$ and $\delta$ homotopies are homotopic as maps into $K^{lf}_1(X\times I;p\times I, R,8\epsilon)$. This is the coherence claimed in 5.4.1, so the proof is complete.

\subhead 5.5 The exactness axiom\endsubhead
The exactness axiom 6.1(3) for $K$-theory is the assertion that given a morphism in the category, i.e\. a commutative diagram 
$$\xymatrix{F\ar[r]\ar[d]^q& E\ar[d]^q\\
Y\ar[r]^f&X
}$$
with $f$ proper, the induced sequence 
$$K^{lf}_1(Y\times\{1\};q,R)\to K^{lf}_1(Y\times I\cup_f X;q\times I\cup p,R)\to K^{lf}_1(Y\times I\cup_f X;q\times [0,1)\cup p,R)$$
is a homotopy fibration. The notation $Y\times I\cup_f X$ indicates the mapping cylinder, with the $Y\times 0$ end identified with its image in $X$. 

In 5.5.1 we describe a homotopy lifting property that implies a sequence is a homotopy fibration. Lemma 5.5.2 is a $K_1$ statement that follows from the homotopy localization lemma 5.3.1, and that together with the homotopy invariance proved in 5.4 shows that the lifting property is satisfied.

\subsubhead 5.5.1 A homotopy lifting property\endsubsubhead
A sequence 
$$\xymatrix{F\ar[r]^r&E\ar[r]^p& B}$$
with $r$ and inclusion and $pr=*$ is a  homotopy fibration if it has the following property: given $M\to E$ and a homotopy $pf\sim *$ then there are 
\roster\item a map $\hat f\:M\to F$;
\item a homotopy $\hat h\:f\sim r\hat f$; and 
\item a homotopy between $h$ and $p\hat h$
\endroster
such that if $W\subset M$ has $f(W)\subset r(F)$ and $h|W=*$ then $\hat f|W=f|W$ and the homotopy of (3) is constant on $W$. 

This is not a general definition, but it implies the general definition and fits the need. Convert $p$ into a Hurewicz fibration (e.g\. using the path space construction), then $pr=*$ defines a map from $F$ to the fiber. The lifting property applied with $M$ the homotopy fiber gives a map in the other direction. The special properties of the lift show it is a homotopy equivalence (in fact a deformation retraction). 

The next statement is an approximation to this lifting property for the spaces of the axiom.
\proclaim{5.5.2 Lemma} Suppose $f\:M\to K^{lf}_1(Y\times I\cup_f X;q\times I\cup p,R)$ is a map, and $h\:pf\sim 0$ is a homotopy, where $p$ is the restriction to $K^{lf}_1(Y\times I\cup_f X;q\times [0,1)\cup p,R)$. Then $h$ is homotopic by a homotopy constant over the complement of $Y\times [1/4,3/4]$ to a composition $\bar h_1\bar h_2$ with $\bar h_1$ constant over $Y\times [3/4,1]$ and $\bar h_2$ constant over $Y\times[0,1/4]\cup_fX$. \endproclaim
This follows from the homotopy localization lemma, more specifically from the Corollary of 5.3.1, where the $Y$ in the Corollary is replaced by $Y\times[3/4,1]$ and $\epsilon$ is chosen so $60\epsilon<1/2$.

This Lemma implies the lifting property as follows: let $\bar f$ be the result of the homotopy $\bar h_1$. Since $\bar h_2$ is constant over $Y\times[0,1/4]\cup_fX$ and the end result is 0, $\bar f$ is 0 over this space. Thus composing $\bar f$ with the projection to $Y\times \{1\}$ gives a map $\hat f\:M\to  K^{lf}_1(Y\times\{1\};q,R)$. The homotopy axiom gives a homotopy between the projection and the inclusion of $Y\times[0,1/4]\cup_fX$. Composing this homotopy with $\bar f$ gives a homotopy $\hat f\sim \bar f$. Composing this with the homotopy $\bar h_1$ gives a homotopy $\hat h\:f\sim \hat f$. Finally the homotopy $\bar h_2$ gives the homotopy-of-homotopies (3) in the lifting property. Notice that none of this changes something that starts out over $Y\times \{1\}$, so $f$, $h$ are unchanged on the subset $W\subset M$ that maps into $K^{lf}_1(Y\times\{1\};q,R)$. This verifies the final condition in the lifting property, and completes the proof that the sequence is a homotopy fibration.

\subhead 5.6 Metric invariance\endsubhead
The last ingredient needed for the axioms of \S6 is independence of the metric. General principles (4.4.2) only give functorality with respect to uniformly continuous maps of metric spaces, while the axioms require functoriality with respect to all proper maps.  We show:

\proclaim{5.6.1 Proposition} Suppose $X$ is locally compact and $d_1\leq d_2$ are metrics on $X$. Then the map 
$$K_1^{lf}((X,d_2);p,R)\to K_1^{lf}((X,d_1);p,R)$$
induced by the identity map is a homotopy equivalence.\endproclaim
Note the condition $d_1\leq d_2$ implies the identity is uniformly continuous, so it does induce a map (in fact an inclusion) of $K$ spaces. Also this is not a version of the stability theorem of \S7: here we get different control by rearranging the data at hand; there we get better control.

This relates to general functoriality as follows: suppose there is a morphism (including reference spaces) over a proper map $f\:X\to Y$ of metric spaces. Define a new metric on $X$ by $d_2(x,y)= \t max \bigl(d_X(x,y), d_Y(f(x),f(y))\bigr)$. The maps 
$$\xymatrix{(X,d_X)&(X,d_2)\ar[l]_{\t id }\ar[r]^f&(Y,d_Y)}$$
are then uniformly continuous and so induce maps on controlled $K$. According to the Proposition the first map induces a homotopy equivalence, and the ``metric independent'' induced map is supposed to be the composition of the second with a homotopy inverse. But since homotopy inverses are not well-defined these ``induced map'' are not well-defined, let alone functorial. 

The technical fix for this is to change the definition. Define ``$\epsilon$ control'' with respect to {\it functions\/} $\epsilon\:X\to (0,\infty)$ (see 2.1.3 and \cite{Quinn 2}), and define $K_1^{lf}$ to be the homotopy inverse limit over all such functions, not just the constants. This inverse system is much larger, so for instance the path model for the homotopy inverse limit does not apply. However the result is metric-independent and fully functorial. In this context the argument of Proposition 5.6.1 shows that the map of the constant-control version into the function-control version is a homotopy equivalence. This means the {\it only\/} need for the function-control version is to provide full functoriality: all of the real mathematics (axiom verifications, etc\.) can be done in the constant-control version. Accordingly we have focused on constant control in this paper, but for axiomatic purposes will assume full functoriality.

The following is a single-metric statement in which $d_2$ will be the ambient metric and constant $d_1$   control follows from non-constant $d_2$ control.
\proclaim{5.6.2 Lemma } Suppose $Y_i$ $i\geq 1$ are increasing compact subspaces of a metric space $X$ as above, $\delta_{j}$ is a positive decreasing sequence, and $Y_i^{60\delta_{i}}\subset Y_{i+1}$. Suppose $f\:M\times[m,\infty)\to K_1^{lf}(X;p,R,\delta_{1})$ is a map so that the restriction to $M\times [k,\infty)$ maps into $K_1(1/k)$. Then there is a homotopy $H\:f\sim \hat f$ in $K_1(9\delta_{1})$ so that 
\roster\item $H$ maps $M\times[k,\infty)\times I$ into $K_1(9/k)$;
\item images of $\hat f(M\times [k,\infty))$ have radius $<9\delta_{i+k}$ over $X-Y_i$; and 
\item if\/ $d_1$ is another metric on $X$ and  $W\subset M$ has images of $f(W\times[k,\infty))$ of $d_1$ radius $<\rho$  then images of $H(W\times[k,\infty)\times I)$ have $d_1$ radius $<9\rho$.
\endroster\endproclaim
\subsubhead 5.6.3 Proof of metric invariance from Lemma 5.6.2\endsubsubhead
Suppose $Y_i\subset X$ an increasing sequence of compact subspaces with  $X=\cup Y_i$. Then compactness implies that there are numbers $\delta_{i,n}>0$ so that for any $n$ and  $x,y\in X$, if $x, y\notin Y_i$ implies $d_2(x,y)<\delta_{i,n}$ (all $i$) then $d_1(x,y)<\frac1n$. Take $\delta_k =\t min \{\delta_{i,j}\mid i+j=k\}/9$, and reduce further if necessary so the conditions of the lemma are also satisfied: it is monotone decreasing  and $Y_i^{60\delta_{i}}\subset Y_{i+1}$. 

Use the path model for the homotopy inverse limit to identify $K_1^{lf}((X,d_2);p,R)$ as a space of maps $[1,\infty)\to K_1^{lf}((X,d_2);p,R,\epsilon)$ so that the radius of the image of $[k,\infty)$ goes to 0 as $k$ goes to $\infty$. Let $M$ be the subset that satisfies the estimate of the lemma: radius $<1/k$ on $M\times[k,\infty)$, and similarly $M_9$ is the subset with radius $<9/k$ on the same subset. It is easily seen that the inclusion of $M$ and $M_9$ into the whole path space are homotopy equivalences.  Let $W$ be the subspace of $M$ on which the $d_1$ radii also go to 0, and $W_9$ the corresponding subset of $M_9$. Then $W$, $W_9$ are homotopy equivalent to $K_1^{lf}((X,d_1);p,R)$. The objective is to show that the inclusion $W\subset M$ is a homotopy equivalence. 

Let $H$ be a homotopy as provided by the lemma, applied to the ``evaluation map'' $M\times[m,\infty)\to K_1^{lf}(X;p,R,\delta_{1})$. Conclusion (1) of the lemma shows this defines a homotopy in $M_9$ beginning with the inclusion $M\subset M_9$ and ending with a map $\hat f$. Conclusion (2) and the choice of $Y_i$, $\delta_{i,j}$ shows that $\hat f$ has image lying in $W$. Finally conclusion (3) shows that points in $W$ stay inside $W_9$ during the homotopy. These together with the facts that $M\subset M_9$ and $W\subset W_9$ are homotopy equivalences show that $\hat f$ gives a homotopy inverse for the inclusion.

\subsubhead 5.6.4 Proof of Lemma 5.6.2\endsubsubhead
The starting data is  $f\:M\times[m,\infty)\to K_1^{lf}(X;p,R,\delta_{1})$  so that the restriction to $M\times [k,\infty)$ maps into $K_1(1/k)$. Choose $n_i$ so that $1/n_i<\delta_i$. Denote by $f_i$ the restriction of $f$ to $M\times \{i\}$.

The objective is to define a homotopy $H\:M\times[m,\infty)\times [0,1]\to K_1(9\delta_1)$. We first describe $H$ on slices $M\times\{k\}\times[0,1]$, then describe how to fill in between these slices. 

Suppose $H$ is defined on $M\times\{k\}\times[0,1-\frac1r]$, and on $M\times\{k\}\times\{1-\frac1r\}$ the restriction over $X-Y_r$ is equal to $f_{n_{k+r}}$. $f$ on $M\times[n_{k+r},n_{k+r+1}]$ gives a homotopy of radius $<\frac1{n_{k+r}}<\delta_r$ from $f_{n_{k+r}}$ to $f_{n_{k+r+1}}$. Use 5.3 to factor this into a homotopy $G_1$ constant over $Y_r$ and $G_2$ constant over $X-Y_r^{60\delta_r}$. $G_1$ also extends by the constant homotopy over $Y_r$ to define a homotopy starting with $H|M\times\{k\}\times\{1-\frac1r\}$. Use this modified version of $G_1$ to define $H$ on $M\times\{k\}\times[1-\frac1r,1-\frac1{r+1}]$. Doing this for all $r$ defines $H$ on the open interval
$M\times\{k\}\times[0,1)$. This extends continuously to a map on the closed interval since it is eventually constant when restricted to any compact set. Specifically the restriction to $Y_r$ is constant (as a function of the last coordinate) on $M\times\{k\}\times[1-\frac1{r+1},1)]$. 

This defines $H$ on slices $M\times\{k\}\times[0,1]$. We check estimates for the final map, on  $M\times \{k\}\times\{1\}$. Over $Y_{r+1}-Y_{r}$ this is a ``splicing'' of $f_{n_{k+r}}$ and $f_{n_{k+r+1}}$, obtained as the result of the first piece of a factoring of the homotopy between them ($G_1$ above). It therefore has radius bounded by 9 times the radius of the homotopy, $\frac1{n_{k+r}}$ in this instance. Since these estimates decrease with $r$ we see that over $X-Y_r$ the radius is bounded by $\frac9{n_{k+r}}$ as required for the lemma.

With this description of $H$ on slices it should be clear how to fill in between them. We describe it explicitly at the $r^{th}$ stage of the construction, i.e\. how to define $H$ on $M\times [k,k+1]\times\{1-\frac1r\}$. Over $X-Y_r$ the $k$ end is $f_{n_{k+r}}$ and the $k+1$ end is $f_{n_{k+r+1}}$, and we  join them with $f|[n_{k+r},n_{k+r+1}]$. If $s\leq r$ then over $Y_s-Y_{s-1}$ the $k$ end is the result of the first factor $G_1$ in a localization factorization $G_1G_2$ of the homotopy $f|[n_{k+s-1},n_{k+s}]$. The $k+1$ end is the result of a similar factorization $G_1'G_2'$ of $f|[n_{k+s},n_{k+s+1}]$. We fill in between these with $G_2G_1'$. These descriptions over $Y_s-Y_{s-1}$ fit together to define a homotopy over all of $X$. This defines $H$ on intervals $M\times [k,k+1]\times\{1-\frac1r\}$, and therefore on boundaries of rectangles $M\times [k,k+1]\times[1-\frac1r,1-\frac1{r+1}]$. But the two paths on this boundary from the $k,1-\frac1r$ corner to the $k+1,1-\frac1{r+1}$ corner are the same, so the rectangle is easily filled in. 

This completes the construction of the homotopy $H$. The size estimate is verified during the proof. Other properties follow from properties of homotopy localization~5.3. 

\head 6. Homology \endhead
This section gives the Characterization Theorem 6.1 needed for the controlled assembly
isomorphism theorem 2.2.1. This is a  sharpened version of Theorem 8.5 of \cite{Quinn 2}, and  more
details are given.  Some of the terms used in the statement are defined and developed in 6.2. Section 6.3 combines
restriction and induced maps in a single functor structure. 6.4 derives the traditional form of excision. The spectrum
structure on $J$ is defined in 6.5, and spectral cosheaves constructed in 6.6. Homology with coefficients in these is
defined in 6.7, and assembly maps constructed in 6.8. The theorem is proved in 6.9.  Finally 6.11 provides an
``iterated homology identity'' that gives a Leray-Serre type spectral sequence for homology of the domain of a morphism.

\proclaim{ 6.1 Characterization theorem}Suppose $\Cal C$ is a category of maps over locally compact 
spaces and proper morphisms, and suppose $J\:\Cal C\to\text{\rm (pointed CW)}$ is a functor
satisfying the following axioms (see notes for explanations):
\roster
\item{\rm (Homotopy)} if\/  $(X,p)$ is an object then the inclusion-induced map 
$$J(X\times\{1\},p)\to J(X\times I,p\times\text{\rm
id})$$ is a homotopy equivalence;
\item{\rm (Restriction)} if\/  $(X,p)$ is a object and\/ $U\subset X$ is open then there is a natural
restriction map $J(X;p)\to J(U,p|U)$ that takes $J(X-U;p|(X-U))$ to the basepoint;
\item{\rm (Exactness)} if\/ $f\:(Y,q)\to (X,p)$ is a morphism in $\Cal C$ then 
$$J(Y\times \{1\},q)@>{\text{\rm incl}}>>J(Y\times [0,1]\cup_{f}X,q\times{\text{\rm id}}\cup p)@>{\text{\rm
rest}}>>J(Y\times[0,1)\cup_fX,q\times{\text{\rm id}}\cup p)$$
is a homotopy fibration sequence;
\item{\rm (Unions)} if\/ $X$ is a disjoint union of open sets $U_{\alpha}$ then the product of the restrictions
$$J(X,p)@>>>\prod_{\alpha}J(U_{\alpha},p|U_{\alpha})$$
is a homotopy equivalence.
\endroster
Then
$J$ has a natural (possibly non-connective) $\Omega$-spectrum structure denoted\/ ${\Bbb J}$. If\/ $p\:E\to X$ is an
object of\/ $\Cal C$ with $X$ a locally compact ANR and\/ $p$  a stratified system of fibrations then there is a
spectral cosheaf\/
$\Bbb J(p)\to X$ defined by applying
$\Bbb J$ ``fiberwise'', and the assembly map
$$\Bbb H^{lf}(X;\Bbb J(p))@>>>\Bbb J(X,p)$$
is an equivalence of spectra.
\endproclaim
\subsubhead 6.1.1 Notes\endsubsubhead\roster\item In the application  the category $\Cal C$ has objects $E\to X$ with
$X$ locally compact and metric. It may be different in other applications: all we need of it is that it is closed under
mapping cylinders, passing to open subsets, and taking products with finite CW complexes, and that it contains the maps
with $X$ a locally compact CW complex filtered by subcomplexes. 
\item The space
$Y\times[0,1]\cup_fX$ in (3) is the mapping cylinder of
$f$ defined by the relation $(y,0)\sim f(y)$ on the disjoint union. The maps are obtained by including and then deleting the outer end of the mapping cylinder. 
\item In (1) note that $X\times
I$ is the mapping cylinder of the identity map. Thus if (3) holds, (1) is equivalent to ``$J(X\times [0,1),p)$ is
contractible''. However (1) as stated seems to be more basic and slightly easier to verify in the application. 
\item The union
axiom (4) gives new information only in the infinite case, since finite unions follow from the other axioms. Note the
map is to the {\it product\/} with the product topology, so in the infinite case homotopy groups are products, not
sums, of homotopy of the factors.
\endroster

\subhead 6.2 Extended naturality\endsubhead
Restriction and naturality with respect to proper maps can conveniently be combined into a single structure. Let $\Cal
C^+$ denote the category of maps $p\:E\to X$ of pointed spaces with $X$ compact and  restriction to the
complement of the basepoint  gives an object in $\Cal C$. Denote the basepoint by ``$\infty$'', then explicitly
this means $p^{-1}(X-\infty)\to (X-\infty)$ is an object in $\Cal C$. In the pointed category ``mapping cylinder'' will
mean the ordinary mapping cylinder with the mapping cylinder over the basepoint identified to a point.
Approximate morphisms in this category are maps of reduced mapping cylinders that preserve basepoints, ranges, and
complements of ranges. 

1-point compactification gives a functor $\Cal C\to
\Cal C^+$, taking $E\to X$ to $E^+\to X^+$. As above $X^+$ is the 1-point compactification. $E^+$ is $E\cup\infty$ with
neighborhoods of $\infty$ given by inverse images of complements of compact sets in $X$. 
\proclaim{6.2.1 Lemma} the functor $J$ extends to $\Cal C^+$
\endproclaim
$J$ is defined on objects simply by applying the original functor to the complement of the basepoint. For morphisms
note that if $f\:X\to Y$ is a basepoint-preserving map of compact spaces then we get
$$X-\infty@<{\text{open subset}}<< f^{-1}(Y-\infty)@>{\text{proper}}>>Y-\infty$$
and similarly for reference maps over these. The inclusion and induced maps give 
$$J(X;p)@>>> J(f^{-1}Y;p|f^{-1}Y)@>>>J(Y;q)$$
and the composition of these is the desired induced map. 

\subsubhead 6.2.2 Caution\endsubsubhead One-point compactification does not make wild behavior of $X$ at
infinity go away, it just packs it into a neighborhood of the basepoint. Consequently the basepoint must be treated with
care. For example if we begin with
$X$ a locally compact CW complex with filtration by skeleta and 1-point compactify it, the result is almost never a CW
complex; the inclusion of the basepoint is not a cofibration; and filtering $X^+$ by adding the basepoint as a new
stratum does not give a filtration even dominated by a CW filtration. Compactification is a useful trick  in
simplifying naturality here, and necessary in defining locally finite  homology in 6.7.3, but does not avoid any of the
real work needed to deal with noncompact spaces.
\subhead 6.3 Excision\endsubhead
Excision relates locally finite homology of a closed subspace to homology of its complement. However this is only
expected for cofibered subspaces, even for
ordinary homology.
\subsubhead 6.3.1 Cofibrations\endsubsubhead
We recall that a subspace $Y\subset X$ is a {\it cofibration\/} if $Y\times I\cup X\times \{0\}$ is a retract of $X\times
I$. Explicitly this is a map $X\times I\to Y\times I\cup X\times \{0\}$ that is the identity on the subspace. To say it
gives a cofibration in a category of maps, as below, means there should be a retraction that lifts to a retraction
of of reference maps. This is a more categorical version of the older term
``$p$-neighborhood deformation retract'' \cite{Quinn 2}.

\proclaim{6.3.2 Lemma} If\/ $(Y,q)\to (X,p)$ is a cofibration in $\Cal C$ then the sequence
$$J(Y,q)@>>>J(X,p)@>>>J(X-Y,p)$$
is a homotopy fibration.\endproclaim
The first map in the sequence is induced by the inclusion, the second by restriction to the complement. The definition
of cofibration provides a morphism $X\times I\to X\times\{0\}\cup Y\times I$ extending the identity. Such a map can be
spread over two variables to give a homotopy of itself to the identity. Thus $X\times\{0\}\cup Y\times I$ is a
deformation retract of
$X\times I$ in the category. 

 Consider the diagram, in which reference maps are suppressed:
$$\xymatrix{J(Y\times\{1\})\ar[r]\ar@2{-}[d]&J(X\times\{1\})\ar[d]\ar[r]&J((X-Y)\times\{1\})\ar[d]\\
J(Y\times\{1\})\ar[r]\ar@2{-}[d]&J(X\times I)\ar[r]&J(X\times I-Y\times\{1\})\\
J(Y\times\{1\})\ar[r]&J(X\times\{0\}\cup Y\times I)\ar[u]\ar[r]&J(X\times\{0\}\cup Y\times[0,1))\ar[u]
}$$
The bottom row is a homotopy fibration according to the exactness axiom. The top row is the excision sequence. The
lemma will therefore follow if we show the vertical maps are equivalences. 

The left verticals are identities, so equivalences. The upper middle is an equivalence by the homotopy axiom. The lower
middle is an equivalence because it is induced by an inclusion that is a deformation retract in the category (the
cofibration hypothesis). Showing the right verticals are equivalences is simplified by using extended functoriality.
Taking 1-point compactifications and dividing out $Y\times \{1\}$ shows the lower right to induced by the inclusion 
$$(X^+\times\{0\}\cup Y^+\times I)/(\infty\times
I\cup Y\times \{1\})\subset  X^+\times I/(\infty\times I\cup
Y\times\{1\}).$$
But the deformation retraction of $X\times I$ passes to the quotient to give a deformation retraction for this
inclusion, so it induces an equivalence. Finally for the upper right, the standard deformation retraction of $X\times
I$ to $X\times \{1\}$ gives a deformation retraction of $X^+\times I/(\infty\times I\cup Y\times \{1\})$ to
$X^+\times \{1\}/(\infty\times I\cup  Y\times \{1\})$. Therefore this also is an equivalence.

\subhead 6.4 Spectrum structure\endsubhead
According to excision (6.3.2) the sequence
$$J(X\times\{0\},p)@>>>J(X\times[0,1),p)@>>>J(X\times (0,1))$$
is a homotopy fibration sequence. But it follows from the homotopy axiom that the middle space is contractible. The
homotopy fiber of the point map is the loop space, so we get a homotopy equivalence $\Omega(J(X\times R))\to J(X)$.
This equivalence is natural up to essentially canonical homotopy.

Recall  that an $\Omega$-spectrum is a sequence of spaces $J_n$ with homotopy equivalences $\Omega J_n\to
J_{n+1}$. We therefore make $J(X,p)$ into a spectrum with $J_n(X,p)=J(X\times R^n,p)$ if $n\geq 0$, and for negative
$n$ $J_n(X,p)=\Omega^{-n}J(X,p)$. Denote this spectrum by $\Bbb J(X,p)$. 
\proclaim{Lemma}This construction extends $J$ to a functor\/ $\Bbb J\:\Cal C^+\to \text{\rm
$\Omega$ spectra}$.\endproclaim There are minor homotopy issues
with the structure maps not quite being  canonical, but we will not go into details.

\subhead 6.5 Homotopy stratified maps\endsubhead
These are slightly more general than  ``stratified systems of
fibrations''. They occur frequently in applications and come up in the next section.
\subsubhead 6.5.1 Definition\endsubsubhead
A map $p\:E\to X$ is {\it homotopy stratified\/}  if there is a closed filtration $X^0\subset\dots\subset X^n\subset
\dots$ so that
\roster\item the restriction to each stratum, $p^{-1}(X^j-X^{j-1})\to (X^j-X^{j-1})$, is an approximate fibration, and
\item the inclusions $X^j\to X$ are approximate cofibrations  in the category of maps.
\endroster
 The following characterization works as a definition of approximate
fibration  if the reader does not recall the standard one:
\proclaim{6.5.2 Lemma} A map is homotopy stratified if and only if the stratified homotopy link in the
mapping cylinder is a stratified system of fibrations.\endproclaim 
Let $\t cyl (p)$ denote the mapping cylinder of $p\:E\to F$. The stratified homotopy link is the space of maps $[0,1]\to
\t cyl (p)$ so that $0$ maps to $X$, $(0,1]$ maps into the complement of $X$ in the inverse image of the stratum
containing the image of 0. The proof in the unstratified case is given in 
\cite{Quinn 5, \S2.7}, and the general case follows directly from this. For stratified homotopy links see also
\cite{Hughes, 1} and \cite{Quinn 8}. 

Evaluation gives a canonical ``approximate morphism" $\t holink (\t cyl (p))\times(0,1]\to E$ that is in the
approximate sense a fiber homotopy equivalence. In these terms an approximate fibration is canonically approximately
fiber homotopic to a genuine fibration, (the homotopy link) and we systematically use this to avoid working on the
homotopy level with approximate things.

\subhead 6.6 Functorial approximation\endsubhead
When working with {\it co\/}homology one constructs a {\it sheaf\/}  easily with no hypotheses: take an open
set $U$ to $J(U)$, and inclusions to restriction morphisms. The {\it co\/}sheaf construction involves applying $J$ to
point inverses rather than open sets. Locating appropriate ``point inverses'' seems to require the stratification
hypotheses and some work. In this section we show that the coefficient maps we work with can be approximated by
realizations of functors. This will make the passage to spectra easy: compose with the functor and realize. 
\S6.6.1 describes stratified maps obtained by realizing functors. 6.6.2 gives a functor version of pullbacks over
simplicial complexes, and 6.6.3 gives approximations of locally compact ANR filtered sets by simplicial complexes. 

\subsubhead 6.6.1 Realization of space-valued functors\endsubsubhead
Suppose $K$ is a simplicial complex, and consider $K$ as a category with objects the simplices and
 morphisms the inclusions $j_i\:\partial_i\sigma\to \sigma$. Realization gives CW complex $|K|$. 

 Realization of complexes extends in a standard way to functors. Let
$\hat E\:*K\to\text{spaces}$ be a contravariant functor. This means we get a space $\hat E(\sigma)$ for every simplex
$\sigma\in K$, and if $\tau\subset \sigma$ then a map $\hat E(\sigma)\to \hat E(\tau)$. The realization is:
$$|\hat E|=\coprod_n\bigl(\coprod_{\sigma^n}\hat E(\sigma)\times \Delta^n\bigr)/\sim$$
where  $\sim$ is given by: if $t\in \partial_i\Delta^n$ then  $(x,j_i(t))\in\hat E(\sigma)\times\Delta^n$  is
identified with $(\hat E(j_i)(x),t)\in\hat E(\partial_i\sigma)\times\Delta^{n-1}$. 

In these terms $|K|$ is the realization of the functor taking each simplex to a point. The natural transformation from
an arbitrary functor $\hat E$ to the point functor gives a map of realizations $|\hat E|\to |K|$. The following lemma
describes when this gives a homotopy stratified map:
\proclaim{Lemma} Suppose $K$ is a simplicial set and\/ $\hat E$ is a functor from $K$ to spaces. The map of
realizations $|\hat E|\to |K|$ is homotopy stratified over a filtration of\/ $|K|$ by subcomplexes if for any simplices
$\tau\subset \sigma$ both in the same stratum of $K$ the map $\hat E(\sigma)\to\hat E(\tau)$ is a homotopy
equivalence.\endproclaim
We often improve the output to be a stratified system of fibrations by applying stratified homotopy links, 6.5.2.

The point to be verified is that the restriction to a stratum is an approximate fibration. We can describe $|\hat E|$ as
built up of iterated mapping cylinders, and the projection to
$K$ is projection to iterated mapping cylinders of point maps. The structure on $K$ corresponds to describing a simplices
$\Delta^n$ as the mapping cylinder of a map $S^{n-1}\to \partial \Delta^n$ union the mapping cylinder of $S^{n-1}$ to a
point. The lemma now follows from the fact that iterated mapping cylinders of homotopy equivalences give approximate
fibrations, c.f\. \cite{Hatcher 1}.

The next result shows stratified systems of fibrations  over simplicial sets are equivalent in a suitable sense to
realizations.
\proclaim{6.6.2 Lemma (Functorialization)} Suppose $p\:E\to |K|$ is a  stratified system of fibrations over a
filtration of $|K|$ by subcomplexes. Then there is
 a functor $\hat E\:K\to \text{\rm spaces}$ so that the realization is homotopy stratified over the given
filtration and there is a fiber homotopy equivalence of stratified systems of fibrations 
$$\xymatrix{{\text{\rm holink}}(\t cyl (|\hat E|\to K),|K|)\ar[r]\ar[d]&E\ar[d]^p\\
K\ar@2{-}[r] &|K|.}$$
\endproclaim
Recall that  the homotopy link is the  stratified system of fibrations canonically associated to the realization, which
itself is only homotopy stratified.
\demo{Proof}
We want a contravariant functor, and to get contravariance we use open stars. The {\it open star\/} of a simplex $\tau$
is the union of interiors of all simplices that contain $\tau$. Inclusions give this the structure of a contravariant
functor of $\tau$: if $\tau\subset
\sigma$ then
$\t openstar (\tau)\supset\t openstar (\sigma)$. 
Define $\hat
E(\sigma)$  to be the restriction of $E$ to $\t openstar (\sigma)$. This inherits a contravariant functor
structure from the stars. The natural transformation to the point functor gives $|\hat E|\to |K|$ as above. The inverse
image of a point $x$ is the total space of the pullback of $E$ over the open star of the simplex $\tau$ that contains
$x$ in its interior.

There is a natural inclusion $E\subset |\hat E|$ defined by taking a point $e$ over $t\in \t int (\sigma)$ to
$(\{\sigma\},e,t)\in \{\sigma\}\times (E|\t openstar (\sigma))\times\Delta^n$. This is a fiber map in the sense that
it commutes with projections to $|K|$. Next we verify it is homotopy equivalences on  fibers. Over a point
$x$ in the interior of
$\tau$ the map of fibers is the inclusion
$E|x\subset  E|\t openstar (\tau)$. There is a (radial) deformation retraction of the open star to $x$
that preserves simplices, and therefore the stratification, until the last moment. Since $ E$ is a fibration over
strata this is covered by a deformation retraction of $ E|\t openstar (\tau)$ to $ E|x$. This is a homotopy
inverse for the inclusion, so the inclusion is a homotopy equivalence. 

Now restrict to a stratum in the filtration of $K$. Over this $E$ is a fibration. The equivalence of point
inverses just above and Lemma 6.6.1 show that $|\hat E|$ is an approximate fibration over this stratum. Thus the map
$E\to |\hat E|$ is an approximate fiber homotopy equivalence. Applying the homotopy link construction to the
approximate fibration turns it into a genuine fibration. A stratified homotopy equivalence of stratified systems of
fibrations has a fiber homotopy inverse, and this provides the map of the lemma.
\enddemo
The final result of the section provides a bridge to the previous lemma from the general case.
\proclaim{6.6.3 Lemma (Simplicial approximation)} Suppose $p\:E\to X$ is a  stratified system of fibrations over a
locally compact metric ANR. Then for every $\epsilon>0$  there is
 a locally compact polyhedron $K$, a stratified systems of fibrations $q\:F\to K$ with simplicial stratification and a
commutative diagram
$$\CD F@>>>E@>>>F\\
@VV{q}V@VV{p}V@VV{q}V\\
K@>f>>X@>g>>K\endCD$$
so that 
\roster\item $f$ is proper and stratum-preserving and the fiber map over it induces homotopy equivalence of fibers;
\item $g$ is filtration-preserving, and
\item the composition $gf$ taking $p$ to itself is filtration-preserving fiber homotopic to the identity, by a
homotopy of radius $<\epsilon$.
\endroster
\endproclaim
Note $g$ and the homotopy are not stratum-preserving. This can be arranged if $X$ is homotopy stratified, \cite{Quinn
5}. The construction is  sketched since it differs from previous versions only in that $X$ is not assumed finite
dimensional, and the only application of this is to shorten theorem statements by two words. 

 Choose neighborhoods of the spaces $X_n$ in the
filtration of
$X$ for which there are almost-strict deformation retractions to $X_n$ covered by homotopies of $p$, as specified in the
definition of stratified system of fibrations. Let $I^{\infty}$ denote  the (compact) Hilbert cube, and choose a proper
embedding $X\to [0,\infty)\times I^P{\infty}$.  Since $X$ is an ANR there is a neighborhood that retracts to $X$,
$r\:U\to X$, and this is proper. We can choose $U$ and neighborhoods $U_n$ of $X_n$ so that
\roster\item there is a non-decreasing sequence $k_j$ so that $U$ and each $U_n$ intersected with $[0,j]\times
I^P{\infty}$ is of the form $V\times I^{\infty-k_j}$ where $V$ is a closed codimension-0 submanifold of $[0,j]\times
I^{k_j}$;
\item $r$ takes $U_n$ into the neighborhood that deforms to $X_n$, and $U-U_n$ into $X-X_n$; and
\item all the deformations, etc\. are small compared to $\epsilon$.
\endroster
Define a stratified system of fibrations over $U$ with filtration $U_n$ by pulling back $p$. Applying the neighborhood
deformations in $X$ gives a small homotopy of $r$ to a stratum-preserving map covered by a map of systems of
fibrations. This is the map $f$ of the Lemma. The composition of $f$ with the inclusion of $X$ is no longer
the identity, but it has a small filtration-preserving homotopy to the identity. 

The last step is to omit $I$ factors from $U$ to get a locally finite complex. Choose a monotone nonnegative sequence
$i_j$ and consider the space $W=\cup_n (U\cap[j,j+1]\times I^{k_j+i_j})$. $W$ and each $W\cap U_n$ are unions of
compact smooth manifolds so can be triangulated, giving $W$ the structure of a locally compact polyhedron with PL
filtration. The restriction $r\:W\to X$ is still proper, stratum-preserving, etc. There is an essentially canonical
homotopy of $U$ into $W$ that is nearly the identity on $W$. In particular this homotopy takes the inclusion $X\to U$
to a map $q\:X\to W$.  
The retraction to $X$ and this homotopy combine to give a homotopy of the composition $rq\:X\to
X$ to the identity. Since the $X_n$ are contained in the interiors of the $U_n$, if the homotopy is small enough this
continues to be true during the homotopy. In other words $q$ and the homotopy are filtration-preserving. To complete the
construction we observe that by choosing the sequence $i_j$ large enough we can arrange the homotopy of $U$ into
$W$ to be arbitrarily small.

\subhead 6.7 Spectral cosheaves and homology\endsubhead
In this section we review spaces and spectra over a space $X$, and homology of these. The discussion is the same as
that in \S8 of \cite{Quinn 2} with a bit  more detail and basepoints. 
\subsubhead 6.7.1 Spectra over $X$\endsubsubhead
 A space {\it over\/} a space with basepoint $(X,*)$ is a map $p\:E\to X$ so that $p^{-1}(*)$ is a single point
(by definition the basepoint of $E$). A morphism over a map $X\to Y$ is a commutative diagram
$$\CD E@>>>F \\
@VVV@VVV\\
X@>>>Y\endCD$$
The map is required to preserve basepoints, but not complements of basepoints. 

The usual notions of topology extend to spaces over other spaces, the only modifications being occasionally dividing
out inverse images of basepoints \cite{James 1}. A {\it based space\/} over $X$ is a space over $X$ with a section:
$X@>s>> E@>p>> X$. If $(W,w_0)$ is a space with basepoint (not over anything) and $E$ is based over $X$ then the smash
product over
$X$,
$E\wedge_XW$, is
$E\times W/\sim$ where $\sim$ is the equivalence relation $(e,w_0)\sim (sp(e),w)$ for any $e\in E$, $w\in W$, and
$(*,w)\sim *$. This defines another based space over $X$. For instance $E\wedge_XS^1\to X$ has point inverses
suspensions of point inverses of $E\to X$. A {\it spectrum over\/} $X$ is a sequence of based spaces $E_n$ over $X$ and
based morphisms $E_n\wedge_XS^1\to E_{n+1}$ over the identity of $X$. 

The simplest example of a spectrum over $X$ is the product of $X$ and an ordinary spectrum. The ones we work with are
elaborations of this:
\proclaim{6.7.2 Lemma}Suppose $K$ is a simplicial set and\/ $\Bbb E\:K\to \text{spectra}$ is a contravariant functor.
Then realizations of the component spaces and maps of\/ $\Bbb E$ define a spectrum over $|K|$.\endproclaim
In each degree $n$ $\Bbb E$ has a space-with-basepoint valued functor $\Bbb E_n$ defined on $K$. Realizations of these
provide a sequence of spaces over $|K|$. Structure maps in $\Bbb E$ provide natural transformations $\Bbb E_n\wedge
S^1\to \Bbb E_{n+1}$. These induce morphisms of realizations $|\Bbb E_n|\wedge_{|K|}S^1\to |\Bbb E_{n+1}|$. These
objects and morphisms  constitute a spectrum over $|K|$.

\subsubhead 6.7.3 Definition\endsubsubhead
A spectrum $\Bbb E\to X$ is a {\it spectral cosheaf\/} over $X$ if there is a filtration so that each $E_n\to X$ is a
stratified system of fibrations over this filtration. 

Most of the cosheaves we work with come from functors and a sharpened form of  Lemma 6.7.2. In practice we mostly use
6.7.2, go directly to homology, and use properties of functors rather than properties of cosheaves. The exception is
the iterated homology identity in 6.10, where a spectral cosheaf is part of the input data.
 
\proclaim{6.7.4 Lemma}Suppose $K$ is a simplicial set with a filtration by subcomplexes, $\Bbb E\:K\to
\text{\rm spectra}$ is a contravariant functor, and if $\sigma\subset \tau$ are in the same stratum then $\Bbb E(\tau)\to
\Bbb E(\sigma)$ is an equivalence of spectra. Then realizations of the component spaces and maps of\/
$\Bbb E$ give an ``approximate'' spectral cosheaf over
$|K|$. Applying stratified homotopy links to mapping cylinders gives a spectral cosheaf.\endproclaim
Here ``approximate'' means the $E_n\to K$ are homotopy stratified (approximate fibrations over strata). This follows from
Lemma 6.6.1, and the improvement using homotopy links comes from 6.5.2.

\subsubhead 6.7.5 Homology\endsubsubhead
Suppose $X$ is a space with basepoint and $\Bbb E$ is a spectrum over $X$. The {homology\/} is the space defined
by the direct limit
$$H(X,*;\Bbb E)= \lim_{n\to\infty} \Omega^n(\Bbb E_n/X).$$
Some details: $\Bbb E_n/X$ is the space obtained by dividing out the image of the ``basepoint'' section $X\to \Bbb E_n$.
Dividing by the section in a smash product over $X$ gives the ordinary smash product, so the structure map of the
spectrum gives 
$$(\Bbb E_n/X)\wedge S^1@>>>\Bbb E_{n+1}/X$$
which is just the structure of an ordinary spectrum (i.e\. over a point). We now turn this into an $\Omega$-spectrum in
the stantard manner: take the adjoint of the structure map and apply $\Omega^n$ to get
$$\Omega^n(\Bbb E_n/X)@>>>\Omega^{n+1}(\Bbb E_{n+1}/X).$$
Finally take the direct limit of this system.

The homology space is in a natural way an $\Omega$-spectrum, denoted $\Bbb H(X;\Bbb E)$. Explicitly the $k^{\text{th}}$
space in the spectrum structure is obtained by changing the number of applications of $\Omega$ in the direct system:
$$\Bbb H(X,*;\Bbb E)_k= \lim_{n\to\infty} \Omega^{n-k}(\Bbb E_n/X).$$
The homology {\it groups\/} are the homotopy groups of this spectrum: $ H_j(X,*;\Bbb E)=\pi_j\Bbb H(X,*;\Bbb E)$. Note
that the homology and spectrum indexings go in opposite directions: $\pi_j(\Bbb H(X,*;\Bbb E)_k)=\pi_{j-k}(\Bbb
H(X,*;\Bbb E)_0)= H_{j-k}(X,*;\Bbb E)$.
\subsubhead 6.7.6 Locally finite homology\endsubsubhead
Locally finite homology is simply reduced homology of the 1-point compactification. If $\Bbb E\to X$ is a spectrum over
a locally compact space then define a spectrum $\Bbb E^+$ over $X^+$ simply by adding a point to each $\Bbb E_n$ over
$\infty$ and defining neighborhoods of the point to be inverse images of complements of compact sets in $X$.  We then
define
$$\Bbb H^{lf}(X;\Bbb E)=\Bbb H(X^+,\infty;\Bbb E^+).$$
See, however, the caution about 1-point compactifications in 6.3.2.

\subhead 6.8 Functors and assembly\endsubhead
Suppose $J$ is a functor satisfying the hypotheses of the Characterization theorem, and let  $\Bbb J\:\Cal C^+\to
\text{spectra}$ be the spectrum-valued enhancement described in 6.5. In this section we construct the assembly map used
in the theorem.
\subsubhead 6.8.1 $\Bbb J$ homology\endsubsubhead
First suppose $K$ is a simplicial complex and $\hat E\:K\to \text{spaces }$ is a contravariant functor. Think of a
space as a map from the space to a point, and compose this
$\Bbb J$ to get a functor $K\to \t spectra $. More explicitly this means the functor $\sigma\mapsto \Bbb
J(\t pt, \hat E(\sigma)\to \t pt )$. Realization defines a spectrum over $|K|$. Denote this $|K|$-spectrum by $\Bbb
J(\hat p)$, where $\hat p$ is the projection $|\hat E|\to |K|$, or by $\Bbb J(\hat p^{-1}(?))$, where $?$ indicates
this is to be thought of as a function of points in $|K|$. Neither notation is completely successful, so we  usually
spell out what they are supposed to mean in each context.

Homology spectra of these spectra are denoted $\Bbb H(|K|;\Bbb J(p))$. In the locally compact case $\Bbb
H^{lf}(|K|;\Bbb J(p))=\Bbb H(|K|^+,*;\Bbb J(p)^+)$. Finally, homology of general homotopy stratified $E\to X$ is defined
to be the direct limit  of homology of functor pullbacks over complexes (6.6.2), locally compact mapping properly to
$X$ in the locally finite case. 

\proclaim{6.8.2 Proposition} Suppose $\Bbb J$ is a homotopy-invariant functor from spaces to spectra. Then the locally
finite homology functor $J(p\:E\to X)=\Bbb H^{lf}(X;\Bbb J(p))$ satisfies the homology axioms 6.1(1)--(4). \endproclaim
Note we are not using the full structure of $J$ of 6.1, but just the ``trivial base'' case $E\to \t pt $ thought of as
a functor of the space $E$. 
 We will not give the proof here because it
is  long and technical and  differs from the standard development of generalized
homology
\cite{Whitehead 1} only in predictable ways. 
Perhaps the most subtle point is the union axiom (4) and the use of the full product topology.
This is where the topology of the 1-point compactification enters crucially. If something is an infinite union then in
the 1-point compactification most of the pieces are very close to the basepoint. This makes it possible for something
compact (e.g\. a sphere) to map nontrivially into all of them at once.

\subsubhead 6.8.3 Assembly\endsubsubhead
In this section we assume $\Bbb J$ satisfies all the hypotheses of the Characterization 
Theorem and construct the comparison map
used in the  theorem:
\proclaim{Lemma} There is a functorial map of spectra 
$$\Bbb H^{lf}(X;\Bbb J(p))@>A_X>> \Bbb J(X;p)$$
defined for $p\:E\to X$ homotopically stratified over a locally compact ANR, and $A_X$ is the identity when $X$ is a
point.\endproclaim The ``identity'' in the last part of the statement refers to the equivalence
$$\Bbb H(\t pt ,\Bbb J(E\to \t pt ))\simeq \Bbb J(E\to \t pt ).$$

We go through a series of reductions to get to  the core  construction. First it is sufficient to consider
the realization of a functor over a locally compact complex, $p\:|\hat E|\to |K|$ because the homology is the
direct limit of homologies of these.

Second, homology is defined using loop spaces of quotients $\Omega^n(\Bbb J_n(p)^+/|K|)$, so it is sufficient to define
maps on the quotients $\Bbb J_n(p)^+/|K|\to \Bbb J_n(|K|,p)$ and then take loops on these. Equivalently we define maps
 $\Bbb J_n(p)\to \Bbb J_n(|K|,p)$ that preserve ``basepoints'' in the sense that $|K|$ is taken to the basepoint of
$\Bbb J_n(|K|,p)$.

Finally, $\Bbb J_n(p)$ is defined by realizing a functor and the right side is natural so it is sufficient to construct
maps on pieces of the realization that fit together. For this we use the dual cones of the triangulation of $K$. Recall
the  {\it dual cone\/} $C_{\sigma}$ of a simplex $\sigma$ is the union of closed simplices in the first barycentric
subdivision that intersect $\sigma$ but not $\partial\sigma$. $C_{\sigma}$ intersects $\sigma$ in the barycenter; it is
the cone, with the barycenter as cone point, on the link of $\sigma$. Note $C_{\sigma}$ is contained in the open star
of
$\sigma$. To describe the data used in the construction of $A$ we need some notation.

Denote  the projection by $q\:\Bbb J_n(p)\to |K|$. Over a cone $C_{\sigma}$ 
$$q^{-1}(C_{\sigma})=\cup_{\alpha\in C_{\sigma}}\Bbb J_n(E|\text{openstar}(\bar \alpha))\times\alpha/\simeq.$$
The union is over simplices in the subdivision lying in the cone, and if $\alpha$ is such a simplex then  $\bar
\alpha$ denotes the smallest simplex in the original triangulation containing $\alpha$. We have also used the
abbreviated notation $\Bbb J_n(E|\text{openstar}(\tau))$ for $\Bbb J_n(E|\text{openstar}(\tau)\to\t pt )$.  Since all
the $\bar \alpha$ in the union contain $\sigma$ there are inclusions $\text{openstar}(\bar \alpha)\subset
\text{openstar}(\sigma)$. Projecting $q^{-1}(C_{\sigma})$ to the $\Bbb J_n$ factors and composing with the inclusion
defines a map $q^{-1}(C_{\sigma})\to \Bbb J_n(E|\text{openstar}(\sigma))$.

Denote the projection by  $\hat p\:|E|\to |K|$. 
The inverse image of a cone is a quotient of a union of products $E|\text{openstar}(\bar \alpha)\times\alpha$ as above.
Projecting to the $E$ factor gives a commutative diagram 
$$\CD \hat p^{-1}(C_{\sigma})@>>>E|\text{openstar}(\sigma)\\
@VVV@VVV\\
C_{\sigma}@>>>\t pt \endCD$$
Applying $\Bbb J_n$ gives a morphism $\Bbb J_n(\hat p^{-1}(C_\sigma)\to C_{\sigma})\longrightarrow \Bbb
J_n(E|\text{openstar}(\sigma))$.

 In these terms the fragments of the assembly map are 
 maps
$A_{\sigma}$ and  homotopies giving a homotopy-commutative diagram 
$${\xymatrix{ q^{-1}(C_{\sigma}) \ar[rr]^{A_{\sigma}}\ar[dr]&& {\Bbb J}_n(\hat p^{-1}(C_{\sigma})\to
C_{\sigma})\ar[dl]\\
&{\Bbb J}_n(E|{\text{openstar}}(\sigma))&}}$$
and these should be natural with respect to inclusions of cones. 

These maps are constructed by induction on the dimension of the cone (so downward with respect to dimension of
simplices, starting with maximal simplices). 0-dimensional cones have $A$ the identity. Suppose these are defined for
cones of dimension $<k$ and let $\sigma$ be a simplex with $C_{\sigma}$ of dimension $k$. The data for the cones in
$\partial C_{\sigma}$ fit together to give a map $A$ and a homotopy of the compositions into the lower right of the
diagram
$${\xymatrix{ q^{-1}(\partial C_{\sigma}) \ar[r]^{A}\ar[d]& {\Bbb J}_n(\hat
p^{-1}(\partial C_{\sigma})\to
\partial C_{\sigma})\ar[d]&\\
{\Bbb J}_n(E|{\text{openstar}}(\sigma))\ar[r]&{\Bbb J}_n(\hat
p^{-1}( C_{\sigma})\to
C_{\sigma})\ar[r]&{\Bbb J}_n(E|{\text{openstar}}(\sigma))}}$$
On the bottom row the first map includes $E|{\text{openstar}}(\sigma)\to \t pt $ over the cone point, and the second
map is induced by the projection of $\hat p^{-1}(C_{\sigma})$ to $E|{\text{openstar}}(\sigma)$. Both of these lower
maps are homotopy equivalences by the homotopy axiom and the fact that the radial deformation retraction of the cone to
the cone point is covered by a deformation of $\hat p^{-1}(C_{\sigma})$. Thus the homotopy of compositions into the
lower right space lifts to a homotopy of compostions into ${\Bbb J}_n(\hat
p^{-1}( C_{\sigma})\to
C_{\sigma})$. Now observe that $q^{-1}(C_{\sigma})$ is the mapping cylinder of the left vertical map. The lifted
homotopy defines a map on this mapping cylinder. This map is $A_{\sigma}$, and the homotopy in the lower right space
gives the homotopy needed to complete the data. 

In this fashion $A$ can be defined on all the $k$-dimensional cones, and by induction on the whole spectral functor
realization. Note this $A$ is not completely canonical: there are choices in lifting the homotopy in the last step.
However any two lifts are homotopic, any two homotopies between lifts are themselves homotopic, etc. This means 
for any practical purpose it can be considered natural and canonical.

\subhead 6.9 Proof of the Characterization Theorem\endsubhead
This follows the traditional (Eilenberg-Steenrod) proof of uniqueness for ordinary homology, working up from points
using excision and homotopy invariance.

\subsubhead Step 1\endsubsubhead Theorem 6.1 holds for projections $F\times R^n\to R^n$.

Let $R^n_+$ denote the upper half space and $p$ the projection $F\times R^n_+\to R^n_+$. The inclusion $R^{n-1}\subset
R^n_+$ gives excision sequences and a homotopy commutative diagram
$$\CD \Bbb H(R^{n-1},\Bbb J(p))@>>> \Bbb H(R^n_+,\Bbb J(p))@>>>\Bbb H(R^n_+-R^{n-1},\Bbb J(p))\\
@VV{A}V@VV{A}V@VV{A}V\\
\Bbb J(R^{n-1};p)@>>> \Bbb J(R^n_+;p)@>>>\Bbb J(R^n_+-R^{n-1};p)\endCD
$$
By excision the rows are homotopy fibrations. Thus if any two of the vertical maps are equivalences the third is also.
The one-point compactification of the half space is contractible, so by homotopy invariance both
functors applied to it give contractible spectra. Thus the center vertical map is trivially an equivalence. The space
on the right, $R^n_+-R^{n-1}$, is homeomorphic to $R^n$.  By
construction the assembly map is an equivalence for a projection to a point ($=R^0$), and this is the left vertical
when $n=1$. The fibration property implies the right vertical, over $R^1$, is an equivalence. This argument shows
inductively that assembly is an equivalence for all $n$.

\subsubhead Step 2\endsubsubhead Theorem 6.1 holds for  $p$ homotopy stratified over a finite complex filtered
by subcomplexes.

Suppose $K$ is such a finite complex and let $K_i$ denote the
$i$-skeleton. Then the inclusions $K_{i}\subset K_{i+1}$ give excision fibration sequences
$$\Bbb H(K_i;\Bbb(p))@>>>\Bbb H(K_{i+1};\Bbb(p))@>>>\Bbb H^{lf}(K_{i+1}-K_i;\Bbb(p))$$
and similarly for $\Bbb J$.  $K_{i+1}-K_i$ is a (finite) union of copies of $R^{i+1}$
over which $p$ is constant, so Step 1 shows the assembly map is  an equivalence for
these terms. This applies to the right-hand term in all the fibration sequences, and the left-hand term when $i=0$. We
can now proceed by induction as above. In a finite number of steps $i$ reaches the dimension of the complex
 so we conclude 6.1 holds for $K$. 
\subsubhead Step 3\endsubsubhead Theorem 6.1 holds for locally finite complexes and $p$ homotopy stratified over a
 filteration by subcomplexes.

For this we decompose a locally
finite 
$|K|$ in a nice way.  Define the {\it closed star\/} of $V\subset |K|$ to be the union of all closed simplices that
intersect
$V$. Note that if $V$ is closed then it is contained in the interior of its closed star. 

Let $V_0$ consist of a vertex in each component of $|K|$, and inductively define $V_n$ to be the closed star of
$V_{n-1}$. Define $A$ to be the union over $n$ of the closures of $V_{2n}-V_{2n-1}$. Define $B$ similarly as the union
of closures of $V_{2n+1}-V_{2n}$. Note each of these pieces is a finite complex, and the unions are disjoint unions.
We get:
\roster\item $|K|=A\cup B$;
\item the components of $A$ and $B$ are finite subcomplexes of $|K|$
\endroster
Since the theorem is true for finite complexes  it is true for disjoint unions of finite complexes by the union axiom
(for both $J$ and homology). Thus it is true for $A$, $B$, and $A\cap B$. Use excision twice, first with $A\cap
B\subset B$ to conclude the theorem holds for $B-A\cap B= A\cup B-A$, and then with $A\subset A\cup B$ to conclude it
holds for $A\cup B$. 
\subsubhead Step 4\endsubsubhead Theorem 6.1 holds in general.

Suppose $p\:E\to X$ is homotopy stratified over the locally compact ANR $X$. According to the simplicial
approximation and functor pullback lemmas 6.6.3 and 6.6.2 $p$ is a proper  retract of the realization
of a functor over the realization of a (locally finite) simplicial complex mapping properly to $X$. Since retracts of
isomorphisms are isomorphisms the ANR case follows from the locally finite complex case.

This completes the proof of the Characterization Theorem.

\subhead 6.10 The iterated homology identity\endsubhead
This is an equivalence of spectra so that application of the Atiyah-Hirzebruch type spectral sequence to its homotopy
groups gives a Leray-Serre type spectral sequence for homology of a map. The compact-support version is stated in
6.10.1, the stratification hypothesis and its locally compact analog are discussed in 6.10.2, and the locally compact
version of the theorem is given as 6.10.3. The remainder of the section contains proofs.
\proclaim{6.10.1 Theorem}{\rm (Iterated homology identity)} Suppose $\Bbb E@>p>> X@>f>>Y$ satisfy:
\roster\item $X$ and $Y$ are filtered ANRs;
\item $f$ is an ``approximate stratified system of stratified fibrations'' (see below) over the filtration of $Y$;
\item $p$ is a spectral cosheaf over the filtration of $X$.
\endroster
Then applying $\Bbb H(-;\Bbb E)$ to point inverses in $X$ gives a spectral cosheaf $\Bbb H(f;\Bbb E))\to Y$, and there is
a natural equivalence of spectra
$$\Bbb H(Y,\Bbb H(f;\Bbb E)))@>\simeq >> \Bbb H(X;\Bbb E).$$
\endproclaim
We discuss the stratification hypothesis.
\subsubhead 6.10.2 Definition\endsubsubhead
A  {\it stratified fibration\/} is a map $U\to V$ with a filtration on the total space that satisfies a
stratum-preserving version of the lifting property for fibrations, see \cite{Hughes 1, 2}. It follows that on each
stratum the restriction
$(U_i-U_{i-1})\to V$ is a fibration, and in the case at hand this is nearly sufficient. We now extend this over strata
in the base: $f\:X\to Y$ is a {\it stratified system\/} of stratified fibrations if the restriction of the filtration of
$X$ to the inverse of each stratum in $Y$, $f^{-1}(Y_i-Y_{i-1})\to (Y_i-Y_{i-1})$ is a stratified fibration, and $f$
satisfies a stratified version of the cofibration condition for the filtration of $Y$. 

This is a lot of data. However it often will come naturally from the setting or can be arranged by general principles.
Examples are:
\roster\item If $X$ and $Y$ are homotopy stratified sets then the projection $X\times Y\to Y$ is a stratified system of
stratified fibrations.
\item Suppose a finite group $G$ acts in a homotopy stratified way on a manifold $M$, suppose $H\subset G$ is a 
subgroup, and let $M^{(H)}$ denote the set of points fixed by some conjugate of $H$. $M^{(H)}/G$ is a homotopy
stratified space with the orbit type filtration \cite{Quinn 5}. The other strata of $M/G$ generally will contain
$M^{(H)}/G$ in their closure. The stratified homotopy link \cite{Hughes 1, 2} provides a model for these strata in a
neighborhood. It is a stratified system of stratified fibrations over $M^{(H)}/G$.
\item More generally if $Y\subset X$ is a pure subset (closed union of components of strata) in a homotopy stratified
space then the stratified homotopy link is a stratified system of stratified fibrations over $Y$.
\item If $K\to L$ is a proper PL map of polyhedra then there is a subdivison and a filtration of $K$  by subcomplexes so
that $K\to L$ becomes a stratified system of stratified fibrations over the filtration of $L$ by skeleta.
\item Generic smooth proper maps are stratified systems of fibrations with filtrations defined by stratifications of jet
bundles and coincidence conditions.
\endroster

We note that the map $f$ in the theorem is only assumed to be an {\it approximate\/} stratified system of stratified
fibrations. This is important in geometric applications because the ones encountered are almost never genuine. The
approximate versions of these definitions are spelled out in \cite{Hughes 1, 2}, but we can dodge the issue by using
the characterization: a map is an approximate stratified system of stratified fibrations if the stratified homotopy
link in the mapping cylinder is a genuine stratified system of stratified fibrations.  

\subsubhead6.10.3 Locally finite iterated homology\endsubsubhead
The same result holds for locally finite homology if the map is proper:
\proclaim{Proposition} Suppose $E\to X\to Y$ are as in 6.10.1, $X, Y$ are locally compact and $f$ is proper. Then
there is an equivalence
$$\Bbb H^{lf}(Y,\Bbb H(f;\Bbb E)))@>\simeq >> \Bbb H^{lf}(X;\Bbb E).$$
\endproclaim
Note the point inverses of $f$ are compact so we do not have to specify locally finite or compact supports for the
homology used in the cosheaf
$\Bbb H(f;\Bbb E)$. There is also a version for nonproper maps and locally finite homology in the cosheaf. This
requires a stratified proper homotopy local triviality hypothesis on the fibers of $f$ over strata in $Y$. 
 
\subsubhead 6.10.4 The map of spectra\endsubsubhead
We begin with the structure map of the spectrum over $X$:
$$S^k\wedge_XE_n@>>>E_{n+k}.$$
A space over $X$ with a section has a partial quotient that gives a space over $Y$ with a section: given $X\to F\to X$
define $F/{\simeq}$ by identifying each of the images $f^{-1}(y)\to F$ to a point. Taking $y$ to this point defines the
section $Y\to F/{\simeq}$. Applying this to the structure map gives spaces over $Y$:
$$S^k\wedge_Y(E_n/{\simeq})=(S^k\wedge_XE_n)/{\simeq} @>>> E_{n+k}/{\simeq}.$$
This map is adjoint to a map to the loop space over $Y$:
$$E_n/{\simeq}@>>> \Omega^k_Y(E_{n+k}/{\simeq}).$$
Take the (homotopy) direct limit $k\to\infty$ and the right side becomes homology, by definition. This gives a map of
spectra  over $Y$
$$E_n/{\simeq}@>>> \Bbb H(f;\Bbb E)_n,$$
To get the associated map of homology divide out $Y$ and apply $\Omega^n$:
$$\Omega^n(E_n/X)=\Omega^n(E_n/{\simeq}/Y)@>>> \Omega^n(\Bbb H(f;\Bbb E)_n/Y)$$
Taking the limit $n\to \infty$ gives $\Bbb H(X;\Bbb E)$ on the left and $\Bbb H(Y;\Bbb H(f;\Bbb E)_n/Y))$ on the right.
This is the map.

We see that the theorem amounts to rearranging the direct limits in taking homology. 

\subsubhead 6.10.5 The proof\endsubsubhead
First note that the classical (constant-coefficient) case follows from the Characterization Theorem. Suppose $\Bbb E$
is a spectrum.   Define a functor on maps $U\to Y$ by  taking the
$0^{th}$ space of the homology of $U$ and ignoring $Y$: $\Bbb H(U;\Bbb E)_0$. This satisfies the axioms. The conclusion
reconstructs the spectrum structure  and for $U\to Y$ homotopy stratified provides an assembly equivalence 
$$\Bbb H(Y;\Bbb H(p^{-1}(?);\Bbb E))@>\simeq>>\Bbb H(U; \Bbb E).$$ 

We now sketch the proof in the general case.
First, by comparing excision fibrations and using the cofibration condition on the filtration of $Y$, reduce to the case
where
$Y$ has one stratum and $X\to Y$ is a stratified fibration.

Next by using ANR and compact support properties reduce to the case $Y$ a finite complex. 

Again by comparing excision fibrations and working by induction on the number of cells we reduce to the relative case
$Y=(D^k,S^{k-1})$. 
Over a disk a stratified fibration is stratified homotopy equivalent to a product, so the disk case is a product
$\Bbb E\times D^n\to W\times D^n\to D^n$, where $\Bbb E\to W$ is a spectral cosheaf. 

Further excisions reduce this to the case $n=0$ where the equivalence is trivial.

\Refs
\noindent Anderson, Douglas R., Munkholm, Hans J{\o}rgen
\ref\no1\paper Geometric modules and Quinn homology theory\jour $K$-Theory\vol  7\yr 1993\pages 443--475\finalinfo MR
95f:55003
\endref

\noindent Aravinda, C. S., Farrell, F. T.; Roushon, S. K.
\ref\no1\paper Algebraic $K$-theory of pure braid groups\jour Asian J. Math.\vol 4 \yr 2000\pages 337--343\finalinfo MR
2002a:19002
\endref

\noindent Bartels, A.;  Reich, H.
\ref\no1\paper ,
On the Farrell-Jones Conjecture for higher algebraic
        K-theory\paperinfo preprint arXiv:math.AT/0308030\endref
        
\noindent Berkove, E.; Farrell, F. T.; Juan-Pineda, D.; Pearson, K.
\ref\no1\paper The Farrell-Jones isomorphism conjecture for finite covolume
            hyperbolic actions and the algebraic $K$-theory of Bianchi
            groups\jour Trans. Amer. Math. Soc.\vol 352\yr 2000\pages 5689-- 5702\finalinfo MR 2001b:57052\endref

\noindent Connolly, Frank, Ko{\'z}niewski, Tadeusz
\ref\no1\paper Rigidity and crystallographic groups. I\jour Invent. Math.\vol 99 \yr 1990\pages 25--48\finalinfo MR
91g:57019
\endref

\noindent Davis, James F.; L{\"u}ck, Wolfgang
\ref\no1\paper Spaces over a category and assembly maps in isomorphism
            conjectures in $K$- and $L$-theory\jour $K$-Theory\vol 15\yr 1998\pages 201-- 252\finalinfo MR 99m:55004\endref

\noindent Dress, Andreas W. M. 
\ref\no1\paper Contributions to the theory of induced representations
     \pages 183-- 240\jour Springer Lecture Notes in Math.\vol 342
 \yr 1973\finalinfo MR 52 5787\endref

\noindent Farrell, F. T., Jones, L. E.
\ref\no 1\paper Isomorphism conjectures in algebraic $K$-theory\jour J. Amer. Math. Soc.\vol 6\yr 1993\pages 249-- 297\finalinfo MR 93h:57032
\endref
\ref\no2\paper Stable pseudoisotopy spaces of compact non-positively curved
            manifolds\jour J. Differential Geom.\vol  34\yr 1991\pages 769--834\finalinfo MR 93b:57020\endref

\noindent Hughes, Bruce
\ref\no1\paper Stratified path spaces and fibrations\jour Proc. Roy. Soc. Edinburgh Sect. A\vol129
\yr 1999\pages 351--384\finalinfo{MR 2000d:57041}\endref
\ref\no2\paper Products and adjunctions of manifold stratified spaces\jour Topology Appl.\vol124
\yr 2002\pages 47--67\finalinfo{MR 1 926 134}\endref
\ref\no3\paper The approximate tubular neighborhood theorem \jour Ann. of Math.\vol 156\yr 2002\pages 867-- 889\finalinfo MR 2004b:57031\endref

\noindent Igusa, Kiyoshi
\ref\no1\book Higher Franz-Reidemeister torsion\jour AMS/IP Studies in Advanced Mathematics\vol 31\publ American Mathematical Society \yr 2002\pages xxii+370\endref

\noindent Pedersen, Erik Kj{\ae}r
\ref\no1\paper Controlled algebraic $K$-theory, a survey\inbook Geometry and topology: Aarhus (1998)\jour AMS Contemp.
Math.\vol 258 \yr2000
\pages 351--368\finalinfo MR 2001i:57034\endref 

\noindent Pedersen, Erik Kj{\ae}r; Yamasaki, Masayuki
\ref\no1\paper Stability in controlled L-theory\paperinfo preprint arXiv:math.GT/04022128\endref

\noindent  Quinn, Frank
\ref\no1\paper Ends of maps. I\jour Ann. of Math. (2)\vol 110 \yr 1979\pages 275--331\finalinfo MR 82k:57009\endref
\ref\no2\paper Ends of maps. II\jour Invent. Math.\vol 68 \yr 1982\pages 353--424\finalinfo MR 84j:57011\endref
\ref\no3\paper Ends of maps. III. Dimensions $4$ and $5$\jour J. Differential Geom.\vol  17\yr \pages 503--521\finalinfo
MR 84j:57012\endref
\ref\no4\paper Ends of maps. IV. Controlled pseudoisotopy\jour Amer. J. Math.\vol  108\yr \pages 1139--1161\finalinfo
MR 88j:57014\endref
\ref\no5\paper Homotopically stratified sets\jour J. Amer. Math. Soc.\vol 1 \yr1988 \pages 441--499\finalinfo {MR
89g:57050}\endref
\ref\no6\paper Geometric algebra\inbook Algebraic and geometric topology (New Brunswick, N.J., 1983)\jour Springer
Lecture Notes in Math.\vol  1126\yr
1985\pages182--198 \finalinfo MR 86m:57023\endref

\ref\no 7\paper Algebraic $K$-theory of poly-(finite or cyclic) groups\jour Bull. Amer. Math. Soc. (N.S.)\vol 12 \yr1985
\pages 221--226
\finalinfo {MR 86e:18015}\endref
\ref\no8\paper Hyperelementary induction for $K$-theory of virtually abelian groups\finalinfo arXiv:math.KT/0509294 preprint 2005\endref

\noindent Roushon, S. K.
\ref\no1\paper $K$-theory of virtually poly-surface groups\jour Algebr. Geom. Topol.\vol 3\yr 2003\pages 103-- 116\endref
\noindent Suslin, A. A.
\ref\no1\paper On the equivalence of $K$-theories\jour Comm. Algebra\vol 9\yr 1981\pages 1559-- 1566\finalinfo MR 82m:18008\endref

\endRefs
\bye